\documentclass[a4paper,11pt]{article}
\usepackage[normalem]{ulem}
\usepackage{cancel}
\usepackage{amsmath}
\usepackage{amssymb}
\usepackage{array}
\usepackage{ragged2e}
\usepackage{caption}
\usepackage{subcaption}
\usepackage{booktabs}
\usepackage[table,xcdraw]{xcolor}
\usepackage{graphicx}
\usepackage{titlesec}
\usepackage{rotating}
\usepackage{changepage}
\usepackage{lineno}
\linenumbers
\nolinenumbers
\usepackage{subcaption}

\usepackage{soul}
\usepackage{rotating}
\usepackage{adjustbox}
\usepackage{bm}

\usepackage{arydshln}
\usepackage{adjustbox, lipsum}
\usepackage{enumerate}
\usepackage[ruled,vlined,linesnumbered]{algorithm2e}
\usepackage{rotating}
\usepackage[utf8]{inputenc}
\usepackage{amsfonts}
\usepackage{pstricks}
\usepackage{pdflscape}
\usepackage{graphicx}
\usepackage{afterpage}
\usepackage{graphics}

\usepackage{hyperref}
\usepackage{rotating}
\usepackage{natbib}
\usepackage{booktabs}

\usepackage{authblk}
\usepackage{multirow}
\usepackage{mathtools}
\usepackage{scalerel}
\newtheorem{prop}{\bf Proposition}[section]

\SetKwRepeat{Do}{do}{while}
\providecommand{\norm}[1]{\lVert#1\rVert}
\DeclarePairedDelimiter\abs{\lvert}{\rvert}%

\newcommand{\dem}{\par \noindent{\bf Proof:} }
\newcommand{\fin}{\hfill $\square$  \par \bigskip}

\usepackage{xcolor}
\definecolor{gr}{RGB}{0,153,0}
\newcommand{\ncon}{NC-OWA-SVR}
\newcommand{\con}{C-OWA-SVR}
\usepackage{todonotes}
\makeatletter 
\@mparswitchfalse%
\makeatother
\normalmarginpar 
\definecolor{cyanazure}{rgb}{0, 0.7, 1}

\newcommand{\msout}[1]{\text{\sout{\ensuremath{#1}}}}
\definecolor{armygreen2}{rgb}{0.99, 0.53, 0.93}

\definecolor{armygreen}{rgb}{0.19, 0.53, 0.43}

\newcommand{\AM}[2]{{\color{armygreen2}#1 \ifmmode\msout{#2}\else\sout{#2}\fi}}

\title{Ordered Weighted Average Support Vector Regression}
\date{}

\author[(a,b)]{Luisa I. Martínez-Merino\footnote{ Corresponding author}}
\author[(b)]{Justo Puerto}
\author[(a)]{Antonio M. Rodríguez-Chía}

\affil[(a)]{\small{Departamento de Estadística e Investigación Operativa,

Universidad de Cádiz, Cádiz, Spain, 

luisa.martinez@uca.es, antonio.rodriguezchia@uca.es}}
\affil[(b)]{\small{Instituto de Matemáticas de la Universidad de Sevilla (IMUS), Sevilla, Spain, puerto@us.es}}

\allowdisplaybreaks
\setcitestyle{authoryear,open={(},close={)}}

\makeatletter
\let\oldabs\abs
\def\abs{\@ifstar{\oldabs}{\oldabs*}}
\makeatother
\begin{document}

%
\maketitle

\begin{abstract}

 This paper introduces a novel Support Vector Regression (SVR) model that incorporates Ordered Weighted Average (OWA) operators to differently penalize deviations of observations from the $\varepsilon$-strip. The penalty is determined based on the position of each deviation in the ordered vector of all deviations. A key contribution of this work is the development of two nonlinear formulations: a continuous formulation for non-decreasing monotone weight vectors and a mixed-integer formulation for general weight vectors. By leveraging dual approaches associated with these formulations, the model accommodates nonlinear kernel functions.
To enhance computational efficiency, two heuristic approaches are proposed for deriving suitable regression hyperplanes in reduced computation times, as compared to exact methods. Additionally, a third heuristic approach is developed to handle nonlinear kernel functions effectively. Computational experiments conducted on both real and synthetic datasets demonstrate that the exact formulations yield remarkable regression functions with respect to the following standard metrics: mean absolute errors (MAE) and mean squared errors (MSE). The heuristic approaches are shown to be particularly efficient for larger datasets, striking a balance between solution quality and computational effort.

\end{abstract}

\textbf{Keywords:} support vector regression,  Ordered Weighted Average operators, mixed integer quadratic programming,  goodness of fitting

\section{Introduction}%
\label{sec:Introduction}%

Support Vector Regression (SVR) is a Mathematical Programming based regression technique which has become popular in recent years, see \cite{Smola2004}. In fact, SVR has shown to be useful in different areas such as  medicine \citep{Ahadian2022}, economy \citep{Wang2023}, mechanics \citep{Li2023} or aeronautics \cite{Liu2021}, to name a few.
Since its introduction, see \cite{Vapnik1994}, many contributions have developed some extensions of the classic SVR which try to improve  its performance in different contexts. 

\subsection{Contributions}

In this paper, we propose a new SVR model including the ordered weighted average (OWA)  to penalize the
 deviations of the observations with respect to the $\varepsilon$-margin. 
The use of OWA operators has been tackled before in classic regression models. Particularly, in \cite{BLANCO2018}, a family of exact models to locate regression hyperplanes considering  different distance-based deviations and OWA operators is developed. Besides, we can find some works focusing on heuristic procedures to include OWA operators in regression, see \cite{FLORES2020} and \cite{FLORES2022}. However, as far as we are aware, SVR including OWA operators has not been  considered before. 

To develop this new model of SVR including OWA, we  have been inspired  by the idea proposed in \cite{Marin2022} where  OWA operators are included in the hinge loss function of the classic soft-margin Support Vector Machine (SVM). Specifically, our contributions in this paper are the following:

\begin{itemize}
 \item We introduce an exact methodology to include OWA operators in classic SVR. First, a general model is proposed where any non-negative  penalty weight vector is considered. Then, we focus on the particular case of monotone non-decreasing weight vectors.
 
 \item We also analyze how non-linear kernel functions can be included in the proposed formulations in order to obtain non-linear regression functions.
 
 \item We develop two different heuristic procedures to obtain accurate linear SVR including OWA operators. One of them is based on the closed coefficient formulas provided by dual classic SVR formulation.  The other heuristic is a two-step method which obtains a feasible solution of the considered problem taking into account the order of the deviations provided by the classic model.

\item  We introduce a heuristic procedure for SVR including OWA operators and non-linear kernel functions. This is a two-step method which considers the order of deviations given by the classic SVR including non-linear kernels and the dual formulation developed for the exact approach.

 \item We report the results of some computational experiments where we  assess the performance of the exact and heuristic approaches proposed along the paper. 
\end{itemize}

The remainder of the paper is organized as follows. First, Section \ref{review} presents some relevant works related to SVR. Section \ref{sec:models} is devoted to the description of the problem and the introduction of an exact methodology for the SVR model including OWA operators. In Section \ref{mon}, the convex case is analyzed.  
Section \ref{sec:kernel} addresses the analysis of the model to accommodate non-linear kernel functions. In Section \ref{sec:heur},  two heuristic approaches are developed to obtain regression hyperplanes which provide good quality performance in smaller computational times than the exact solution method. Besides, a heuristic approach allowing the use of non-linear kernel functions is developed. Section \ref{sec:comp} focuses on the analysis of the computational experiments carried out on several real-life datasets already used in related literature. Finally, in  Section \ref{sec:conclusions}, the main conclusions of our paper and some future research lines are detailed.

\subsection{ Literature review}\label{review}

Some SVR models have extended the idea of quantile regression to SVR. \cite{Takeuchi2006} introduces the first model considering quantile regression in SVR proposing a non-parametric quantile regression model. Then, \cite{Anand2020} presents a support vector quantile regression model which considers an asymmetric $\varepsilon$-insensitive pinball loss function for measuring the empirical risk. More contributions in this sense can be seen in \cite{Yang2019} or \cite{Ye2021}.

Other extensions of SVR focus their attention on providing SVR-based models which are computationally faster. Particularly, the twin support vector regression determines two non-parallel functions including an $\varepsilon$-insensitive margin, see \cite{Peng2010}. This involves solving two quadratic models smaller than the quadratic model of the classic SVR. Besides, in \cite{Ye2024}, a model combining quantile and twin SVR is developed to capture the characteristics of more complex data.

Classic SVR is sensitive to outliers and noise in datasets. This performance has tried to be partially mitigated using different proposals, not all of them directly comparable as in \cite{Hu2015}, \cite{Anand2020_2}, \cite{Fu2023} and \cite{Zhu2023}.

 However, to the best of our knowledge, Ordered Weighted Average (OWA) operators have not been introduced before in  SVR. The use of OWA operators allows to differently penalize the deviations according to their sizes, \cite{BLANCO2018,Marin2022}. This approach is more flexible than the original one in the sense that the decision maker can determine the  penalties of the deviations based on their sizes. For instance, assigning small weights to big deviations could be a way to limit the influence of outlier data.

\section{A general SVR model including OWA operators}
\label{sec:models}

Let $N=\{1,\ldots, n\}$ be a set of observations of a training sample. Each observation, $i\in N$, is represented by a vector $(\bm{x}_i,y_i)\in\mathbb{R}^d\times\mathbb{R}$, where $\bm{x}_i\in \mathbb{R}^d$ is the vector of feature values, $y_i\in\mathbb{R}$ is the value of the feature which must be estimated and let  $D=\{1,\ldots,d\}$. In regression models, the goal is to obtain a function $f:\mathbb{R}^d\rightarrow \mathbb{R}$ that fits the data, thus $y$-values can be estimated by using $f(\bm{x})$.

Particularly, classic Support Vector Regression  model seeks for the optimal hyperplane ($f(\bm{x})=\bm{w}^T\cdot \bm{x}+b$, with $\bm{w}\in \mathbb{R}^d$, $b\in \mathbb{R}$) that provides at most an $\varepsilon$-deviation from $y$-values, see \cite{Vapnik1994}. This classic model is given by:
\begin{alignat*}{3}
& \min_{\bm{w},b}&  \quad & \frac{1}{2}\norm{\bm{w}}^2\nonumber\\
&\mbox{s.t.}&&y_i-\bm{w}^T\bm{x}_i-b\leq \varepsilon,&\quad& i\in N,\\
&&&\bm{w}^T\bm{x}_i+b-y_i\leq \varepsilon,&& i\in N,\\
&&&\bm{w}\in\mathbb{R}^d,b\in\mathbb{R}.
\end{alignat*}

In this formulation, $\bm{w}$ are the variables related to the hyperplane coefficients and $b$ refers to the independent term. Besides, $\varepsilon$ is a parameter determined by the decision maker that establishes the tolerance against misfitting.

Since, depending on the data, it is not possible to obtain a hyperplane satisfying the above constraints, a soft margin version was also introduced inspired in similar insights used in classification for Support Vector Machine (SVM), see \cite{CortesVapnik1995}. Thus, the soft-margin Support Vector Regression model (SVR) can be stated as follows:

\begin{alignat}{3}
\mbox{(SVR)}& \min_{\bm{w},b,\bm{\xi},\bm{\xi}^*}&  \quad & \frac{1}{2}\norm{\bm{w}}^2+C\sum_{i=1}^n (\xi_i+\xi^*_i)\nonumber\\
&\mbox{s.t.}&&y_i-\bm{w}^T\bm{x}_i-b\leq \varepsilon+\xi_i,&&i\in N,\label{SVR1}\\
&&&\bm{w}^T\bm{x}_i+b-y_i\leq \varepsilon+\xi_i^*,&&i\in N,\label{SVR2}\\
&&&\xi_i,\xi_i^*\geq 0,&&i\in N,\label{SVR3}\\
&&&\bm{w}\in\mathbb{R}^d,b\in\mathbb{R}.\label{SVR4}
\end{alignat}

Note that, in this formulation, two new families of non-negative variables $\bm{\xi}$ and $\bm{\xi}^*$ are included. These new variables represent the deviations of the observations with respect to the established $\varepsilon$ margin. Now, the objective function is a balance between the minimization of the norm associated with $\bm{w}$ and the minimization of the  deviations associated with the observations. The trade-off between both objectives is obtained using a parameter $C>0$.

It must be emphasized that non-linear regression functions can also be obtained based on SVR model. In fact, the implicit mapping provided by kernel functions together with the dual formulation of SVR allow us to find adequate non-linear regression functions. For more details, see \cite{Smola2004}.

OWA operators have been successfully applied to the location of fitting hyperplanes in the context of multilinear regression, see \cite{BLANCO2018,BLANCO2021} where a solid theoretical framework is introduced. For Support Vector Regression, this idea has not been addressed in the literature yet. The goal of including OWA operators in SVR is to minimize a balance between the norm of hyperplane coefficients (the inverse of the margin between the two parallel hyperplanes) and the weighted sum of the sorted vector of deviations. Given a vector of weights $\bm{\lambda}=(\lambda_1,\ldots,\lambda_n)$ and denoting by 
$(\bm{\xi+\xi^*})_{()}=\Big((\xi+\xi^*)_{(1)},\ldots,(\xi+\xi^*)_{(n)}\Big)$
the vector of absolute deviations sorted in non-decreasing order, the mentioned model could be described as follows,
\begin{alignat*}{3}
& \min_{\bm{w},b,\bm{\xi},\bm{\xi}^*}&  \quad & \frac{1}{2}\norm{\bm{w}}^2+C\sum_{i=1}^n \lambda_i(\xi+\xi^*)_{(i)} \\ 
&\mbox{s.t.}&& \eqref{SVR1}-\eqref{SVR4}.
\end{alignat*}

Note that the main benefit with respect to the classic model is the different penalization of absolute deviations. This model including OWA-operators is more flexible since it allows the decision maker to establish a penalization (weight) depending on the relative position of the deviation of the observation with respect to the others. Observe that for the particular case of $\lambda_i=1$ for $i\in N$, the resulting model is the classic SVR and for $\lambda_p=1$ and $\lambda_i=0$ for $i\neq p$ results in the $p$-quantile SVR model.

We introduce a mixed integer quadratic formulation to include OWA in SVR for any $\bm{\lambda}$ such that $\lambda_i\geq 0$ for $i\in N$. For this aim, it is necessary to consider previous variables ($\bm{w},b,\bm{\xi},\bm{\xi}^*$) and two additional families of variables,
\begin{eqnarray*}
z_{ij}&=&\begin{cases}
1,&\mbox{if deviation of $i$ is in the $j$-th position of the vector}\\
&\mbox{of deviations sorted in non-decreasing order,}\\
0,&\mbox{otherwise.}
\end{cases}\mbox{for }i,j\in N,\\
\theta_k&=&
\mbox{deviation of the observation in position $k$ of the sorted}\\ 
&&\mbox{vector of deviations in non-decreasing order, for $k\in N$.}
\end{eqnarray*}

The obtained formulation is the following,
\begin{alignat}{3}
\mbox{(\ncon)}\quad & \min_{\bm{w},b,\bm{\xi},\bm{\xi}^*,\bm{z},\bm{\theta}}&  \quad & \frac{1}{2}\norm{\bm{w}}^2+C\sum_{k=1}^n\lambda_k\theta_k\nonumber\\
&\quad \quad \mbox{s.t.}&&\eqref{SVR1}-\eqref{SVR4},\nonumber\\
&&&\sum_{i=1}^n z_{ij}=1,&&j\in N,\label{rest_z1}\\
&&&\sum_{j=1}^n z_{ij}=1,&&i\in N,\label{rest_z2}\\
&&&\theta_k\geq \xi_i+\xi^*_i-M(1-\sum_{j=1}^k z_{ij}),&&\: \: i,k\in N,\label{nc_rest}\\
&&&z_{ij}\in\{0,1\},&&i,j\in N,\label{z_dom}\\
&&&\theta_k\geq 0,&&k\in N. \label{theta_dom}
\end{alignat}

Observe that the objective function of \ncon\xspace is a balance between the minimization of the norm associated with $\bm{w}$ and the ordered weighted average of the deviations. Constraints \eqref{SVR1} to \eqref{SVR4} are the classic constraints related to SVR which determine the deviation associated with each observation depending on the hyperplane. Constraints \eqref{rest_z1} and \eqref{rest_z2} ensure that each position in the sorted vector of deviations is associated with just one observation and viceversa. Besides, family of constraints \eqref{nc_rest} restricts each $\theta_k$-variable to take the value of deviation which is in position $k$, for $k\in K$, of the vector of deviations sorted in non-decreasing order. Note that $M>0$ appearing in constraints \eqref{nc_rest} is a large enough constant. Constraints \eqref{z_dom} and \eqref{theta_dom} are the domains of the variables for this formulation.

Note that the obtained formulation is denoted by \ncon \xspace since it considers non-convex OWA operators, i.e., the components of the $\lambda$-verctor are not given in increasing order. The resulting general formulation is a mixed integer quadratic formulation that includes a family of big-M constraints.

Given a fixed value of $\bm{z}$-variables vector denoted by $\hat{z}$, we consider below another quadratic continuous formulation defined as:

\begin{alignat}{3}
(\mbox{\ncon($\widehat{\bm{z}}$)})\quad& \min&  \quad & \frac{1}{2}\norm{\bm{w}}^2+C\sum_{k=1}^n \lambda_k 
\theta_k\nonumber\\
&\mbox{s.t.}&&\eqref{SVR1}-\eqref{SVR4},
\eqref{theta_dom},\nonumber\\
&&& 
\theta_k\geq \xi_i+\xi^*_i-M(1-\sum_{j=1}^k \widehat{z}_{ij}),&&\quad 
i,k\in N.
\end{alignat}

\begin{prop}
\label{Prop21}
The dual formulation of \ncon($\widehat{\bm{z}}$) is
\begin{alignat}{3}
(\mbox{\ncon$_D(\widehat{\bm{z}})$}) \:\:& \max_{\bm{\alpha},\bm{\alpha}^*, \bm{\mu}}&  \quad & -\frac{1}{2}\sum_{i=1}^n\sum_{j=1}^n(\alpha_i-\alpha_i^*)(\alpha_j-\alpha_j^*)\bm{x}_i \cdot \bm{x}_j\nonumber\\ 
&&&-\varepsilon\sum_{i=1}^n(\alpha_i+\alpha_i^*)+\sum_{i=1}^n y_i(\alpha_i-\alpha_i^*)\nonumber\\
&\quad \: \mbox{s.t.}&&\sum_{i=1}^n(\alpha_i-\alpha_i^*)=0,\label{dual_1}\\
&&&\alpha_i\leq C \sum_{k=1}^n \lambda_k \hat{z}_{ik},&& i\in N,\label{dual_2}\\
&&&\alpha_i^*\leq C \sum_{k=1}^n \lambda_k \hat{z}_{ik}&&i\in N, \label{dual_3}\\
&&&\alpha_i,\alpha_i^*\geq 0,&&i\in N, \label{dual_5} 
\\
&&&\mu_{ik}\geq 0,&&i,k\in N.\label{dual_6} 
\end{alignat}
\end{prop}

\dem
Formulation \ncon($\widehat{\bm{z}})$ can be rewritten in an equivalent way as:
\begin{alignat}{3}
\quad& \min&  \quad & \frac{1}{2}\norm{\bm{w}}^2+C\sum_{k=1}^n \lambda_k 
\sum_{i=1}^n (\xi_i+\xi_i^*)\hat{z}_{ik}\nonumber\\
&\mbox{s.t.}&&\eqref{SVR1}-\eqref{SVR4},\nonumber
\end{alignat}
The Lagrangian function related to this problem is
\begin{eqnarray*}
L(\bm{w},b,\bm{\xi},\bm{\xi}^*)&=&\frac{1}{2}\norm{\bm{w}}^2+C\sum_{k=1}^n \lambda_k
\sum_{i=1}^n (\xi_i+\xi_i^*)\hat{z}_{ik}+\sum_{i=1}^n\alpha_i(y_i-\bm{w}^T\bm{x}_i-b-\varepsilon-\xi_i)\\
&&+\sum_{i=1}^n\alpha_i^*(\bm{w}^T\bm{x}_i+b-y_i-\varepsilon-\xi_i^*)\\
&&-\sum_{i=1}^n\eta_i\xi_i-\sum_{i=1}^n\eta_i^*\xi_i^* ,
\end{eqnarray*}
where $\bm{\alpha}\geq \bm{0}$, $\bm{\alpha}^*\geq \bm{0}$, 
$\bm{\eta}\geq\bm{0}$ and $\bm{\gamma} \geq \bm{0}$ are the lagrangian multipliers. Using the necessary and sufficient optimality conditions, next statements must hold:

\begin{eqnarray}
\frac{\partial L(w,b,\xi,\xi^*)}{\partial w_{\ell}}&=&w_{\ell}-\sum_{i=1}^n \alpha_ix_{i\ell}+\sum_{i=1}^n\alpha_i^*x_{i\ell}=0,\, \ell\in D,\hspace*{0.5cm}\label{lag1}\\
\frac{\partial L(w,b,\xi,\xi^*)}{\partial b}&=&-\sum_{i=1}^n\alpha_i+\sum_{i=1}^n\alpha_i^*=0,\label{lag2}\\
\frac{\partial L(w,b,\xi,\xi^*)}{\partial \xi_i}&=&-\alpha_i+C\sum_{k=1}^{n} \lambda_k \hat{z}_{ik}-\eta_i=0,\: \,i\in N,\label{lag3}\\
\frac{\partial L(w,b,\xi,\xi^*)}{\partial \xi_i^*}&=&-\alpha_i^*+C\sum_{k=1}^{n} \lambda_k \hat{z}_{ik}-\eta_i^*=0,\: \,i\in N,\label{lag4}\\
\alpha_i(y_i-w^Tx_i-b-\varepsilon-\xi_i)&=&0,\: \,i\in N,\label{lag6} \nonumber\\
\alpha_i^*(w^Tx_i+b-y_i-\varepsilon-\xi_i^*&=&0,\: \,i\in N, \label{lag7} \nonumber\\
\eta_i\xi_i&=&0, \: \,i\in N,\label{lag9} \nonumber\\
\eta_i^*\xi_i^*&=&0,\: \,i\in N.\label{lag10} \nonumber
\end{eqnarray}
Particularly, from \eqref{lag1}, it can be derived that $w_{\ell}=\sum_{i\in N}(\alpha_i-\alpha_i^*)x_{i\ell}$ for any $\ell \in D$. This expression of $\bm{w}$ together with the use of complementary slackness conditions in the Lagrangian function results in
\begin{eqnarray*}
L(\bm{w},b,\bm{\xi},\bm{\xi}^*)&=&
-\frac{1}{2}\sum_{i=1}^n\sum_{j=1}^n(\alpha_i-\alpha_i^*)(\alpha_j-\alpha_j^*)\bm{x}_i \cdot \bm{x}_j
-\varepsilon\sum_{i=1}^n(\alpha_i+\alpha_i^*)\\
&&+\sum_{i=1}^n y_i(\alpha_i-\alpha_i^*).
\end{eqnarray*}

The above expression constitutes the objective function of \mbox{\ncon$_D(\widehat{\bm{z}})$}. Besides note that constraint \eqref{dual_1}-\eqref{dual_3} are obtained from \eqref{lag2}-\eqref{lag4} since $\eta_i, \eta_i^* \ge 0$ for any $i\in N$.

\fin

Using the properties of \ncon($\widehat{\bm{z}}$) and \ncon$_D(\widehat{\bm{z}})$ formulations, an alternative general formulation for SVR with OWA is derived. This formulation is given by

\begin{alignat}{3}
(\mbox{\ncon$_D$})& \min_{\bm{\alpha},\bm{\alpha}^*,\bm{\xi},\bm{\xi}^*}&  \quad & \frac{1}{2}\sum_{i=1}^n\sum_{j=1}^n(\alpha_i-\alpha_i^*)(\alpha_j-\alpha_j^*)\bm{x}_i \cdot \bm{x}_j+C\sum_{k=1}^n\lambda_k\theta_k\nonumber\\
&\mbox{s.t.}&&\eqref{SVR3},\eqref{rest_z1}-\eqref{theta_dom},\eqref{dual_1},\nonumber\\
&&&y_i-\sum_{j=1}^n (\alpha_j-\alpha_j^*)\bm{x}_i\cdot \bm{x}_j-b\leq \varepsilon+\xi_i,&&\hspace{-1.5cm}i\in N, \label{xi1}\\
&&&\sum_{j=1}^n (\alpha_j-\alpha_j^*)\bm{x}_i\cdot \bm{x}_j+b-y_i\leq \varepsilon+\xi_i^*,&&\hspace{-1.5cm}i\in N, \label{xi2}\\
&&&0\le\alpha_i \le C \sum_{k=1}^n\lambda_k z_{ik},&&\hspace{-1.5cm}i\in N,\label{alpha_dom}\\
&&&0\le\alpha^*_i \le C \sum_{k=1}^n\lambda_k z_{ik} ,&&\hspace{-1.5cm}i\in N.\label{alpha2_dom}
\end{alignat}

The resulting formulation depends on scalar products of $\bm{x}$-vectors. This is useful to accommodate non-linear kernel functions in the formulation. More details about this aspect are specified in Subsection \eqref{sec:kernel}.

\begin{prop}
\ncon\xspace and \ncon$_D$ are equivalent in the sense that optimal
solutions to \ncon\xspace can be built from optimal solutions of \ncon$_D$.
\end{prop}
\dem
We first prove that given  an optimal solution of NC-OWA-SVR,
$(\widehat{\bm{w}},\widehat{b},\widehat{\bm{\xi}},\widehat{\bm{\xi}}^*,\widehat{\bm{\theta}},\bm{\widehat{z}} )$, it can be built a feasible solution of NC-OWA-SVR$_D$, $(\widehat{\bm{\alpha}},\widehat{\bm{\alpha}}^*,\widehat{b},\widehat{\bm{\xi}},\widehat{\bm{\xi}}^*,\widehat{\bm{\theta}},\bm{\widehat{z}})$, with the same objective value.

Given a solution of NC-OWA-SVR, $(\widehat{\bm{w}},\widehat{\bm{b}},\widehat{\bm{\xi}},\widehat{\bm{\xi}}^*,\widehat{\bm{\theta}},\bm{\widehat{z}})$, then
$(\widehat{\bm{w}},\widehat{b},\widehat{\bm{\xi}},\widehat{\bm{\xi}}^*,\widehat{\bm{\theta}})$ is an optimal solution of $\mbox{\ncon}(\widehat{\bm{z}})$. 
By Proposition \ref{Prop21}, we have that exists $(\widehat{\bm{\alpha}},\widehat{\bm{\alpha}}^*)$ 
satisfying \eqref{dual_1}-\eqref{dual_6} with the same optimal value that
$\mbox{\ncon}(\widehat{\bm{z}})$. 

Hence, using \eqref{lag1}   we obtain that:

\begin{enumerate}
\item
$\displaystyle \frac{1}{2}\norm{\widehat{\bm{w}}}^2+C\sum_{k=1}^n\lambda_k\widehat{\theta}_k
=\frac{1}{2}\sum_{i=1}^n\sum_{j=1}^n(\widehat{\alpha}_i-\widehat{\alpha}_i^*)(\widehat{\alpha}_j-\widehat{\alpha}_j^*)x_i \cdot x_j+C\sum_{k=1}^n\lambda_k\widehat{\theta}_k.$
\item  Constraints \eqref{xi1} and \eqref{xi2} are derived from \eqref{SVR1} and \eqref{SVR2}, respectively.
\end{enumerate}

Moreover, using \eqref{lag3} and \eqref{lag4} we obtain \eqref{alpha_dom} and \eqref{alpha2_dom}, respectively. The remaining constraints of NC-OWA-SVR$_D$ are trivially fulfilled.
Therefore, we have that $(\widehat{\bm{\alpha}},\widehat{\bm{\alpha}}^*,\widehat{b},\widehat{\bm{\xi}},\widehat{\bm{\xi}}^*,\widehat{\bm{\theta}},\bm{\widehat{z}})$ is a feasible solution of NC-OWA-SVR$_D$ with the same objective value  as $(\widehat{\bm{w}},\widehat{\bm{b}},\widehat{\bm{\xi}},\widehat{\bm{\xi}}^*,\widehat{\bm{\theta}},\bm{\widehat{z}})$ for NC-OWA-SVR.

Conversely, given an optimal solution of \ncon$_D$,	$(\widehat{\bm{\alpha}},\widehat{\bm{\alpha}}^*,\widehat{b},\widehat{\bm{\xi}},\widehat{\bm{\xi}}^*,\widehat{\bm{\theta}},\widehat{\bm{z}})$, a feasible solution of \ncon\xspace with the
	same objective value, $(\widehat{\bm{w}},\widehat{b},\widehat{\bm{\xi}},\widehat{\bm{\xi}}^*,\widehat{\bm{\theta}},\widehat{\bm{z}})$, can be built just by considering  $\widehat{w}_{\ell}=\sum_{i=1}^n(\widehat{\alpha}_i-\widehat{\alpha}_i^*)x_{i\ell}$ for any 
    $\ell \in D$. 
    \mbox{\hspace*{0.5cm}} \fin

Observe  that \ncon$_D$ depends on the scalar product of $\bm{x}$-vectors which allows the use of kernel functions to obtain non-linear regression models as we will see in Section \ref{sec:kernel}. 
Both, \ncon\xspace and \ncon$_D$, can be applied to any given $\bm{\lambda}$ weight such that $\lambda_i\geq 0$ for $i\in N$. However, when considering non-decreasing $\bm{\lambda}$-weights, a simpler quadratic continuous formulation can be derived as we will analyze in the next subsection.


\section{OWA-SVR considering convex OWA operators}
\label{mon}

Let $\bm{\lambda}$ be a vector of non-decreasing weights, i.e, $\bm{\lambda}=(\lambda_1,\ldots,\lambda_n)$ where $0 \le \lambda_1\leq\lambda_2\leq\ldots\leq\lambda_n$. For a fixed vector of deviations $\overline{\bm{\xi}+\bm{\xi}^*}=(\overline{\xi}_1+\overline{\xi}^*_1,\ldots,\overline{\xi}_n+\overline{\xi}^*_n)$, the ordered weighted average of these deviations using $\bm{\lambda}$-weight can be modelled as follows:
\begin{alignat}{3}
\sum_{i=1}^n \lambda_i(\overline{\xi}+\overline{\xi}^*)_{(i)}=& \max_{\bm{z}}&  \quad & \sum_{i=1}^n\sum_{j=1}^n\lambda_j(\xi_i+\xi^*_i)z_{ij}\nonumber\\
&\mbox{s.t.}&&\eqref{rest_z1},\eqref{rest_z2},&& \nonumber \\
&&&z_{ij}\geq 0,&&i,j\in N.\label{z_relax}
\end{alignat}

Recall that $\bm{z}$-variables are the allocation of each individual with its corresponding position in the vector of deviations sorted in non-decreasing order. Note that integrality constraints for $\bm{z}$ variables can be relaxed due to the total unimodularity property of constraints \eqref{rest_z1},\eqref{rest_z2}. If we embed this model of the OWA deviations in the classic SVR model, it results in the following,
\begin{alignat}{3}
& \min_{\bm{w},b,\bm{\xi},\bm{\xi}^*}&  \quad & \frac{1}{2}\norm{\bm{w}}^2+\max_{z}\sum_{i=1}^n\sum_{j=1}^n C\lambda_j(\xi_i+\xi_i^*)z_{ij}\nonumber\\
&\mbox{s.t.}&&\eqref{SVR1}-\eqref{rest_z2},\eqref{z_relax}.\nonumber
\end{alignat}
Note that \eqref{SVR1} to \eqref{SVR4} are the original constraints of SVR to determine the deviations of the individuals with respect to the $\varepsilon$-strip. Constraints \eqref{rest_z1} and \eqref{rest_z2} are related to the order of the deviation vector. Using the dual form of the inner model, the resulting formulation is

\begin{alignat}{3}
\mbox{(\con)}& \min_{\bm{w},b,\bm{\xi},\bm{\xi}^*,\bm{u},\bm{v}}&  \quad & \frac{1}{2}\norm{w}^2+\sum_{i=1}^n u_i+\sum_{j=1}^n v_j\nonumber\\
&\:\:\mbox{s.t.}&&\eqref{SVR1}-\eqref{SVR4},\nonumber\\
&&&u_i+v_j\geq C\lambda_j(\xi_i+\xi_i^*),&&i,j\in N,\nonumber\\
&&&u_i,v_j\in\mathbb{R},&&i,j\in N.\nonumber
\end{alignat}

Observe that $\bm{v}$- and $\bm{u}$-variables are dual variables associated with constraints \eqref{rest_z1} and \eqref{rest_z2}.The resulting model is a convex quadratic continuous model.

\begin{prop}
The dual formulation of \con\xspace is given by

\begin{alignat}{3}
(\mbox{\con}_D) \: \:& \max_{\bm{\alpha},\bm{\alpha}^*,\bm{\beta}}&  & -\frac{1}{2}\sum_{i=1}^n\sum_{j=1}^n(\alpha_i-\alpha_i^*)(\alpha_j-\alpha_j^*)\bm{x}_i\cdot \bm{x}_j\nonumber\\ 
&&&-\varepsilon\sum_{i=1}^n(\alpha_i+\alpha_i^*)+\sum_{i=1}^n y_i(\alpha_i-\alpha_i^*)\nonumber\\
&\mbox{s.t.}&&\eqref{dual_1},\nonumber\\
&&&\alpha_i\leq C\sum_{j=1}^n\beta_{ij}\lambda_j,&&i\in N,\label{cdual2}\\
&&&\alpha_i^*\leq C\sum_{j=1}^n\beta_{ij}\lambda_j,&&i\in N,\label{cdual3}\\
&&&\sum_{j=1}^n\beta_{ij}=1,&&i\in N,\label{cdual4}\\
&&&\sum_{i=1}^n\beta_{ij}=1,&&j\in N,\label{cdual5}\\
&&&\alpha_i,\alpha_i^*\geq 0,&&i\in N,\label{cdual6}\\
&&&\beta_{ij}\geq 0,&&i,j\in N.\label{cdual7}
\end{alignat}
\end{prop}

\dem
The Lagrangian function related to \con\xspace can be expressed as:
\begin{eqnarray*}
L(\bm{w},b,\bm{\xi},\bm{\xi}^*,\bm{u},\bm{v})&=&\frac{1}{2}\norm{w}^2+\sum_{i=1}^n u_i+\sum_{j=1}^nv_j+\sum_{i=1}^n\sum_{j=1}^n \beta_{ij}(C\lambda_j(\xi_i+\xi_i^*)-u_i-v_j)\\
&&+\sum_{i=1}^n\alpha_i(y_i-w^Tx_i-b-\varepsilon-\xi_i)+\sum_{i=1}^n\alpha_i^*(w^Tx+b-y_i-\varepsilon-\xi_i^*)\\
&&-\sum_{i=1}^n\eta_i\xi_i-\sum_{i=1}^n\eta_i^*\xi_i^*,
\end{eqnarray*}
where $\bm{\beta}\geq \bm{0}$, $\bm{\alpha}\geq \bm{0}$, $\bm{\alpha}\geq \bm{0}$, $\bm{\eta}\geq \bm{0}$ are the lagrangian multipliers. Using the necessary and sufficient optimality conditions, it must hold:

\begin{eqnarray}
	\frac{\partial L(w,b,\xi,\xi^*,u,v)}{\partial w_{\ell}}&=& w_{\ell}-\sum_{i=1}^n\alpha_ix_{i\ell}+\sum_{i=1}^n\alpha_i^*x_{i\ell}=0,\,\ell \in D,
    \label{clag1}\\
		\frac{\partial L(w,b,\xi,\xi^*,u,v)}{\partial b}&=&-\sum_{i=1}^n\alpha_i+\sum_{i=1}^n\alpha_i^*=0,\label{clag2}\\
		\frac{\partial L(w,b,\xi,\xi^*,u,v)}{\partial \xi_i}&=&C\sum_{j=1}^n\beta_{ij}\lambda_j-\alpha_i-\eta_i=0,\,i\in N,\label{clag3}\\
		\frac{\partial L(w,b,\xi,\xi^*,u,v)}{\partial \xi_i^*}&=&C \sum_{j=1}^n\beta_{ij}\lambda_j-\alpha_i^*-\eta_i^*=0,\,i\in N,\label{clag35}\\
		\frac{\partial L(w,b,\xi,\xi^*,u,v)}{\partial u_i}&=&1-\sum_{j=1}^n\beta_{ij}=0,\,i\in N,\label{clag4}\\
		\frac{\partial L(w,b,\xi,\xi^*,u,v)}{\partial v_j}&=&1-\sum_{i=1}^n\beta_{ij}=0,\,i\in N,\label{clag5}\\	
		\beta_{ij}(C\lambda_j(\xi_i+\xi_i^*)-u_i-v_j)&=&0,\,i,j\in N,\label{clag6} \nonumber\\
		\alpha_i(y_i-w^Tx_i-b-\varepsilon-\xi_i)&=&0,\,i\in N,\label{clag7} \nonumber\\
		\alpha_i^*(w^Tx+b-y_i-\varepsilon-\xi_i^*)&=&0,\,i\in N,\label{clag8} \nonumber\\
		\eta_i\xi_i&=&0,\,i\in N,\label{clag9}\\
		\eta_i^*\xi_i^*&=&0,\,i\in N.\label{clag10} \nonumber 
\end{eqnarray}
From \eqref{clag1}, $\bm{w}$ can be expressed as $w_{\ell}=\sum_{i\in N}(\alpha_i-\alpha_i^*)x_{i\ell},\ell \in D$. Besides, in this case, constraint \eqref{dual_1} is obtained from \eqref{clag2}. Expressing \eqref{clag3} and \eqref{clag35} as inequalities provide constraints \eqref{cdual2} and \eqref{cdual3}. Finally, constraints \eqref{cdual4} and \eqref{cdual5} are inferred from \eqref{clag4} and \eqref{clag5}. Using the derived expression of $\bm{w}$ and the complementary conditions, the objective function of \con\xspace is obtained. The optimal solution, $(\widehat{\bm{\alpha}},\widehat{\bm{\alpha}}^*,\widehat{\bm{\beta}})$, of \con$_D$ provides the optimal values of $\bm{w}$ by calculating $\sum_{i\in N}(\widehat{\alpha}_i-\widehat{\alpha}_i^*)x_{i\ell},\ell\in D$. 
\fin
The following result describes a way to calculate the independent term of the hyperplane equation.

\begin{prop}\label{prop:b}
Let $(\widehat{\bm{\alpha}},\widehat{\bm{\alpha}}^*,\widehat{\bm{\beta}})$ be a solution of \con$_D$. Then, the independent term, $b$, of the hyperplane equation can be calculated as:
\begin{eqnarray}
\label{label_b}
\max_{i=1,\ldots,n}\{-\epsilon+y_i-\bm{w}^t\bm{x}_i:\alpha_i<C \sum_{j=1}^n\widehat{\beta}_{ij}\mbox{ or } \alpha_i^*>0\}&\le& b \nonumber\\\
\le 
\min_{i=1,\ldots,n} \{-\epsilon+y_i-\bm{w}^t\bm{x}_i: \alpha_i >0 \mbox{ or } \alpha_i^* < C\sum_{j=1}^n\widehat{\beta}_{ij}\}.&&
\end{eqnarray}

Moreover,if exists $i\in I$ such that 
$0<\widehat{\alpha}_i<C\sum_{j=1}^n \widehat{\beta}_{ij}\lambda_j$ or
$0<\widehat{\alpha}_i^*<C\sum_{j=1}^n \widehat{\beta}_{ij}\lambda_j$, then
$$b=y_i-\sum_{j=1}^n(\widehat{\alpha}_j-\widehat{\alpha}^*_j)x_i\cdot x_j-\epsilon.$$
\end{prop}
\dem
Expression \eqref{label_b} is straightforward derived form \eqref{SVR1} and \eqref{SVR2}. Moreover,  if it exists $i\in I$ such that $0<\widehat{\alpha}_i<C\sum_{j=1}^n \widehat{\beta}_{ij}\lambda_j$, then
by \eqref{clag3}, $\widehat{\eta}_i=C\sum_{j=1}^n \widehat{\beta}_{ij}  \lambda_j-\widehat{\alpha}_i >0$ and  using the complementary condition \eqref{clag9}, $\widehat{\xi}_i=0$. Then, the independent term, $b$, of the hyperplane equation can be calculated as:
$$b=y_i-\sum_{j=1}^n(\widehat{\alpha}_j-\widehat{\alpha}^*_j)x_i\cdot x_j-\epsilon.$$
The same value can be obtained if  it exists $i\in I$ such that $0<\widehat{\alpha}_i^*<C\sum_{j=1}^n \widehat{\beta}_{ij}\lambda_j$. 
\hfill \mbox{} \fin

Note that \con$_D$ depends on the scalar product of $\bm{x}$-vectors which allows one the use of kernel functions to obtain non-linear regression models as we will see in Subsection \ref{sec:kernel}. 

\section{Introducing kernel functions in OWA-SVR} 
\label{sec:kernel} 
In classic SVR, see \cite{Smola2004}, implicit mapping via kernel functions are used in the dual formulation to obtain non-linear regression functions. This idea is commonly known as kernel trick. Regarding OWA-SVR using general weights, \ncon$_D$ can accommodate kernel functions. Given a kernel function, $K(\bm{x}_i,\bm{x}_j)$ and a $\bm{\lambda}$-weight with $\lambda_i\geq 0,i\in N$, the formulation for OWA-SVR considering this kernel is,
\begin{alignat}{3}
(\mbox{\ncon$_K$})& \hspace{-0.1cm}\min_{\bm{\alpha},\bm{\alpha}^*,b,\bm{z},\bm{\xi},\bm{\xi}^*}&  \quad & \frac{1}{2}\sum_{i=1}^n\sum_{j=1}^n(\alpha_i-\alpha_i^*)(\alpha_j-\alpha_j^*)K(\bm{x}_i,\bm{x}_j)+C\sum_{k=1}^n\lambda_k\theta_k\nonumber\\
&\mbox{s.t.}&&\eqref{SVR3},\eqref{rest_z1}-\eqref{theta_dom},\eqref{dual_1},\eqref{alpha_dom},\eqref{alpha2_dom}\nonumber\\
&&&\hspace{-0.2cm}y_i-\sum_{j=1}^n (\alpha_j-\alpha_j^*)K(\bm{x}_i,\bm{x}_j)-b\leq \varepsilon+\xi_i,&&\hspace{-1.7cm}i\in N,\label{err1}\\
&&&\hspace{-0.2cm}\sum_{j=1}^n (\alpha_j-\alpha_j^*)K(\bm{x}_i,\bm{x}_j)+b-y_i\leq \varepsilon+\xi_i^*,&&\hspace{-1.7cm}i\in N.\label{err2}
\end{alignat}

Let $\bm{x}$ be a feature vector associated with certain individual. Then, the estimation of its $y$-value is calculated by 
$$y=\sum_{i \in N}(\widehat{\alpha}_i-\widehat{\alpha}_i^*)K(\bm{x}_i,\bm{x})+\widehat{b},$$
where $\widehat{\bm{\alpha}}$, $\widehat{\bm{\alpha}}^*$  and $\widehat{b}$ are the optimal values associated with $\bm{\alpha}$-, $\bm{\alpha}^*$  and $b$-variables. 

For the particular case of considering non-decreasing $\bm{\lambda}$ weights, the following model can be derived by using \con$_D$ formulation:

\begin{alignat}{3}
(\mbox{\con}_K)& \max_{\bm{\alpha},\bm{\alpha}^*,\bm{\beta}}&  \quad & -\frac{1}{2}\sum_{i=1}^n\sum_{j=1}^n(\alpha_i-\alpha_i^*)(\alpha_j-\alpha_j^*)K(\bm{x}_i,\bm{x}_j)\nonumber\\ 
&&&-\varepsilon\sum_{i=1}^n(\alpha_i+\alpha_i^*)+\sum_{i=1}^n y_i(\alpha_i-\alpha_i^*)\nonumber\\
&\mbox{s.t.}&&\eqref{dual_1}, \eqref{cdual2}-\eqref{cdual7}.\nonumber
\end{alignat}

In contrast with \ncon$_K$, no integrality constraints are needed in \con$_K$, i.e., the resulting model is a quadratic continuous formulation. In the next section, we present some heuristic procedures to provide good quality feasible solutions for OWA-SVR.

\section{Heuristic approaches for OWA-SVR}
\label{sec:heur}
In the literature, we can find different heuristic approaches which try to include OWA operators in classification and regression models, see for instance \cite{FLORES2020,FLORES2022,Maldonado2018,Maldonado2020}. In this paper, we introduce two different heuristic approaches to consider OWA operators in SVR.  Besides, we present a third heuristic procedure which is able to deal with non-linear kernel functions providing good quality non-linear regression functions.

The first heuristic approach, denoted as HEUR1, is based on the same rationale of a heuristic procedure presented in \cite{FLORES2022} to take into account OWA weight vectors in classic multilinear regression. As we will see in Subsection \ref{sec:heur1}, HEUR1 introduces OWA operators in the expressions associated with the classic coefficients of SVR hyperplanes. 
Obviously, this method is heuristic,
since no loss function minimization is considered in this method. In spite of that, as we will also observe in the experiments, some interesting results can be obtained with our first heuristic.

The second heuristic approach, HEUR2, consists in a two-step procedure. In the first step, classic SVR model is applied and the individuals are sorted considering the order given by the SVR solution. In the second step, it is solved a SVR model where the deviation of each individual is penalized according to the first step. This kind of heuristic approach has been developed before for SVM in the context of classification, see \cite{Maldonado2018}. Once again, this second proposed method gives remarkable good results.

The third heuristic procedure, HEUR$_K$, tries to approach the solutions of OWA-SVR with a certain kernel function, $K(\cdot,\cdot)$. To this aim, a two-step algorithm is proposed. First, classic SVR considering $K(\cdot,\cdot)$ is applied. Then, the obtained deviations of the individuals are sorted in non-decreasing order.  This order determines the penalizations of the individuals in the second step, where a formulation derived from  \ncon$_K$ is solved.

Next subsections detail the  proposed heuristic approaches: HEUR1, HEUR2 and HEUR$_K$.  Besides, the algorithms of these heuristic approaches are presented by flowcharts in Appendix \ref{App:alg}.

\subsection{Heuristic method based on closed OWA formulas (HEUR1)}
\label{sec:heur1}

In this heuristic, we apply the closed form expansion of $\bm{w}$ and $b$ in standard SVR replacing standard moments by OWA moments.
Recall that in classic SVR
$$w=\sum_{i\in N} (\alpha_i-\alpha_i^*)x_i$$ 
where $\bm{\alpha}$ and $\bm{\alpha}^*$ are the dual multipliers of the constraints \eqref{SVR1} and \eqref{SVR2}. Besides, analogously to Proposition \ref{prop:b}, $b$ can be computed using the complementary slackness conditions as:
$$
\max\{-\epsilon+y_i-\bm{w}^t\bm{x}_i:\alpha_i<C \mbox{ or } \alpha_i^*>0\} \le b \le \min \{-\epsilon+y_i-\bm{w}^t\bm{x}_i: \alpha_i >0 \mbox{ or } \alpha_i^* < C\}
$$    
If  $\alpha_i\in (0,C)$ or $\alpha_i^*\in (0,C)$ for some $i=1,\ldots,n$ the inequalities become equalities and determine $b$, see \cite{Smola2004} for more details.

Let $\widehat{\bm{\alpha}}$ and $\widehat{\bm{\alpha}}^*$ be the solution values of $\bm{\alpha}$- and $\bm{\alpha}^*$-variables in the classic dual formulation of SVR. Besides, consider $$\bm{(\widehat{\alpha}-\widehat{\alpha}^*)x_j}=((\widehat{\alpha_1}-\widehat{\alpha}_1^*)x_{1j},(\widehat{\alpha}_2-\widehat{\alpha_2}^*)x_{2j},\ldots,(\widehat{\alpha}_n-\widehat{\alpha}_n^*)x_{nj}),\mbox{ for }j\in D,$$
and $\bm{((\widehat{\alpha}-\widehat{\alpha}^*)x_j)}_{()}$ is the same vector but sorting its values in non-decreasing order.
The resulting regression coefficients, $\bm{w}^{OWA}$, can be computed as:
\begin{eqnarray*}
w^{OWA}_j&=&\sum_{i\in N}\lambda_i\bm{((\widehat{\alpha}-\widehat{\alpha}^*)x_j)}_{(i)}
\end{eqnarray*}

Furthermore, $b^{OWA}$ values can also be computed following the same idea. A detailed description of this heuristic procedure is detailed in Algorithm \ref{alg:heur1}.

\begin{algorithm}[htbp]
{\small
	\KwData{Training sample composed by a set of $n$ individuals with $d$ features.}
	\KwResult{Hyperplane whose coefficients introduce OWA operators.}
 Solve the dual form of classic SVR.\\
 Let $\widehat{\bm{\alpha}}$ and $\widehat{\bm{\alpha}}^*$  be optimal values of the dual model.\\
 Compute  $\bm{(\widehat{\alpha}-\widehat{\alpha}^*)x_{\ell}}=((\widehat{\alpha}_1-\widehat{\alpha}_1^*)x_{1\ell},(\widehat{\alpha}_2-\widehat{\alpha}_2^*)x_{2\ell},\ldots,(\widehat{\alpha}_n-\widehat{\alpha}_n^*)x_{n\ell}),\mbox{ for }\ell\in D.$\\
 Compute  $\bm{((\widehat{\alpha}-\widehat{\alpha}^*)x_\ell)}_{()}$ sorting $\bm{(\widehat{\alpha}-\widehat{\alpha}^*)x_j}$ vector in non-decreasing order, such that, $((\widehat{\alpha}-\widehat{\alpha}^*)x_\ell)_{(i)}= (\widehat{\alpha}-\widehat{\alpha}^*)x_{v(i),\ell}$.\\
 Define: $w^{OWA}_j=\sum_{i\in N}\lambda_i\bm{((\widehat{\alpha}-\widehat{\alpha}^*)x_j)}_{(i)}$ for $j\in D$.\\
 Define $k=0$ \\
\For{ $i\in N$}{
\If{$0<\widehat{\alpha}_i<C$}{
$k=1$\\
$b^{OWA}=y_i-\bm{w}^{OWA}\bm{x}_i-\epsilon.$\\
break
}
\If{$0<\widehat{\alpha}_i^*<C$}{
$k=1$\\
$b^{OWA}=y_i-\bm{w}^{OWA}\bm{x}_i+\epsilon.$\\
break
}
}
\If{$k=0$}{
Let $i_0\in N$ such that $i_0=argmax\{-\epsilon+y_i-\bm{w}^t\bm{x}_i:\alpha_i<C \lambda_{v(i)}\mbox{ or } \alpha_i^*>0\} $, where $\bm{w}$ is optimal for classic SVR.\\
Compute $b_1=-\epsilon+y_{i_0}-\bm{w}^{OWA}\bm{x}_{i_0}$.\\
Let $i_0^{\prime}\in N$ such that $i_0^{\prime}=argmin \{-\epsilon+y_i-\bm{w}^t\bm{x}_i: \alpha_i >0 \mbox{ or } \alpha_i^* < C \lambda_{v(i)}\}$, where $\bm{w}$ is optimal for classic SVR.\\
Compute $b_2=-\epsilon+y_{i_0}-\bm{w}^{OWA}\bm{x}_{i_0}$.\\
Compute $b^{OWA}=\frac{b_1+b_2}{2}$.
}
Retrieve estimator hyperplane:
$y=\bm{w}^{OWA}\bm{x}+b^{OWA}.$
 \caption{Heuristic method based on OWA-formulas (HEUR1)
 \label{alg:heur1}
 }
 }
\end{algorithm}

\subsection{Two-step heuristic procedure for OWA-SVR (HEUR2)}
\label{sec:heur2}
The second presented heuristic method is a two-step approach.
The main characteristics of these two steps are the following:
\begin{itemize}
\item {\bf Step 1.} The primal formulation of SVR is applied obtaining solution $(\widehat{\bm{w}},\widehat{b},\widehat{\bm{\xi}},\widehat{\bm{\xi}}^*)$. Calculate the vector of deviations,
$$\widehat{\bm{\xi}+\bm{\xi}^*}=(\widehat{\xi}_1+\widehat{\xi}_1^*,\widehat{\xi}_2+\widehat{\xi}_2^*,\ldots, \widehat{\xi}_n+\widehat{\xi}_n^*).$$
Sort this vector in non-decreasing order and denote it by $\widehat{\bm{\xi}+\bm{\xi}^*}_{()}$.
\item {\bf Step 2.} Solve the SVR minimizing a balance between the norm of the hyperplane coefficients and the weighted sum of deviations based on order provided by vector $\widehat{\bm{\xi}+\bm{\xi}^*}_{()}$.
\end{itemize}

 A detailed description of mentioned steps can be found in Algorithm \ref{alg:heur2}. 

\begin{algorithm}[htbp]
{\small
	\KwData{Training sample composed by a set of $n$ individuals with $d$ features.}
	\KwResult{Feasible solution for OWA-SVR}
 Solve the primal form of classsic SVR.\\
 Let $(\widehat{\bm{w}},\widehat{b},\widehat{\bm{\xi}},\widehat{\bm{\xi}}^*)$ be the optimal solution.\\
 Compute: $\widehat{\bm{\xi}+\bm{\xi}^*}=(\widehat{\xi}_1+\widehat{\xi}_1^*,\widehat{\xi}_2+\widehat{\xi}_2^*,\ldots, \widehat{\xi}_n+\widehat{\xi}_n^*)$.\\
Define $\widehat{\bm{\xi}+\bm{\xi}^*}_{()}$ as vector $\widehat{\bm{\xi}+\bm{\xi}^*}$ sorted in non-decreasing order.\\
Compute vector of positions, $v$, related to $\widehat{\bm{\xi}+\bm{\xi}^*}_{()}$ in the sense that $\widehat{\bm{\xi}+\bm{\xi}^*}_{(i)}=\widehat{\xi}_{v(i)}+\widehat{\xi}_{v(i)}^*$.\\
Solve: 
\begin{alignat*}{3}
(\mbox{HEUR2}_{model})& \min_{\bm{w},b,\bm{\xi},\bm{\xi}^*}&  \quad & \frac{1}{2}\norm{\bm{w}}^2+C\sum_{i=1}^n \lambda_i(\xi_{v(i)}+\xi^*_{v(i)})\nonumber\\
&\mbox{s.t.}&&\eqref{SVR1}-\eqref{SVR4}.
\end{alignat*}
Let$(\widetilde{\bm{w}},\widetilde{b},\widetilde{\bm{\xi}},\widetilde{\bm{\xi}}^*)$ be the solution of the above model.\\
Retrieve regression hyperplane:
$y=\widetilde{\bm{w}}^T\bm{x}+\widetilde{b}.$
 \caption{Two-step heuristic procedure (HEUR2)
 \label{alg:heur2}
 }
 }
\end{algorithm}

\subsection{Two-step heuristic procedure for OWA-SVR considering non-linear kernels (HEUR$_K$) }
\label{2-step-kernel}

In this two-step approach, it is  necessary to take into account the theoretical development presented in Section \ref{sec:models} to determine the penalization associated with the deviation of each observation. The two steps of this procedure are:

\begin{itemize}
\item {\bf Step 1.} Given a certain kernel function, $K(\cdot,\cdot)$, the dual formulation of classic SVR is applied to obtain an optimal non-linear regression function, see \cite{Smola2004} for more details. As a result, a vector of deviations of the individuals can be calculated, $\widehat{\bm{\xi}+\bm{\xi}^*}$. This vector is sorted in non-decreasing order and denoted as $\widehat{\bm{\xi}+\bm{\xi}^*}_{()}=\widehat{\xi}_{v(i)}+\widehat{\xi}_{v(i)}^*$, where $v(\cdot)$ is the position of $\widehat{\xi}_{i}+\widehat{\xi}_{i}^*$ in the ordered vector.

\item  {\bf Step 2.} The order determined in Step 1 is used to penalized the deviations in the following model derived from NC-OWA-SVR$_K$.

\begin{alignat}{3}
(\mbox{OWA-SVR-HEUR$_K$})& \hspace{-0.1cm}\min_{\bm{\alpha},\bm{\alpha}^*,b,\bm{\xi},\bm{\xi}^*}&  \quad & \frac{1}{2}\sum_{i=1}^n\sum_{j=1}^n(\alpha_i-\alpha_i^*)(\alpha_j-\alpha_j^*)K(\bm{x}_i,\bm{x}_j)\nonumber\\
&&&+C\sum_{k=1}^n\lambda_k(\xi_{v(k)}+\xi^*_{v(k)})\nonumber\\
&\mbox{s.t.}&&\eqref{SVR3},\eqref{dual_1},
\eqref{err1},\eqref{err2}.\nonumber \\
&&&  0 \le \alpha_{i}\le C\lambda_{v(i)}, \quad i\in N, \nonumber\\
&&&  0 \le \alpha_{i}^* \le C\lambda_{v(i)}, \quad i\in N.  \nonumber
\end{alignat}
\end{itemize}

 Note that, in the above formulation, $v$ is the index vector
related to $\widehat{\bm{\xi}+\bm{\xi}^*}_{()}$ in the sense that $\widehat{\bm{\xi}+\bm{\xi}^*}_{(i)}=\widehat{\xi}_{v(i)}+\widehat{\xi}_{v(i)}^*$.

 A detailed description of this procedure can be found in Algorithm \ref{alg:heur3}. In Section \ref{sec:comp}, we show the performance of this heuristic procedure. 

\begin{algorithm}[htbp]
{\small
	\KwData{Training sample composed by a set of $n$ individuals with $d$ features.}
	\KwResult{Feasible solution for OWA-SVR considering kernel function $K(\cdot,\cdot)$}
 Solve the dual form of classic SVR with kernel $K(\cdot,\cdot)$.\\
 Given the optimal solution of this model,
 compute the vector of deviations: $\widehat{\bm{\xi}+\bm{\xi}^*}$.\\
Define $\widehat{\bm{\xi}+\bm{\xi}^*}_{()}$ as vector $\widehat{\bm{\xi}+\bm{\xi}^*}$ sorted in non-decreasing order.\\
Compute vector of positions, $v$, such that $\widehat{\bm{\xi}+\bm{\xi}^*}_{(i)}=\widehat{\xi}_{v(i)}+\widehat{\xi}_{v(i)}^*$.\\
Solve the \mbox{OWA-SVR-HEUR$_K$} model obtaining a solution 
$ (\widehat{\bm{\alpha}},\widehat{\bm{\alpha}}^*,\widehat{b},\widehat{\bm{\xi}},\widehat{\bm{\xi}^*})$.\\
Retrieve the regression function:
$y=\sum_{i \in N}(\widehat{\alpha}_i-\widehat{\alpha}_i^*)K(\bm{x}_i,\bm{x})+\widehat{b}$.
 \caption{Two-step heuristic procedure considering a given kernel $K(\cdot,\cdot)$, (HEUR$_K$).
 \label{alg:heur3}
 }
 }
\end{algorithm}

 Once explained the detailed procedures of HEUR1, HEUR2 and HEUR$_K$, some limits and potentials of these heuristics can be analyzed. Regarding HEUR1, we can appreciate that OWA weights are introduced in the  expression of the dual form of the classic SVR. Then, it is not possible to ensure that the resulting coefficients derive from any feasible solution of the exact approach. The possible limits of this aspect will be seen in Section \ref{sec:comp}.

Focusing on HEUR2 and HEUR$_K$, the strategy is to determine a certain order of the deviation vector based on  a solution of the classic approach, and then they solve an auxiliar subproblem which is simpler than the original ordered weighted average exact approach. This allows us to obtain an optimal solution for  the classical SVR regression where the penalty weight vector is defined by the  order of the deviations of the solution previously obtained. The potentials and limits of these heuristic methods will be appreciated in next section.

\section{Computational experiments } \label{sec:comp}
In this section, we present several computational studies comparing the performance of the proposed exact OWA-SVR formulations with respect to the introduced heuristic approaches.
These computational experiments
were implemented using Gurobi 10.0 in Python on an Intel(R)
Xeon(R) W-2245 CPU 256 GB RAM computer. Note that ``addGenConstrIndicator" utility from Gurobi was used to deal with big-M constraints appearing in some of the formulations.

Regarding this aim, we consider  eight datasets previously addressed in literature for regression purposes:
\begin{itemize}
    \item {\bf Quandt}. This dataset is composed by 20 
     bidimensional points $(x_i,y_i)$ for $i=1,\ldots,20$. Note that $x$-values are generated as a random permutation of the first 20 natural numbers. Then, $y$-values are generated as follows:
     \begin{eqnarray*}
     y_i&=&2.5+0.7x_i+\delta_i,\mbox{ for }i=1,\ldots,12,\\
     y_i&=&5+0.5x_i+\delta_i,\mbox{ for } i=13,\ldots,20,
     \end{eqnarray*}
    where $\delta_i$ is randomly generated as a $N(0,1)$. This dataset was initially introduced in \cite{Quandt1958} and analyzed in many regression works. 
    \item {\bf Boston}. This real dataset is related to stock exchanges, see \citep{McGee1970}. It consists on bidimensional points $(x_i,y_i)$ for $i=1,\ldots, 35$, where $x$-values are New York and American stock exchanges and $y$-values are Boston stock exchanges. This data corresponds to values from January 1967 to November 1969.
    \item {\bf CYGOB1}. The information contained in this dataset is associated with stars of the CYG OB1 cluster in the direction of Cygnus. Each $x$-value is the logarithm of the effective temperature at the star surface and each $y$-value is the logarithm of its light intensity. This dataset is composed by 47 points. For more details, see \cite{Humphreys1978}.
    \item {\bf Eilon}. The 50 bidimensional points appearing in this dataset were  fist introduced in \cite{Eilon1971} in the context of network analysis.
    \item {\bf Durbin-Watson}. This dataset, introduced in \cite{Durbin1951}, is a real sample related to the annual consumption of spirits from 1870 to 1938. It is a tridimentional dataset, $(x_{1i},x_{2i},y_i)$, for $i=1,\ldots,69$. In this case, $x_1$-values are logarithms of real income per head, $x_2$-values are logarithms of the relative price of spirits, and $y$-values are the logarithms of spirit consumption per head. This dataset is analyzed twice in the reported computational experiments. First, we analyzed the results considering the bidimensional dataset composed by $(x_1,y)$. This dataset will be denoted by \textit{DW2}. Then, the complete tridimensional data $(x_1,x_2,y)$, denoted by \textit{DW3} is studied.
    \item {\bf Alcohol}. This multivariate dataset is composed by physicochemical properties of $44$ aliphatic alcohols and the values of the following variables: $x_1$ is related to the molecular volume, $x_2$ is the logarithm of the octanol-water partitions coefficient, $x_3$ correspods to the polarizability, $x_4$ is the molar refractivity, $x_5$ is the mass, $x_6$ is the volume and $y$ represents the logarithm of the solubility. See \cite{Maronna2006} for more details.
    \item {\bf HBK}. This dataset was generated by \cite{MR0770368} and it is composed by $75$ $4$-dimension observations (three explanatory variables and one response). Thus, we consider $(x_{1i},x_{2i},x_{3i},y_{i})$ for $i=1,\ldots,75$. This dataset has also been considered for the analysis performed in \cite{Bertsimas2014}.

    \item {\bf Friedman}. This is a synthetic dataset that can be generated with different sample size and variables. See \cite{Friedman1991} for a detailed description of this family of datasets. Note that we generated three instances from this sythetic structure using {\sl sklearn} module\footnote{\url{https://scikit-learn.org/dev/modules/generated/sklearn.datasets.make_friedman1.html}} of python:
    \begin{itemize}
    \item Friedman (150): It has 8 explanatory variables and the sample size is 150.
        \item Friedman (300): It has 8 explanatory variables and the sample size is 300.
            \item Friedman (500): It has 8 explanatory variables and the sample size is 500.
        
    \end{itemize}
    
\end{itemize}

All the mentioned datasets can be found in the following repository:\\  \url{https://github.com/orel-uca/OWASVR}

Recall that, in contrast with previous regression models, OWA-SVR allows to impose a different penalization to the deviations depending on their sizes. Consequently, different weights vectors can be taken into consideration depending on the purpose or application of the decision maker. For the computational experiments, we consider the weight vectors described explicitly in Table \ref{tab:2}. These weights correspond to the following deviation aggregation criteria, already considered in \cite{BLANCO2018,Marin2022}:
\begin{itemize}
\item {\bf SUM}: Sum of all deviations, i.e., the resulting model is equivalent to the classic SVR, \cite{Smola2004}.
\item {\bf kC}: Sum of the $k$ largest deviations.
\item {\bf MAX}: It only considers the largest deviation.
\item {\bf akC}: Sum of the $k$ smallest deviations.
\item {\bf MED}: It only considers the median deviation.
\item {\bf TM} (trimmed-mean): Sum of the central deviation values. In this weight vector smallest and largest deviations are not taken into account.
\end{itemize}

\begin{table}[htbp]
  \centering
  \caption{Weights used in the computational experiments}
    \begin{tabular}{ll}
    \hline
    \hline
    Name  & Weight \\
    \hline
    \hline
    SUM  & $\lambda_1=\lambda_2=\ldots=\lambda_n=1$  \\
    kC   &$\lambda=\ldots=\lambda_{\lfloor\frac{n}{2}\rfloor}=0, \: \lambda_{\lceil\frac{n}{2}\rceil}=\ldots=\lambda_n=1 $\\
    MAX&$\lambda_1=\ldots=\lambda_{n-1}=0, \: \lambda_n=1$\\
    akC&$\lambda_1=\ldots=\lambda_{\lfloor\frac{n}{2}\rfloor}=1,\lambda_{\lceil\frac{n}{2}\rceil}=\ldots=\lambda_n=0 $\\
    MED&$\lambda_1=\ldots=\lambda_{\lfloor\frac{n}{2}\rfloor}=0,\: \lambda_{\lceil\frac{n}{2}\rceil}=1, \:\:\lambda_{\lceil\frac{n}{2}\rceil+1}=\ldots=\lambda_n=0 $\\
    TM&$\lambda_1=\ldots=\lambda_{\lfloor\frac{n}{4}\rfloor}=0, \: \lambda_{\lceil\frac{n}{4}\rceil}=\ldots\lambda_{\lfloor\frac{3n}{4}\rfloor-1}=1,\: \lambda_{\lfloor\frac{3n}{4}\rfloor}=\ldots=\lambda_n=0 $\\
    \hline
    \hline
    \end{tabular}%
      \label{tab:2}%
\end{table}%

Table \ref{tab:2} details that, for kC weight, one half of the deviations with the largest values are considered. Similarly, we consider one half of individuals with smallest values for the akC weight. Finally, in the case of TM weight, one half of the central deviation values are taken into account.

 This computational section is divided in five main studies. In the first four sections (\ref{gof}, \ref{performance}, \ref{graphical}, \ref{time}), we focus on Quandt, Boston, CYGOB1, Eilon, Durbin-Watson, Alcohol and HBK, respectively. First, we extend the goodness of fitting index definition introduced in \cite{BLANCO2018} to evaluate the deviations of these approaches with respect the classic SVR.  
Then, a fold-cross validation is performed  on each dataset comparing different classic performance measures. Besides, some graphical insights about the evaluation metrics are reported and a computational time analysis is included. Finally, in Section \ref{friedman} we focus on Friedman datasets which present the biggest sample sizes. 

 Table \ref{tab:1} includes a sum up of the different considered datasets,  the size of each fold in the FCV and the size of the test set.

\begin{table}[htbp]
  \centering
        \caption{Characteristics of datasets considered for the computational experience}
    \begin{tabular}{lrrrr}
    \hline
    \hline
    Data  & \multicolumn{1}{l}{Number of } & \multicolumn{1}{l}{Explanatory} & \multicolumn{1}{l}{Folds size } & \multicolumn{1}{l}{Test set} \\
    name  & \multicolumn{1}{l}{observations} & \multicolumn{1}{l}{variables} & \multicolumn{1}{l}{in FCV} & \multicolumn{1}{l}{size}\\
    \hline
    \hline
    Quandt & 20    & 1     & 5&5 \\
    Boston & 35    & 1     & 5/6&9 \\
    CYGOB1 & 47    & 1     & 7 &12\\
    Eilon & 50    & 1      & 7/8   &13  \\
    DW2 & 69    & 1        & 10/11 &18 \\
    Alcohol & 44    & 6     & 6/7&11 \\
    DW3 & 69    & 2         & 10/11  & 18 \\
    HBK   & 75    & 3     & 11/12 &19\\
    \hline
    \hline
    \end{tabular}%
  \label{tab:1}%
\end{table}%

\subsection{Goodness of fitting index}\label{gof}

 The goodness of fitting index (GoF) was introduced in \cite{BLANCO2018}. It extends the coefficient of determination definition for the case where OWA regression is considered. It evaluates the optimal hyperplane under an specific criterion with respect to the best horizontal hyperplane. Generalizing this notion to the SVR with OWA, GoF could be defined as follows:
 
  $$GoF=1-\frac{\phi^*}{\phi^0},$$
where:
   \begin{eqnarray*}
\phi^*&=&\sum_{k=1}^n\lambda_k\theta_k^*,\\ \phi^0&=&\sum_{k=1}^n\lambda_k\theta_{k}^0.\\ 
   \end{eqnarray*}
In this definition, $\theta^*$ is the ordered vector of deviations with respect to the strip defined by $\epsilon$ in the evaluated model. Besides $\theta^0$ is the sorted vector of deviations to the strip induced by $\epsilon$ obtained using the hyperplane given by the optimal solution of SVR considering only the independent term ($b$).

Note that this measure is well defined if $\phi^0>0$. Particularly,
when GoF is close to $0$, the model does not provide  a meaningful improvement with respect to the best horizontal hyperplane, since $\phi^*\simeq \phi^0$. However, when GoF is close to one, the model is significantly better than the best horizontal hyperplane  derived from SVR, since it means that $\phi^*\ll\phi^0$. Furthermore, if GoF$=1$, then $\phi^*=0$ and the model provides no deviations. Then, we can conclude that we obtain a good performance of the model when this measure is  close to one.

Table \ref{tab:gof2} reports the GoF index for all developed approaches along this work. First column of this table presents the names of the analyzed datasets and the second one reports the considered weights previously described in Table \ref{tab:2}.  Then, we report the best results provided by each approach for these weights. Note that column block HEUR1 reports the results for the first proposed heuristic procedure described in Section \ref{sec:heur1} and column block HEUR2 reports the results related to the heuristic in Section \ref{sec:heur2}. Column block named as OWA-SVR shows the results of the exact formulation. Note that, depending on the considered weight a different formulation  (\con) or (\ncon)  must be chosen to solve the problem. Besides, column block HEUR$_K$ reports the results of the heuristic approach described in Section \ref{2-step-kernel} considering the Gaussian Kernel, which is given by:

$$K(\bm{x}_i,\bm{x}_j)=exp\left(\frac{\norm{\bm{x}_i-\bm{x}_j}}{2\sigma^2}\right),\,\,\,\sigma>0.$$

Finally, column block OWA-SVR$_K$ shows the results related with the exact approaches \ncon$_K$ and \con$_K$ described in  Section \ref{sec:kernel}  and considering the Gaussian Kernel.

In each column block, we observe two columns. The first column consists in a tuple with the hyperparameter combination that provides the best GoF and the computing time ($t$).  We set a time limit of one hour that is denoted by $T$ in Table \ref{tab:gof2}. The second column of the  block, reports the GoF. In our computational experience, we have considered the following hyperparameter values  described in Table \ref{tab:hyperparameters}.

\begin{table}[htbp]
  \centering
  \caption{Considered hyperparameter values.}
    \begin{tabular}{cl}
    \hline
    \hline
    Hyperparameter & Values\\
        \hline
    \hline
    $C$     &  $\{0.01,0.1,1,10,100\}$\\
    $\varepsilon$ & $\{0,0.01,0.1,1,10\}$ \\
    $\sigma$ & $\{0.01,0.1,1,10,100\}$\\
    \hline
    \hline
    \end{tabular}%
  \label{tab:hyperparameters}%
\end{table}%

The color codes in these columns inform about the best results per row. Gray color reports the best results among the different approaches using linear Kernel, i.e., HEUR1, HEUR2 and OWA-SVR. Blue color reports the best results between two Gaussian kernel approaches, HEUR$_K$ and OWA-SVR$_K$, when these improve the ones of the linear approaches.

From Table \ref{tab:gof2} the reader can realize that in many cases OWA-SVR$_K$ provides the best GoF, with the exception of few cases where the best results are attained by HEUR$_K$. In addition, in all but one case (Quandt dataset with the TM criterion), OWA-SVR provides at least the second best GoF among the studied approaches, which is consistently an assessment of the quality of the proposed method. In any case, either OWA-SVR or OWA-SVR$_K$ consistently provide GoF indices close to 1, which shows the robustness of our two exact approaches for solving the OWA-SVR problem in the 8 considered datasets.  On the negative side, HEUR1 is the method that performs worse in terms of the GoF index.

The reader may note that the GoF index should rank in $[0,1]$. Nevertheless, in some cases where the problems related to $\phi^*$ are not solved to optimality, we face cases with $\phi^0<\phi^*$ which results in values smaller than 0. This is often the case of HEUR1 and HEUR2 which fail to provide meaningful results, as compared with the exact OWA-SVR methods.

\begin{table}[htbp]
  \centering
  \caption{GoF for the different datasets and solution approaches}
   \begin{adjustwidth}{-3.5cm}{0cm}
   \scalebox{0.75}{
  
    \begin{tabular}{|c|l|lr|lr|lr|lr|lr|}
    \hline
    \multirow{2}[2]{*}{Data} & \multicolumn{1}{c|}{\multirow{2}[2]{*}{Weight}} & \multicolumn{2}{c|}{HEUR1} & \multicolumn{2}{c|}{HEUR2} & \multicolumn{2}{c|}{OWA-SVR} & \multicolumn{2}{c|}{HEUR$_K$} & \multicolumn{2}{c|}{OWA-SVR$_K$} \\
          &       & $(C,\epsilon,t)$ & \multicolumn{1}{l|}{GoF} & $( C,\epsilon , t )$ & \multicolumn{1}{l|}{GoF} & $(C,\epsilon,t)$ & \multicolumn{1}{l|}{GoF} & $(C,\epsilon,\sigma,t)$ & \multicolumn{1}{l|}{GoF} & $(C,\epsilon,\sigma,t)$ & \multicolumn{1}{l|}{GoF} \\
    \hline
    \multirow{6}[2]{*}{\begin{sideways}Quandt\end{sideways}} &  SUM  & \cellcolor[rgb]{ .788,  .788,  .788}\textbf{(0.1,1,0)} & \cellcolor[rgb]{ .788,  .788,  .788}\textbf{0.89} & \cellcolor[rgb]{ .788,  .788,  .788}\textbf{(0.1,1,0)} & \cellcolor[rgb]{ .788,  .788,  .788}\textbf{0.89} & \cellcolor[rgb]{ .788,  .788,  .788}\textbf{(0.1,1,0)} & \cellcolor[rgb]{ .788,  .788,  .788}\textbf{0.89} & \cellcolor[rgb]{ .851,  .882,  .949}\textbf{(10,0,0.01,0)} & \cellcolor[rgb]{ .851,  .882,  .949}\textbf{1.00} & \cellcolor[rgb]{ .851,  .882,  .949}\textbf{(10,0,0.01,0)} & \cellcolor[rgb]{ .851,  .882,  .949}\textbf{1.00} \\
          &  kC   & (0.01,1,0) & 0.71  & \cellcolor[rgb]{ .788,  .788,  .788}\textbf{(0.1,1,0)} & \cellcolor[rgb]{ .788,  .788,  .788}\textbf{0.87} & \cellcolor[rgb]{ .788,  .788,  .788}\textbf{(0.1,1,0)} & \cellcolor[rgb]{ .788,  .788,  .788}\textbf{0.87} & (100,1,10,0) & 0.88  & \cellcolor[rgb]{ .851,  .882,  .949}\textbf{(10,0,0.01,0)} & \cellcolor[rgb]{ .851,  .882,  .949}\textbf{1.00} \\
          &  MAX  & (0.01,0,0) & 0.22  & (0.1,0,0) & -0.34 & \cellcolor[rgb]{ .788,  .788,  .788}\textbf{(1,1,0)} & \cellcolor[rgb]{ .788,  .788,  .788}\textbf{0.79} & (10,0,10,0) & 0.00  & \cellcolor[rgb]{ .851,  .882,  .949}\textbf{(100,0,0.01,0)} & \cellcolor[rgb]{ .851,  .882,  .949}\textbf{1.00} \\
          &  akC  & (0.01,1,0) & -1.42 & \cellcolor[rgb]{ .788,  .788,  .788}\textbf{(0.1,1,0)} & \cellcolor[rgb]{ .788,  .788,  .788}\textbf{1.00} & \cellcolor[rgb]{ .788,  .788,  .788}\textbf{(0.1,1,0.3)} & \cellcolor[rgb]{ .788,  .788,  .788}\textbf{1.00} & (1,1,1,0) & 0.85  & \cellcolor[rgb]{ .851,  .882,  .949}\textbf{(100,0.1,10,5.8)} & \cellcolor[rgb]{ .851,  .882,  .949}\textbf{0.92} \\
          &  MEDIAN & (0.01,0.1,0) & -0.12 & (0.01,1,0) & -0.02 & \cellcolor[rgb]{ .788,  .788,  .788}\textbf{(1,1,0.8)} & \cellcolor[rgb]{ .788,  .788,  .788}\textbf{1.00} & (10,0.1,0.01,0) & 0.02  & (100,0.1,10,0.2) & 0.56 \\
          & TM    & (0.01,0.01,0) & 0.39  & \cellcolor[rgb]{ .788,  .788,  .788}\textbf{(1,1,0)} & \cellcolor[rgb]{ .788,  .788,  .788}\textbf{1.00} & (0.1,1,1.6) & 0.99  & \cellcolor[rgb]{ .851,  .882,  .949}\textbf{(100,1,10,0)} & \cellcolor[rgb]{ .851,  .882,  .949}\textbf{1.00} & (100,1,10,2.3) & 0.98 \\
    \hline
    \multirow{6}[2]{*}{\begin{sideways}Boston\end{sideways}} &  SUM  & \cellcolor[rgb]{ .788,  .788,  .788}\textbf{(10,0.1,0)} & \cellcolor[rgb]{ .788,  .788,  .788}\textbf{0.41} & \cellcolor[rgb]{ .788,  .788,  .788}\textbf{(10,0.1,0)} & \cellcolor[rgb]{ .788,  .788,  .788}\textbf{0.41} & \cellcolor[rgb]{ .788,  .788,  .788}\textbf{(10,0.1,0.1)} & \cellcolor[rgb]{ .788,  .788,  .788}\textbf{0.41} & \cellcolor[rgb]{ .851,  .882,  .949}\textbf{(100,0.1,0.01,0)} & \cellcolor[rgb]{ .851,  .882,  .949}\textbf{1.00} & \cellcolor[rgb]{ .851,  .882,  .949}\textbf{(100,0.1,0.01,0)} & \cellcolor[rgb]{ .851,  .882,  .949}\textbf{1.00} \\
          &  kC   & (0.1,0.1,0) & 0.37  & \cellcolor[rgb]{ .788,  .788,  .788}\textbf{(10,0.1,0)} & \cellcolor[rgb]{ .788,  .788,  .788}\textbf{0.40} & \cellcolor[rgb]{ .788,  .788,  .788}\textbf{(10,0.1,0.1)} & \cellcolor[rgb]{ .788,  .788,  .788}\textbf{0.40} & (1,0.1,0.01,0) & 0.72  & \cellcolor[rgb]{ .851,  .882,  .949}\textbf{(100,0.1,0.01,0)} & \cellcolor[rgb]{ .851,  .882,  .949}\textbf{1.00} \\
          &  MAX  & (1,0,0) & 0.18  & (0.01,0,0) & -0.32 & \cellcolor[rgb]{ .788,  .788,  .788}\textbf{(100,0.1,0.1)} & \cellcolor[rgb]{ .788,  .788,  .788}\textbf{0.59} & (1,0,0.1,0) & 0.21  & \cellcolor[rgb]{ .851,  .882,  .949}\textbf{(100,0.1,0.01,0)} & \cellcolor[rgb]{ .851,  .882,  .949}\textbf{1.00} \\
          &  akC  & (0.01,0.1,0) & 0.37  & \cellcolor[rgb]{ .788,  .788,  .788}\textbf{(0.1,0.1,0)} & \cellcolor[rgb]{ .788,  .788,  .788}\textbf{1.00} & \cellcolor[rgb]{ .788,  .788,  .788}\textbf{(0.01,0.1,0.4)} & \cellcolor[rgb]{ .788,  .788,  .788}\textbf{1.00} & \cellcolor[rgb]{ .851,  .882,  .949}\textbf{(0.1,0.1,0.01,0)} & \cellcolor[rgb]{ .851,  .882,  .949}\textbf{1.00} & \cellcolor[rgb]{ .851,  .882,  .949}\textbf{(0.01,0.1,0.01,0.4)} & \cellcolor[rgb]{ .851,  .882,  .949}\textbf{1.00} \\
          &  MEDIAN & (0.01,0.1,0) & 0.39  & (0.1,0.1,0) & 0.98  & \cellcolor[rgb]{ .788,  .788,  .788}\textbf{(0.01,0.1,0.5)} & \cellcolor[rgb]{ .788,  .788,  .788}\textbf{1.00} & \cellcolor[rgb]{ .851,  .882,  .949}\textbf{(0.1,0.1,1,0)} & \cellcolor[rgb]{ .851,  .882,  .949}\textbf{1.00} & \cellcolor[rgb]{ .851,  .882,  .949}\textbf{(0.01,0.1,0.1,0.1)} & \cellcolor[rgb]{ .851,  .882,  .949}\textbf{1.00} \\
          & TM    & (0.01,0.01,0) & 0.02  & (10,0.1,0) & 0.77  & \cellcolor[rgb]{ .788,  .788,  .788}\textbf{(100,0.1,137.8)} & \cellcolor[rgb]{ .788,  .788,  .788}\textbf{0.78} & \cellcolor[rgb]{ .851,  .882,  .949}\textbf{(1,0.1,0.01,0)} & \cellcolor[rgb]{ .851,  .882,  .949}\textbf{1.00} & (10,0,0.01,T) & 0.97 \\
    \hline
    \multirow{6}[2]{*}{\begin{sideways}CYGOB1\end{sideways}} &  SUM  & \cellcolor[rgb]{ .788,  .788,  .788}\textbf{(1,1,0)} & \cellcolor[rgb]{ .788,  .788,  .788}\textbf{1.00} & \cellcolor[rgb]{ .788,  .788,  .788}\textbf{(1,1,0)} & \cellcolor[rgb]{ .788,  .788,  .788}\textbf{1.00} & \cellcolor[rgb]{ .788,  .788,  .788}\textbf{(1,1,0.3)} & \cellcolor[rgb]{ .788,  .788,  .788}\textbf{1.00} & \cellcolor[rgb]{ .851,  .882,  .949}\textbf{(1,1,0.01,0)} & \cellcolor[rgb]{ .851,  .882,  .949}\textbf{1.00} & \cellcolor[rgb]{ .851,  .882,  .949}\textbf{(1,1,0.01,0)} & \cellcolor[rgb]{ .851,  .882,  .949}\textbf{1.00} \\
          &  kC   & (0.01,0,0) & -0.07 & \cellcolor[rgb]{ .788,  .788,  .788}\textbf{(1,1,0)} & \cellcolor[rgb]{ .788,  .788,  .788}\textbf{1.00} & \cellcolor[rgb]{ .788,  .788,  .788}\textbf{(1,1,0.3)} & \cellcolor[rgb]{ .788,  .788,  .788}\textbf{1.00} & \cellcolor[rgb]{ .851,  .882,  .949}\textbf{(1,1,0.01,0)} & \cellcolor[rgb]{ .851,  .882,  .949}\textbf{1.00} & \cellcolor[rgb]{ .851,  .882,  .949}\textbf{(1,1,0.01,0)} & \cellcolor[rgb]{ .851,  .882,  .949}\textbf{1.00} \\
          &  MAX  & (0.01,0,0) & -0.04 & (1,1,0) & -0.46 & \cellcolor[rgb]{ .788,  .788,  .788}\textbf{(10,1,0.2)} & \cellcolor[rgb]{ .788,  .788,  .788}\textbf{1.00} & (1,0.01,0.01,0) & 0.02  & \cellcolor[rgb]{ .851,  .882,  .949}\textbf{(1,1,0.01,0)} & \cellcolor[rgb]{ .851,  .882,  .949}\textbf{1.00} \\
          &  akC  & (0.01,0.01,0) & 0.16  & (10,0.1,0) & 0.34  & \cellcolor[rgb]{ .788,  .788,  .788}\textbf{(1,0.1,T)} & \cellcolor[rgb]{ .788,  .788,  .788}\textbf{0.40} & \cellcolor[rgb]{ .851,  .882,  .949}\textbf{(1,0.1,0.01,0)} & \cellcolor[rgb]{ .851,  .882,  .949}\textbf{1.00} & (100,0.1,0.1,T) & 0.98 \\
          &  MEDIAN & (10,0.1,0) & 0.03  & (10,0,0) & 0.02  & \cellcolor[rgb]{ .788,  .788,  .788}\textbf{(100,0.1,377.6)} & \cellcolor[rgb]{ .788,  .788,  .788}\textbf{0.47} & (0.1,0.1,0.1,0) & 0.14  & \cellcolor[rgb]{ .851,  .882,  .949}\textbf{(100,0,0.01,1400.8)} & \cellcolor[rgb]{ .851,  .882,  .949}\textbf{0.82} \\
          & TM    & (0.1,0,0) & 0.08  & (1,0.1,0) & 0.13  & \cellcolor[rgb]{ .788,  .788,  .788}\textbf{(10,0.1,T)} & \cellcolor[rgb]{ .788,  .788,  .788}\textbf{0.38} & \cellcolor[rgb]{ .851,  .882,  .949}\textbf{(10,0.1,0.01,0)} & \cellcolor[rgb]{ .851,  .882,  .949}\textbf{0.99} & \cellcolor[rgb]{ .851,  .882,  .949}\textbf{(10,0.1,0.01,T)} & \cellcolor[rgb]{ .851,  .882,  .949}\textbf{0.99} \\
    \hline
    \multirow{6}[2]{*}{\begin{sideways}Eilon\end{sideways}} &  SUM  & \cellcolor[rgb]{ .788,  .788,  .788}\textbf{(1,0.1,0)} & \cellcolor[rgb]{ .788,  .788,  .788}\textbf{0.01} & \cellcolor[rgb]{ .788,  .788,  .788}\textbf{(1,0.1,0)} & \cellcolor[rgb]{ .788,  .788,  .788}\textbf{0.01} & \cellcolor[rgb]{ .788,  .788,  .788}\textbf{(1,0.1,0.4)} & \cellcolor[rgb]{ .788,  .788,  .788}\textbf{0.01} & \cellcolor[rgb]{ .851,  .882,  .949}\textbf{(100,0.1,0.01,0)} & \cellcolor[rgb]{ .851,  .882,  .949}\textbf{1.00} & \cellcolor[rgb]{ .851,  .882,  .949}\textbf{(100,0.1,0.01,0)} & \cellcolor[rgb]{ .851,  .882,  .949}\textbf{1.00} \\
          &  kC   & (0.01,0.1,0) & 0.00  & (1,0.1,0) & -0.36 & \cellcolor[rgb]{ .788,  .788,  .788}\textbf{(0.01,0.1,0.4)} & \cellcolor[rgb]{ .788,  .788,  .788}\textbf{0.04} & (10,0.1,0.01,0) & 0.37  & \cellcolor[rgb]{ .851,  .882,  .949}\textbf{(100,0.1,0.01,0)} & \cellcolor[rgb]{ .851,  .882,  .949}\textbf{1.00} \\
          &  MAX  & (0.01,0.01,0) & 0.03  & (0.01,0.01,0) & -0.91 & \cellcolor[rgb]{ .788,  .788,  .788}\textbf{(10,0.01,0.4)} & \cellcolor[rgb]{ .788,  .788,  .788}\textbf{0.06} & (100,0.1,0.01,0) & -0.09 & \cellcolor[rgb]{ .851,  .882,  .949}\textbf{(100,0.1,0.01,0)} & \cellcolor[rgb]{ .851,  .882,  .949}\textbf{0.67} \\
          &  akC  & (0.01,0,0) & -0.02 & (0.01,0.1,0) & 0.09  & \cellcolor[rgb]{ .788,  .788,  .788}\textbf{(1,0.1,T)} & \cellcolor[rgb]{ .788,  .788,  .788}\textbf{0.20} & \cellcolor[rgb]{ .851,  .882,  .949}\textbf{(1,0,0.01,0)} & \cellcolor[rgb]{ .851,  .882,  .949}\textbf{1.00} & (100,0.1,0.1,T) & 0.97 \\
          &  MEDIAN & (0.1,0.1,0) & 0.11  & (0.01,0,0) & 0.00  & \cellcolor[rgb]{ .788,  .788,  .788}\textbf{(10,0.1,248.2)} & \cellcolor[rgb]{ .788,  .788,  .788}\textbf{0.40} & (100,0,0.1,1800.1) & 0.15  & \cellcolor[rgb]{ .851,  .882,  .949}\textbf{(10,0,0.01,75.9)} & \cellcolor[rgb]{ .851,  .882,  .949}\textbf{0.61} \\
          & TM    & (0.1,0.1,0) & 0.03  & \cellcolor[rgb]{ .788,  .788,  .788}\textbf{(0.01,0.1,0)} & \cellcolor[rgb]{ .788,  .788,  .788}\textbf{0.12} & \cellcolor[rgb]{ .788,  .788,  .788}\textbf{(0.01,0.1,T)} & \cellcolor[rgb]{ .788,  .788,  .788}\textbf{0.12} & (10,0.1,0.01,0) & 0.95  & \cellcolor[rgb]{ .851,  .882,  .949}\textbf{(1,0.01,0.01,T)} & \cellcolor[rgb]{ .851,  .882,  .949}\textbf{0.97} \\
    \hline
    \multirow{6}[2]{*}{\begin{sideways}DW2\end{sideways}} &  SUM  & \cellcolor[rgb]{ .788,  .788,  .788}\textbf{(100,0.1,0)} & \cellcolor[rgb]{ .788,  .788,  .788}\textbf{0.21} & \cellcolor[rgb]{ .788,  .788,  .788}\textbf{(100,0.1,0)} & \cellcolor[rgb]{ .788,  .788,  .788}\textbf{0.21} & \cellcolor[rgb]{ .788,  .788,  .788}\textbf{(100,0.1,0.2)} & \cellcolor[rgb]{ .788,  .788,  .788}\textbf{0.21} & \cellcolor[rgb]{ .851,  .882,  .949}\textbf{(100,0.1,0.01,0.1)} & \cellcolor[rgb]{ .851,  .882,  .949}\textbf{0.83} & \cellcolor[rgb]{ .851,  .882,  .949}\textbf{(100,0.1,0.01,0)} & \cellcolor[rgb]{ .851,  .882,  .949}\textbf{0.83} \\
          &  kC   & (0.01,0.1,0) & -0.06 & (100,0.1,0) & 0.11  & \cellcolor[rgb]{ .788,  .788,  .788}\textbf{(100,0.1,0.7)} & \cellcolor[rgb]{ .788,  .788,  .788}\textbf{0.25} & (100,0.1,0.1,0.1) & 0.71  & \cellcolor[rgb]{ .851,  .882,  .949}\textbf{(100,0.1,0.01,0.1)} & \cellcolor[rgb]{ .851,  .882,  .949}\textbf{0.83} \\
          &  MAX  & (0.01,0.1,0) & 0.19  & (0.01,0.1,0) & -0.16 & \cellcolor[rgb]{ .788,  .788,  .788}\textbf{(1,0.1,0.2)} & \cellcolor[rgb]{ .788,  .788,  .788}\textbf{0.41} & (1,0.1,0.01,0) & 0.32  & (100,0.1,0.01,0) & 0.40 \\
          &  akC  & (0.01,0.01,0) & 0.12  & \cellcolor[rgb]{ .788,  .788,  .788}\textbf{(1,0.1,0)} & \cellcolor[rgb]{ .788,  .788,  .788}\textbf{1.00} & \cellcolor[rgb]{ .788,  .788,  .788}\textbf{(0.01,0.1,7.3)} & \cellcolor[rgb]{ .788,  .788,  .788}\textbf{1.00} & \cellcolor[rgb]{ .851,  .882,  .949}\textbf{(0.01,0.1,0.01,0)} & \cellcolor[rgb]{ .851,  .882,  .949}\textbf{1.00} & \cellcolor[rgb]{ .851,  .882,  .949}\textbf{(0.01,0.1,0.01,4.5)} & \cellcolor[rgb]{ .851,  .882,  .949}\textbf{1.00} \\
          &  MEDIAN & (1,0.01,0) & 0.29  & (1,0.1,0) & 0.51  & \cellcolor[rgb]{ .788,  .788,  .788}\textbf{(0.01,0.1,0.7)} & \cellcolor[rgb]{ .788,  .788,  .788}\textbf{1.00} & \cellcolor[rgb]{ .851,  .882,  .949}\textbf{(1,0.1,1,0.1)} & \cellcolor[rgb]{ .851,  .882,  .949}\textbf{1.00} & \cellcolor[rgb]{ .851,  .882,  .949}\textbf{(0.01,0.1,0.01,0.6)} & \cellcolor[rgb]{ .851,  .882,  .949}\textbf{1.00} \\
          & TM    & (0.1,0,0) & 0.14  & (10,0.1,0) & 0.51  & \cellcolor[rgb]{ .788,  .788,  .788}\textbf{(100,0.1,T)} & \cellcolor[rgb]{ .788,  .788,  .788}\textbf{0.52} & \cellcolor[rgb]{ .851,  .882,  .949}\textbf{(1,0.1,0.01,0)} & \cellcolor[rgb]{ .851,  .882,  .949}\textbf{1.00} & (0.1,0.01,0.01,T) & 0.56 \\
    \hline
    \multirow{6}[2]{*}{\begin{sideways}alcohol\end{sideways}} &  SUM  & (100,0.1,0) & 0.91  & \cellcolor[rgb]{ .788,  .788,  .788}\textbf{(0.01,1,0)} & \cellcolor[rgb]{ .788,  .788,  .788}\textbf{1.00} & \cellcolor[rgb]{ .788,  .788,  .788}\textbf{(0.01,1,0.2)} & \cellcolor[rgb]{ .788,  .788,  .788}\textbf{1.00} & (10,0,0.01,0) & 0.99  & (10,0,0.01,0) & 0.99 \\
          &  kC   & (0.1,1,0) & -15.00 & \cellcolor[rgb]{ .788,  .788,  .788}\textbf{(0.01,1,0)} & \cellcolor[rgb]{ .788,  .788,  .788}\textbf{1.00} & \cellcolor[rgb]{ .788,  .788,  .788}\textbf{(0.01,1,0.2)} & \cellcolor[rgb]{ .788,  .788,  .788}\textbf{1.00} & \cellcolor[rgb]{ .851,  .882,  .949}\textbf{(10,1,10,0)} & \cellcolor[rgb]{ .851,  .882,  .949}\textbf{1.00} & (10,0.01,1,0) & 0.99 \\
          &  MAX  & (0.1,1,0) & -10.97 & (0.01,0,0) & 0.01  & \cellcolor[rgb]{ .788,  .788,  .788}\textbf{(0.01,1,0.2)} & \cellcolor[rgb]{ .788,  .788,  .788}\textbf{1.00} & (100,0.01,1,0) & 0.30  & (100,0.1,1,0) & 0.99 \\
          &  akC  & (0.1,1,0) & -22.38 & \cellcolor[rgb]{ .788,  .788,  .788}\textbf{(0.01,1,0)} & \cellcolor[rgb]{ .788,  .788,  .788}\textbf{1.00} & \cellcolor[rgb]{ .788,  .788,  .788}\textbf{(0.01,1,10.7)} & \cellcolor[rgb]{ .788,  .788,  .788}\textbf{1.00} & \cellcolor[rgb]{ .851,  .882,  .949}\textbf{(1,1,100,0)} & \cellcolor[rgb]{ .851,  .882,  .949}\textbf{1.00} & (10,0,10,1800.1) & 0.99 \\
          &  MEDIAN & (0.01,1,0) & 0.40  & (1,0,0) & 0.33  & \cellcolor[rgb]{ .788,  .788,  .788}\textbf{(0.01,1,11.8)} & \cellcolor[rgb]{ .788,  .788,  .788}\textbf{1.00} & (10,1,0.1,0) & 0.58  & (100,0.01,10,280.2) & 0.89 \\
          & TM    & (0.01,1,0) & 0.19  & \cellcolor[rgb]{ .788,  .788,  .788}\textbf{(0.01,1,0)} & \cellcolor[rgb]{ .788,  .788,  .788}\textbf{1.00} & \cellcolor[rgb]{ .788,  .788,  .788}\textbf{(0.01,1,53.7)} & \cellcolor[rgb]{ .788,  .788,  .788}\textbf{1.00} & \cellcolor[rgb]{ .851,  .882,  .949}\textbf{(1,1,100,0)} & \cellcolor[rgb]{ .851,  .882,  .949}\textbf{1.00} & \cellcolor[rgb]{ .851,  .882,  .949}\textbf{(10,1,100,40.6)} & \cellcolor[rgb]{ .851,  .882,  .949}\textbf{1.00} \\
    \hline
    \multirow{6}[2]{*}{\begin{sideways}DW3\end{sideways}} &  SUM  & \cellcolor[rgb]{ .788,  .788,  .788}\textbf{(10,0.1,0)} & \cellcolor[rgb]{ .788,  .788,  .788}\textbf{0.99} & \cellcolor[rgb]{ .788,  .788,  .788}\textbf{(10,0.1,0)} & \cellcolor[rgb]{ .788,  .788,  .788}\textbf{0.99} & \cellcolor[rgb]{ .788,  .788,  .788}\textbf{(10,0.1,0.2)} & \cellcolor[rgb]{ .788,  .788,  .788}\textbf{0.99} & \cellcolor[rgb]{ .851,  .882,  .949}\textbf{(1,0,0.01,0.1)} & \cellcolor[rgb]{ .851,  .882,  .949}\textbf{1.00} & \cellcolor[rgb]{ .851,  .882,  .949}\textbf{(1,0,0.01,0)} & \cellcolor[rgb]{ .851,  .882,  .949}\textbf{1.00} \\
          &  kC   & (0.01,0.1,0) & -0.57 & \cellcolor[rgb]{ .788,  .788,  .788}\textbf{(10,0.1,0)} & \cellcolor[rgb]{ .788,  .788,  .788}\textbf{0.99} & \cellcolor[rgb]{ .788,  .788,  .788}\textbf{(10,0.1,0.2)} & \cellcolor[rgb]{ .788,  .788,  .788}\textbf{0.99} & \cellcolor[rgb]{ .851,  .882,  .949}\textbf{(1,0.1,0.1,0.1)} & \cellcolor[rgb]{ .851,  .882,  .949}\textbf{1.00} & \cellcolor[rgb]{ .851,  .882,  .949}\textbf{(1,0,0.01,0)} & \cellcolor[rgb]{ .851,  .882,  .949}\textbf{1.00} \\
          &  MAX  & (0.1,0,0) & 0.20  & (1,0,0) & 0.19  & \cellcolor[rgb]{ .788,  .788,  .788}\textbf{(100,0.1,0.3)} & \cellcolor[rgb]{ .788,  .788,  .788}\textbf{0.96} & (1,0,0.01,0) & 0.21  & \cellcolor[rgb]{ .851,  .882,  .949}\textbf{(10,0.1,0.01,0)} & \cellcolor[rgb]{ .851,  .882,  .949}\textbf{1.00} \\
          &  akC  & (0.01,0.01,0) & 0.20  & \cellcolor[rgb]{ .788,  .788,  .788}\textbf{(1,0.1,0)} & \cellcolor[rgb]{ .788,  .788,  .788}\textbf{1.00} & \cellcolor[rgb]{ .788,  .788,  .788}\textbf{(0.01,0.1,8.9)} & \cellcolor[rgb]{ .788,  .788,  .788}\textbf{1.00} & \cellcolor[rgb]{ .851,  .882,  .949}\textbf{(0.01,0.1,0.01,0)} & \cellcolor[rgb]{ .851,  .882,  .949}\textbf{1.00} & \cellcolor[rgb]{ .851,  .882,  .949}\textbf{(0.01,0.1,0.01,7)} & \cellcolor[rgb]{ .851,  .882,  .949}\textbf{1.00} \\
          &  MEDIAN & \cellcolor[rgb]{ .788,  .788,  .788}\textbf{(0.1,0.1,0)} & \cellcolor[rgb]{ .788,  .788,  .788}\textbf{1.00} & (100,0,0) & 0.27  & \cellcolor[rgb]{ .788,  .788,  .788}\textbf{(0.01,0.1,0.7)} & \cellcolor[rgb]{ .788,  .788,  .788}\textbf{1.00} & (1,0.1,0.1,0.1) & 0.83  & \cellcolor[rgb]{ .851,  .882,  .949}\textbf{(0.01,0.1,0.01,0.5)} & \cellcolor[rgb]{ .851,  .882,  .949}\textbf{1.00} \\
          & TM    & (0.1,0.01,0) & 0.29  & \cellcolor[rgb]{ .788,  .788,  .788}\textbf{(10,0.1,0)} & \cellcolor[rgb]{ .788,  .788,  .788}\textbf{1.00} & \cellcolor[rgb]{ .788,  .788,  .788}\textbf{(10,0.1,1670.9)} & \cellcolor[rgb]{ .788,  .788,  .788}\textbf{1.00} & \cellcolor[rgb]{ .851,  .882,  .949}\textbf{(1,0.1,0.1,0.1)} & \cellcolor[rgb]{ .851,  .882,  .949}\textbf{1.00} & \cellcolor[rgb]{ .851,  .882,  .949}\textbf{(10,0,0.01,T)} & \cellcolor[rgb]{ .851,  .882,  .949}\textbf{1.00} \\
    \hline
    \multirow{6}[2]{*}{\begin{sideways}HBK\end{sideways}} &  SUM  & \cellcolor[rgb]{ .788,  .788,  .788}\textbf{(0.1,1,0)} & \cellcolor[rgb]{ .788,  .788,  .788}\textbf{0.53} & \cellcolor[rgb]{ .788,  .788,  .788}\textbf{(0.1,1,0)} & \cellcolor[rgb]{ .788,  .788,  .788}\textbf{0.53} & \cellcolor[rgb]{ .788,  .788,  .788}\textbf{(0.1,1,0.3)} & \cellcolor[rgb]{ .788,  .788,  .788}\textbf{0.53} & \cellcolor[rgb]{ .851,  .882,  .949}\textbf{(10,0,0.01,0)} & \cellcolor[rgb]{ .851,  .882,  .949}\textbf{1.00} & \cellcolor[rgb]{ .851,  .882,  .949}\textbf{(10,0,0.01,0)} & \cellcolor[rgb]{ .851,  .882,  .949}\textbf{1.00} \\
          &  kC   & (0.01,1,0) & -12.19 & \cellcolor[rgb]{ .788,  .788,  .788}\textbf{(0.1,1,0)} & \cellcolor[rgb]{ .788,  .788,  .788}\textbf{0.53} & \cellcolor[rgb]{ .788,  .788,  .788}\textbf{(0.1,1,0.9)} & \cellcolor[rgb]{ .788,  .788,  .788}\textbf{0.53} & \cellcolor[rgb]{ .851,  .882,  .949}\textbf{(10,1,1,0.1)} & \cellcolor[rgb]{ .851,  .882,  .949}\textbf{1.00} & \cellcolor[rgb]{ .851,  .882,  .949}\textbf{(10,0,0.01,0.1)} & \cellcolor[rgb]{ .851,  .882,  .949}\textbf{1.00} \\
          &  MAX  & (0.01,0.01,0) & -0.63 & (0.01,1,0) & 0.03  & \cellcolor[rgb]{ .788,  .788,  .788}\textbf{(10,1,0.3)} & \cellcolor[rgb]{ .788,  .788,  .788}\textbf{0.53} & (10,1,0.1,0) & 0.11  & \cellcolor[rgb]{ .851,  .882,  .949}\textbf{(100,1,1,0)} & \cellcolor[rgb]{ .851,  .882,  .949}\textbf{1.00} \\
          &  akC  & (0.01,0.1,0) & -6.85 & (1,0.01,0.0156) & 0.28  & \cellcolor[rgb]{ .788,  .788,  .788}\textbf{(10,0.1,T)} & \cellcolor[rgb]{ .788,  .788,  .788}\textbf{0.42} & \cellcolor[rgb]{ .851,  .882,  .949}\textbf{(1,0,0.01,0)} & \cellcolor[rgb]{ .851,  .882,  .949}\textbf{1.00} & (100,0,1,T) & 0.99 \\
          &  MEDIAN & (0.1,0.1,0) & 0.20  & (0.01,0.1,0) & 0.20  & \cellcolor[rgb]{ .788,  .788,  .788}\textbf{(10,0.1,T)} & \cellcolor[rgb]{ .788,  .788,  .788}\textbf{0.40} & (0.01,0.1,0.01,0) & 0.20  & \cellcolor[rgb]{ .851,  .882,  .949}\textbf{(100,0,1,T)} & \cellcolor[rgb]{ .851,  .882,  .949}\textbf{0.62} \\
          & TM    & (0.01,0.1,0) & -0.05 & (1,0.1,0) & 0.07  & \cellcolor[rgb]{ .788,  .788,  .788}\textbf{(0.01,0.1,T)} & \cellcolor[rgb]{ .788,  .788,  .788}\textbf{0.11} & (1,0.1,0.01,0) & 0.88  & \cellcolor[rgb]{ .851,  .882,  .949}\textbf{(100,0,1,T)} & \cellcolor[rgb]{ .851,  .882,  .949}\textbf{0.98} \\
    \hline
    \end{tabular}%
    }
    \end{adjustwidth}
  \label{tab:gof2}%
\end{table}%

\subsection{Comparison of classic performance measures}
\label{performance}


To carry out this analysis, each dataset was divided into two parts: training set (75\% of the individuals) and test set (25\%). For Boston, CYGOB1, Eilon, DW2, DW3, Alcohol and HBK datasets, a 5-fold cross-validation (5-FCV) was performed in each training set to determine the most appropiate parameter values in terms of different performance metrics. In the training sample we stablished a time limit of 900 seconds and for test set we allowed a time limit of 1800 seconds. In the case of Quandt dataset, a 3-fold-cross validation (3-FCV) was performed due to the small number of individuals associated to this dataset. Then, the best parameter combination is applied to the training set and evaluated in the test set.  The implemented code can be found in \url{https://github.com/orel-uca/OWASVR}.

To evaluate the performance of the regression models we use two of the most common metrics in this field, namely MAE (mean absolute error) and MSE (mean squared error):

\begin{eqnarray*}
MAE &=&\frac{\sum_{i=1}^n |y_i-\widehat{y}_i|}{n},\\
MSE &=&\frac{\sum_{i=1}^n (y_i-\widehat{y}_i)^2}{n}.\\
\end{eqnarray*}

Note that, for these measures, $\widehat{y}_i$ denotes the predicted value of individual $i$ following the corresponding regression model. The first measure is the average of the absolute deviations and the second one is the average of the quadratic deviations obtained when the model is applied.

Table \ref{tab:3MAE} reports the best MAE results of each approach obtained in the test set after determining the best hyperparameter combination given by the fold cross validation. Observe that Table \ref{tab:3MAE} follows a similar structure to Table \ref{tab:gof2}. Note that, in   Table \ref{tab:3MAE}, each reported time corresponds 
 to the solving time required by the considered approach when applied to $75\%$ of the sample using the best hyperparameter values,  i.e., the time to compute the regression hyperplane to be evaluated by the test sample. 

Regarding the results of the first three column blocks of Table \ref{tab:3MAE}, we can observe that the exact approach (OWA-SVR block) provides the best MAE values in 35 out of the 48 studied cases. In 16 of these 35 results there are ties with the heuristic results. If we focus on the first heuristic (HEUR1 block), it only obtains the best results in 15 cases. Finally, the second heuristic (HEUR2 block) presents the best MAE results in 25 cases. Consequently, we can conclude that, in general, it is worthwhile applying the exact formulation to solve this problem. On the other hand, we can see that the required computing times in all cases are smaller when we use the heuristic approaches. 

If we focus on the last two column blocks, we observe that the use of Gaussian Kernel improves the MAE results in some datasets.  Specifically, we highlight in blue the cases when the performance is better than using the linear Kernel. Note that OWA-SVR$_K$ provides the best results for almost all considered weights in the  CYGOB1, Eilon and DW2 datasets.  We could derive that, for these datasets, it is advisable to use a non linear regression function.

Table \ref{tab:4MSE} follows the same structure as Table \ref{tab:3MAE}, but, in this case, the best results related to MSE are reported.  Note that, regarding the first three columns, the exact approach (OWA-SVR) shows the best average measure values in most of the cases. Particularly, in 32 out of the 48 studied cases, the exact approach provides the best results. Observe also that, HEUR1 provides the best MAE results in 21 cases and HEUR2 in 19 cases.

The results of two particular datasets deserve a special attention. Note that in the case of Quandt dataset for weights kc, MAX, akC, MED and TM, the average MSE values of the exact approach OWA-SVR and the heuristic HEUR2 are much better than the MSE provided by HEUR1. This implies that the use of the exact approach or HEUR2 is more adequate than using HEUR1 for this dataset.

A big contrast with respect to HEUR1 can also be observed for the Alcohol dataset. If we focus on weights kC, MAX and akC, we can see that the exact approach and HEUR2 perform much better for this measure.

If we focus on the results considering the Gaussian Kernel (HEUR$_K$ and OWA-SVR$_K$ blocks), we can observe that both approaches provide better MSE results for some datasets. In particular, HEUR$_K$ seems to be the best approach in terms of MSE for Boston datasets. In addition, MSE values of GYGOB1, Eilon and DW2 are improved in general when using the exact approach OWA-SVR$_K$.

Regarding computing times, Table \ref{tab:4MSE} shows that heuristic approaches are much faster than the exact approach. This difference is even bigger when we consider non-convex OWA operators.

 In Tables \ref{tab:3MAE} and \ref{tab:4MSE}, we can also compare the classic SVR (SUM $\lambda$-weight) with other weights. In general, we can observe that the classic SUM weight is not always the one providing best MAE and MSE performances. For instance, if we focus on Eilon dataset and OWA-SVR$_K$ approach, MAX, akC and  MED show better MAE than SUM. For the same dataset, in Table \ref{tab:4MSE}, we can see that MAX, akC and  MED provide the same MSE as SUM weight. In the case of Alchohol dataset and OWA-SVR approach, the MAE associated with classic SVM is improved by MAX and MED weights. Similarly, kC, MAX, akC and MED allow to obtain a better MSE than SUM weight.

\begin{table}[htbp]
  \centering
  \caption{Mean Absolute Error (MAE) for the different datasets and solution methods}
    \begin{adjustwidth}{-3.5cm}{0cm}
  \scalebox{0.75}{
    \begin{tabular}{|c|l|lr|lr|lr|lr|lr|}
        \hline
    \multirow{2}[2]{*}{Data} & \multicolumn{1}{c|}{\multirow{2}[2]{*}{Weight}} & \multicolumn{2}{c|}{HEUR1} & \multicolumn{2}{c|}{HEUR2} & \multicolumn{2}{c|}{OWA-SVR} & \multicolumn{2}{c|}{HEUR$_K$} & \multicolumn{2}{c|}{OWA-SVR$_K$} \\
          &       & $(C,\varepsilon,t)$ & \multicolumn{1}{l|}{MAE} & $(C,\varepsilon,t)$ & \multicolumn{1}{l|}{MAE} & $(C,\varepsilon,t)$ & \multicolumn{1}{l|}{MAE} & $(C,\varepsilon,\sigma, t)$ & \multicolumn{1}{l|}{MAE} & $(C,\varepsilon,\sigma,t)$& \multicolumn{1}{l|}{MAE} \\
        \hline
    \multirow{6}[2]{*}{\begin{sideways}Quandt\end{sideways}} & SUM   & \cellcolor[rgb]{  .859,  .859,  .859}\textbf{(0.1,1,0)} & \cellcolor[rgb]{  .859,  .859,  .859}\textbf{0.75} & \cellcolor[rgb]{  .859,  .859,  .859}\textbf{(0.1,1,0)} & \cellcolor[rgb]{  .859,  .859,  .859}\textbf{0.75} & \cellcolor[rgb]{  .859,  .859,  .859}\textbf{(0.1,1,0)} & \cellcolor[rgb]{  .859,  .859,  .859}\textbf{0.75} & (10,1,10,0) & 1.04  & (100,1,10,0) & 1.42 \\
          & kC    & (0.01,0.1,0) & 3.51  & \cellcolor[rgb]{  .859,  .859,  .859}\textbf{(0.1,1,0)} & \cellcolor[rgb]{  .859,  .859,  .859}\textbf{0.75} & \cellcolor[rgb]{  .859,  .859,  .859}\textbf{(0.1,1,0)} & \cellcolor[rgb]{  .859,  .859,  .859}\textbf{0.75} & (100,1,10,0) & 1.32  & (100,1,10,0) & 1.42 \\
          & MAX   & (0.01,0.1,0) & 2.40  & (1,0.1,0) & 3.36  & \cellcolor[rgb]{  .859,  .859,  .859}\textbf{(10,0.1,0)} & \cellcolor[rgb]{  .859,  .859,  .859}\textbf{0.73} & (100,0.1,10,0) & 3.84  & (100,0.01,1,0) & 4.01 \\
          & akC   & (0.01,0.1,0) & 6.14  & \cellcolor[rgb]{  .859,  .859,  .859}\textbf{(100,1,0)} & \cellcolor[rgb]{  .859,  .859,  .859}\textbf{0.79} & (100,0,0.7) & 1.59  & (100,0.01,1,0) & 2.75  & (100,0.01,10,0.6) & 1.65 \\
          & MED   & (0.1,0.1,0) & 10.93 & (0.1,0.1,0) & 5.04  & \cellcolor[rgb]{  .859,  .859,  .859}\textbf{(10,0.1,0.5)} & \cellcolor[rgb]{  .859,  .859,  .859}\textbf{1.02} & (100,0,10,0) & 7.05  & (100,0.01,10,0.1) & 2.90 \\
          & TM    & (0.01,0.1,0) & 3.71  & (10,0.1,0) & 1.05  & \cellcolor[rgb]{  .859,  .859,  .859}\textbf{(1,0.1,3.8)} & \cellcolor[rgb]{  .859,  .859,  .859}\textbf{0.90} & \cellcolor[rgb]{ .851,  .882,  .949}\textbf{(100,1,10,0)} & \cellcolor[rgb]{ .851,  .882,  .949}\textbf{0.82} & (10,0.01,10,1.2) & 3.41 \\
        \hline
    \multirow{6}[2]{*}{\begin{sideways}Boston\end{sideways}} & SUM   & \cellcolor[rgb]{  .859,  .859,  .859}\textbf{(100,0.1,0)} & \cellcolor[rgb]{  .859,  .859,  .859}\textbf{0.16} & \cellcolor[rgb]{  .859,  .859,  .859}\textbf{(100,0.1,0)} & \cellcolor[rgb]{  .859,  .859,  .859}\textbf{0.16} & \cellcolor[rgb]{  .859,  .859,  .859}\textbf{(100,0.1,0)} & \cellcolor[rgb]{  .859,  .859,  .859}\textbf{0.16} & \cellcolor[rgb]{ .851,  .882,  .949}\textbf{(100,0.1,10,0)} & \cellcolor[rgb]{ .851,  .882,  .949}\textbf{0.16} & \cellcolor[rgb]{ .851,  .882,  .949}\textbf{(1,0.01,1,0)} & \cellcolor[rgb]{ .851,  .882,  .949}\textbf{0.16} \\
          & kC    & (0.1,0.1,0) & 0.20  & \cellcolor[rgb]{ .859,  .859,  .859}\textbf{(100,0.1,0)} & \cellcolor[rgb]{ .859,  .859,  .859}\textbf{0.16} & \cellcolor[rgb]{ .859,  .859,  .859}\textbf{(100,0.1,0)} & \cellcolor[rgb]{ .859,  .859,  .859}\textbf{0.16} & \cellcolor[rgb]{ .851,  .882,  .949}\textbf{(100,0.1,10,0)} & \cellcolor[rgb]{ .851,  .882,  .949}\textbf{0.16} & (100,0.1,10,0) & 0.18 \\
          & MAX   & \cellcolor[rgb]{  .859,  .859,  .859}\textbf{(1,0.1,0)} & \cellcolor[rgb]{  .859,  .859,  .859}\textbf{0.20} & (0.01,0,0) & 0.63  & (100,0,0) & 0.25  & \cellcolor[rgb]{ .851,  .882,  .949}\textbf{(1,0.01,0.01,0)} & \cellcolor[rgb]{ .851,  .882,  .949}\textbf{0.17} & (1,0.1,0.01,0) & 0.20 \\
          & akC   & (0.01,0.1,0) & 0.26  & \cellcolor[rgb]{  .859,  .859,  .859}\textbf{(1,0,0)} & \cellcolor[rgb]{  .859,  .859,  .859}\textbf{0.19} & (100,0,144.6) & 0.22  & \cellcolor[rgb]{ .851,  .882,  .949}\textbf{(1,0.01,1,0)} & \cellcolor[rgb]{ .851,  .882,  .949}\textbf{0.18} & (10,0,1,114.6) & 0.22 \\
          & MED   & (10,0.1,0) & 0.23  & (0.1,0.1,0) & 0.23  & \cellcolor[rgb]{  .859,  .859,  .859}\textbf{(10,0,7.2)} & \cellcolor[rgb]{  .859,  .859,  .859}\textbf{0.22} & \cellcolor[rgb]{ .851,  .882,  .949}\textbf{(0.1,0.01,0.1,0)} & \cellcolor[rgb]{ .851,  .882,  .949}\textbf{0.18} & (100,0,1,0.5) & 0.22 \\
          & TM    & (10,0.1,0) & 0.23  & \cellcolor[rgb]{ .859,  .859,  .859}\textbf{(100,0.1,0)} & \cellcolor[rgb]{ .859,  .859,  .859}\textbf{0.16} & (10,0.01,1330.4) & 0.20  & \cellcolor[rgb]{ .851,  .882,  .949}\textbf{(10,0.1,1,0)} & \cellcolor[rgb]{ .851,  .882,  .949}\textbf{0.16} & (1,0.01,1,160) & 0.20 \\
        \hline
    \multirow{6}[2]{*}{\begin{sideways}CYGOB1\end{sideways}} & SUM   & (0.01,0.01,0) & 0.46  & \cellcolor[rgb]{  .859,  .859,  .859}\textbf{(0.01,1,0)} & \cellcolor[rgb]{  .859,  .859,  .859}\textbf{0.45} & \cellcolor[rgb]{  .859,  .859,  .859}\textbf{(0.01,1,0.1)} & \cellcolor[rgb]{  .859,  .859,  .859}\textbf{0.45} & \cellcolor[rgb]{ .851,  .882,  .949}\textbf{(1,0.1,0.01,0)} & \cellcolor[rgb]{ .851,  .882,  .949}\textbf{0.24} & \cellcolor[rgb]{ .851,  .882,  .949}\textbf{(1,0.1,0.01,0)} & \cellcolor[rgb]{ .851,  .882,  .949}\textbf{0.24} \\
          & kC    & (0.01,0.1,0) & 0.47  & \cellcolor[rgb]{  .859,  .859,  .859}\textbf{(0.01,1,0)} & \cellcolor[rgb]{  .859,  .859,  .859}\textbf{0.45} & \cellcolor[rgb]{ .859,  .859,  .859}\textbf{(10,10,0.1)} & \cellcolor[rgb]{ .859,  .859,  .859}\textbf{0.45} & \cellcolor[rgb]{ .851,  .882,  .949}\textbf{(1,0.1,0.01,0)} & \cellcolor[rgb]{ .851,  .882,  .949}\textbf{0.22} & (1,0.1,0.01,0) & 0.24 \\
          & MAX   & (0.01,0.01,0) & 0.46  & (0.1,0.1,0) & 0.98  & \cellcolor[rgb]{ .859,  .859,  .859}\textbf{(0.01,1,0.1)} & \cellcolor[rgb]{ .859,  .859,  .859}\textbf{0.45} & (10,0,0.01,0) & 0.70  & \cellcolor[rgb]{ .851,  .882,  .949}\textbf{(10,0.1,0.01,0)} & \cellcolor[rgb]{ .851,  .882,  .949}\textbf{0.26} \\
          & akC   & (0.01,0,0) & 0.46  & \cellcolor[rgb]{ .859,  .859,  .859}\textbf{(10,10,0)} & \cellcolor[rgb]{ .859,  .859,  .859}\textbf{0.45} & \cellcolor[rgb]{ .859,  .859,  .859}\textbf{(10,10,0.1)} & \cellcolor[rgb]{ .859,  .859,  .859}\textbf{0.45} & \cellcolor[rgb]{ .851,  .882,  .949}\textbf{(100,0.1,1,0)} & \cellcolor[rgb]{ .851,  .882,  .949}\textbf{0.26} & (10,0,0.1,T) & 0.46 \\
          & MED   & \cellcolor[rgb]{  .859,  .859,  .859}\textbf{(0.1,0.1,0)} & \cellcolor[rgb]{  .859,  .859,  .859}\textbf{0.45} & (10,0,0) & 0.70  & \cellcolor[rgb]{ .859,  .859,  .859}\textbf{(0.1,1,0.1)} & \cellcolor[rgb]{ .859,  .859,  .859}\textbf{0.45} & (1,0,0.01,0) & 0.48  & \cellcolor[rgb]{ .851,  .882,  .949}\textbf{(100,0,0.1,26.3)} & \cellcolor[rgb]{ .851,  .882,  .949}\textbf{0.43} \\
          & TM    & \cellcolor[rgb]{  .859,  .859,  .859}\textbf{(0.01,0.01,0)} & \cellcolor[rgb]{  .859,  .859,  .859}\textbf{0.45} & \cellcolor[rgb]{ .859,  .859,  .859}\textbf{(0.01,1,0)} & \cellcolor[rgb]{ .859,  .859,  .859}\textbf{0.45} & \cellcolor[rgb]{ .859,  .859,  .859}\textbf{(0.01,1,0.1)} & \cellcolor[rgb]{ .859,  .859,  .859}\textbf{0.45} & (100,0.01,0.01,0) & 0.36  & \cellcolor[rgb]{ .851,  .882,  .949}\textbf{(10,0,0.01,T)} & \cellcolor[rgb]{ .851,  .882,  .949}\textbf{0.24} \\
       \hline
    \multirow{6}[2]{*}{\begin{sideways}Eilon\end{sideways}} & SUM   & \cellcolor[rgb]{  .859,  .859,  .859}\textbf{(0.01,0.01,0)} & \cellcolor[rgb]{  .859,  .859,  .859}\textbf{0.22} & \cellcolor[rgb]{  .859,  .859,  .859}\textbf{(1,1,0)} & \cellcolor[rgb]{  .859,  .859,  .859}\textbf{0.22} & \cellcolor[rgb]{  .859,  .859,  .859}\textbf{(10,1,0.1)} & \cellcolor[rgb]{  .859,  .859,  .859}\textbf{0.22} & \cellcolor[rgb]{ .851,  .882,  .949}\textbf{(10,1,0.01,0)} & \cellcolor[rgb]{ .851,  .882,  .949}\textbf{0.22} & \cellcolor[rgb]{ .851,  .882,  .949}\textbf{(0.01,0,100,0)} & \cellcolor[rgb]{ .851,  .882,  .949}\textbf{0.22} \\
          & kC    & \cellcolor[rgb]{  .859,  .859,  .859}\textbf{(0.01,0.01,0)} & \cellcolor[rgb]{  .859,  .859,  .859}\textbf{0.21} & (1,1,0) & 0.23  & (1,1,0.1) & 0.23  & \cellcolor[rgb]{ .851,  .882,  .949}\textbf{(100,1,0.01,0)} & \cellcolor[rgb]{ .851,  .882,  .949}\textbf{0.22} & (1,0.1,0.1,T) & 0.30 \\
          & MAX   & \cellcolor[rgb]{  .859,  .859,  .859}\textbf{(0.01,0.01,0)} & \cellcolor[rgb]{  .859,  .859,  .859}\textbf{0.22} & (10,0,0) & 0.44  & \cellcolor[rgb]{ .859,  .859,  .859}\textbf{(1,0,0.1)} & \cellcolor[rgb]{ .859,  .859,  .859}\textbf{0.22} & (100,0.01,0.01,0) & 0.21  & \cellcolor[rgb]{ .851,  .882,  .949}\textbf{(100,0.1,0.01,0)} & \cellcolor[rgb]{ .851,  .882,  .949}\textbf{0.20} \\
          & akC   & (0.01,0,0) & 0.23  & \cellcolor[rgb]{ .859,  .859,  .859}\textbf{(100,10,0)} & \cellcolor[rgb]{ .859,  .859,  .859}\textbf{0.21} & \cellcolor[rgb]{ .859,  .859,  .859}\textbf{(100,10,0.1)} & \cellcolor[rgb]{ .859,  .859,  .859}\textbf{0.21} & (100,1,0.1,T) & 0.82  & \cellcolor[rgb]{ .851,  .882,  .949}\textbf{(10,1,0.01,0.1)} & \cellcolor[rgb]{ .851,  .882,  .949}\textbf{0.21} \\
          & MED   & \cellcolor[rgb]{  .859,  .859,  .859}\textbf{(0.01,0.01,0)} & \cellcolor[rgb]{  .859,  .859,  .859}\textbf{0.22} & (0.1,0.1,0) & 0.37  & \cellcolor[rgb]{ .859,  .859,  .859}\textbf{(0.01,0.1,12.9)} & \cellcolor[rgb]{ .859,  .859,  .859}\textbf{0.22} & (1,0.1,0.1,T) & 0.37  & \cellcolor[rgb]{ .851,  .882,  .949}\textbf{(1,0.1,0.01,2.6)} & \cellcolor[rgb]{ .851,  .882,  .949}\textbf{0.20} \\
          & TM    & (0.1,0,0) & 0.23  & \cellcolor[rgb]{ .859,  .859,  .859}\textbf{(10,1,0)} & \cellcolor[rgb]{ .859,  .859,  .859}\textbf{0.22} & \cellcolor[rgb]{ .859,  .859,  .859}\textbf{(100,1,0.1)} & \cellcolor[rgb]{ .859,  .859,  .859}\textbf{0.22} & (10,10,0.1,T) & 9.55  & (100,10,100,0.2) & 0.56 \\
        \hline
    \multirow{6}[2]{*}{\begin{sideways}DW2\end{sideways}} & SUM   & \cellcolor[rgb]{  .859,  .859,  .859}\textbf{(10,0.01,0)} & \cellcolor[rgb]{  .859,  .859,  .859}\textbf{0.28} & \cellcolor[rgb]{  .859,  .859,  .859}\textbf{(10,0.01,0)} & \cellcolor[rgb]{  .859,  .859,  .859}\textbf{0.28} & \cellcolor[rgb]{  .859,  .859,  .859}\textbf{(10,0.01,0.3)} & \cellcolor[rgb]{  .859,  .859,  .859}\textbf{0.28} & (100,0,0.1,T) & 0.30  & \cellcolor[rgb]{ .851,  .882,  .949}\textbf{(100,0,0.1,T)} & \cellcolor[rgb]{ .851,  .882,  .949}\textbf{0.09} \\
          & kC    & (0.01,0.01,0) & 0.31  & (10,0.1,0) & 0.29  & \cellcolor[rgb]{ .859,  .859,  .859}\textbf{(10,0,0.4)} & \cellcolor[rgb]{ .859,  .859,  .859}\textbf{0.28} & \cellcolor[rgb]{ .851,  .882,  .949}(1,0.1,0.1,T) & \cellcolor[rgb]{ .851,  .882,  .949}0.10 & \cellcolor[rgb]{ .851,  .882,  .949}\textbf{(10,0.01,0.1,T)} & \cellcolor[rgb]{ .851,  .882,  .949}\textbf{0.10} \\
          & MAX   & (0.01,0.01,0) & 0.29  & (0.01,1,0) & 0.60  & \cellcolor[rgb]{ .859,  .859,  .859}\textbf{(100,0.01,0.3)} & \cellcolor[rgb]{ .859,  .859,  .859}\textbf{0.28} & (100,0.1,0.01,0) & 0.30  & \cellcolor[rgb]{ .851,  .882,  .949}\textbf{(100,0.1,0.1,T)} & \cellcolor[rgb]{ .851,  .882,  .949}\textbf{0.14} \\
          & akC   & (0.01,0,0) & 0.30  & \cellcolor[rgb]{ .859,  .859,  .859}\textbf{(10,0,0)} & \cellcolor[rgb]{ .859,  .859,  .859}\textbf{0.28} & (100,0.01,T) & 0.31  & (10,0.01,0.1,T) & 0.10  & \cellcolor[rgb]{ .851,  .882,  .949}\textbf{(100,0,0.1,T)} & \cellcolor[rgb]{ .851,  .882,  .949}\textbf{0.09} \\
          & MED   & \cellcolor[rgb]{  .859,  .859,  .859}\textbf{(0.1,0.01,0)} & \cellcolor[rgb]{  .859,  .859,  .859}\textbf{0.29} & \cellcolor[rgb]{ .859,  .859,  .859}\textbf{(10,0,0)} & \cellcolor[rgb]{ .859,  .859,  .859}\textbf{0.29} & (100,0.01,126.3) & 0.32  & (10,0,1,0) & 0.31  & (100,0.01,0.1,T) & 0.34 \\
          & TM    & (0.1,0,0) & 0.30  & \cellcolor[rgb]{ .859,  .859,  .859}\textbf{(100,0.1,0)} & \cellcolor[rgb]{ .859,  .859,  .859}\textbf{0.28} & (100,0.1,244.3) & 0.29  & (10,0,0.1,T) & 0.11  & \cellcolor[rgb]{ .851,  .882,  .949}\textbf{(10,0,0.1,T)} & \cellcolor[rgb]{ .851,  .882,  .949}\textbf{0.10} \\
        \hline
    \multirow{6}[2]{*}{\begin{sideways}alcohol\end{sideways}} & SUM   & \cellcolor[rgb]{  .859,  .859,  .859}\textbf{(10,0,0)} & \cellcolor[rgb]{  .859,  .859,  .859}\textbf{0.49} & \cellcolor[rgb]{  .859,  .859,  .859}\textbf{(10,0,0)} & \cellcolor[rgb]{  .859,  .859,  .859}\textbf{0.49} & \cellcolor[rgb]{  .859,  .859,  .859}\textbf{(10,0,0.1)} & \cellcolor[rgb]{  .859,  .859,  .859}\textbf{0.49} & (100,0,100,0) & 0.51  & (100,0.01,100,0) & 0.53 \\
          & kC    & (100,10,0) & 11.42 & \cellcolor[rgb]{  .859,  .859,  .859}\textbf{(100,0.01,0)} & \cellcolor[rgb]{  .859,  .859,  .859}\textbf{0.48} & (10,0.01,0.1) & 0.49  & (10,1,100,0) & 0.82  & (100,0.1,100,0) & 0.72 \\
          & MAX   & (100,10,0) & 11.42 & (1,0,0) & 3.07  & \cellcolor[rgb]{  .859,  .859,  .859}\textbf{(1,0.01,0.1)} & \cellcolor[rgb]{  .859,  .859,  .859}\textbf{0.39} & (100,0,100,0) & 3.04  & (100,1,100,0) & 0.82 \\
          & akC   & (10,10,0) & 11.42 & (10,0.01,0) & 0.61  & \cellcolor[rgb]{  .859,  .859,  .859}\textbf{(100,0,T)} & \cellcolor[rgb]{  .859,  .859,  .859}\textbf{0.59} & \cellcolor[rgb]{ .851,  .882,  .949}\textbf{(100,0.1,100,0)} & \cellcolor[rgb]{ .851,  .882,  .949}\textbf{0.59} & (100,0.01,100,T) & 1.73 \\
          & MED   & (0.01,1,0) & 3.17  & (0.01,0.01,0) & 3.00  & \cellcolor[rgb]{  .859,  .859,  .859}\textbf{(0.1,0.01,753.9)} & \cellcolor[rgb]{  .859,  .859,  .859}\textbf{0.45} & (1,0.1,100,0) & 3.06  & (100,0.01,100,5.4) & 1.69 \\
          & TM    & (0.01,1,0) & 3.17  & (100,0.1,0) & 0.66  & \cellcolor[rgb]{  .859,  .859,  .859}\textbf{(100,0.01,T)} & \cellcolor[rgb]{  .859,  .859,  .859}\textbf{0.61} & (10,0.1,100,0) & 0.65  & (100,0.1,100,T) & 0.67 \\
        \hline
    \multirow{6}[2]{*}{\begin{sideways}DW3\end{sideways}} & SUM   & \cellcolor[rgb]{  .859,  .859,  .859}\textbf{(1,0.01,0)} & \cellcolor[rgb]{  .859,  .859,  .859}\textbf{0.06} & \cellcolor[rgb]{  .859,  .859,  .859}\textbf{(1,0.01,0)} & \cellcolor[rgb]{  .859,  .859,  .859}\textbf{0.06} & \cellcolor[rgb]{  .859,  .859,  .859}\textbf{(1,0.01,0.4)} & \cellcolor[rgb]{  .859,  .859,  .859}\textbf{0.06} & \cellcolor[rgb]{ .851,  .882,  .949}\textbf{(1,0,0.1,0)} &\cellcolor[rgb]{ .851,  .882,  .949} \textbf{0.04} & (10,0.01,0.1,0) & 0.05 \\
          & kC    & (0.01,0.01,0) & 0.41  & (10,0.1,0) & 0.07  & \cellcolor[rgb]{ .859,  .859,  .859}\textbf{(10,0,0.4)} & \cellcolor[rgb]{ .859,  .859,  .859}\textbf{0.05} & \cellcolor[rgb]{ .851,  .882,  .949}\textbf{(10,0.01,0.1,0.1)} & \cellcolor[rgb]{ .851,  .882,  .949}\textbf{0.05} & \cellcolor[rgb]{ .851,  .882,  .949}\textbf{(1,0.01,0.1,0)} & \cellcolor[rgb]{ .851,  .882,  .949}\textbf{0.05} \\
          & MAX   & (0.01,0.01,0) & 0.30  & (1,0,0) & 0.34  & \cellcolor[rgb]{ .859,  .859,  .859}\textbf{(10,0.1,0.3)} & \cellcolor[rgb]{ .859,  .859,  .859}\textbf{0.08} & (100,0,10,0) & 0.34  & (10,0.1,1,0) & 0.09 \\
          & akC   & (0.01,0,0) & 0.19  & \cellcolor[rgb]{  .859,  .859,  .859}\textbf{(100,0,0)} & \cellcolor[rgb]{  .859,  .859,  .859}\textbf{0.06} & \cellcolor[rgb]{ .859,  .859,  .859}\textbf{(100,0,T)} & \cellcolor[rgb]{ .859,  .859,  .859}\textbf{0.06} & \cellcolor[rgb]{ .851,  .882,  .949}\textbf{(10,0,0.1,0)} & \cellcolor[rgb]{ .851,  .882,  .949}\textbf{0.05} & \cellcolor[rgb]{ .851,  .882,  .949}\textbf{(100,0,1,T)} & \cellcolor[rgb]{ .851,  .882,  .949}\textbf{0.05} \\
          & MED   & (0.1,0,0) & 0.31  & (0.01,0,0) & 0.34  & \cellcolor[rgb]{ .859,  .859,  .859}\textbf{(100,0,T)} & \cellcolor[rgb]{ .859,  .859,  .859}\textbf{0.06} & (0.1,0,1,0) & 0.32  & \cellcolor[rgb]{ .851,  .882,  .949}\textbf{(100,0,0.1,T)} & \cellcolor[rgb]{ .851,  .882,  .949}\textbf{0.05} \\
          & TM    & (0.1,0,0) & 0.23  & (10,0,0) & 0.06  & \cellcolor[rgb]{ .859,  .859,  .859}\textbf{(10,0,T)} & \cellcolor[rgb]{ .859,  .859,  .859}\textbf{0.05} & \cellcolor[rgb]{ .851,  .882,  .949}\textbf{(100,0.01,0.1,0.1)} & \cellcolor[rgb]{ .851,  .882,  .949}\textbf{0.05} & (10,0.01,0.1,T) & 0.06 \\
        \hline
    \multirow{6}[2]{*}{\begin{sideways}HBK\end{sideways}} & SUM   & (0.1,1,0) & 0.74  & \cellcolor[rgb]{  .859,  .859,  .859}\textbf{(0.01,1,0)} & \cellcolor[rgb]{  .859,  .859,  .859}\textbf{0.65} & \cellcolor[rgb]{  .859,  .859,  .859}\textbf{(0.01,1,0.5)} & \cellcolor[rgb]{  .859,  .859,  .859}\textbf{0.65} & \cellcolor[rgb]{ .851,  .882,  .949}\textbf{(100,1,10,0.1)} & \cellcolor[rgb]{ .851,  .882,  .949}\textbf{0.61} & \cellcolor[rgb]{ .851,  .882,  .949}\textbf{(100,1,10,0)} & \cellcolor[rgb]{ .851,  .882,  .949}\textbf{0.61} \\
          & kC    & (10,10,0) & 5.81  & (0.01,0.1,0) & 1.37  & \cellcolor[rgb]{ .859,  .859,  .859}\textbf{(0.01,1,0.6)} & \cellcolor[rgb]{ .859,  .859,  .859}\textbf{0.65} & \cellcolor[rgb]{ .851,  .882,  .949}\textbf{(100,1,10,0.1)} & \cellcolor[rgb]{ .851,  .882,  .949}\textbf{0.61} & \cellcolor[rgb]{ .851,  .882,  .949}\textbf{(100,1,10,0)} & \cellcolor[rgb]{ .851,  .882,  .949}\textbf{0.61} \\
          & MAX   & \cellcolor[rgb]{ .859,  .859,  .859}\textbf{(0.01,0.01,0)} & \cellcolor[rgb]{ .859,  .859,  .859}\textbf{1.05} & (0.01,0,0) & 1.46  & (0.01,10,0.5) & 1.96  & (100,0.01,1,0) & 1.86  & \cellcolor[rgb]{ .851,  .882,  .949}\textbf{(100,1,10,0)} & \cellcolor[rgb]{ .851,  .882,  .949}\textbf{0.69} \\
          & akC   & (100,10,0) & 5.81  & \cellcolor[rgb]{ .859,  .859,  .859}\textbf{(0.1,0.1,0)} & \cellcolor[rgb]{ .859,  .859,  .859}\textbf{0.66} & (1,1,2.2) & 1.48  & (100,0.01,10,0.1) & 0.70  & (100,1,0.01,4) & 1.61 \\
          & MED   & \cellcolor[rgb]{ .859,  .859,  .859}\textbf{(0.01,0.1,0)} & \cellcolor[rgb]{ .859,  .859,  .859}\textbf{1.51} & (0.1,0,0) & 1.57  & (10,1,0.7) & 1.92  & (10,0,10,0) & 1.85  & (10,1,0.1,0.2) & 1.85 \\
          & TM    & (1,1,0) & 8.62  & \cellcolor[rgb]{ .859,  .859,  .859}\textbf{(1,1,0)} & \cellcolor[rgb]{ .859,  .859,  .859}\textbf{0.65} & (0.01,1,4.4) & 1.51  & (10,0.01,1,0) & 1.17  & (100,0.01,10,T) & 0.70 \\
      \hline
    \end{tabular}%
        }
       \end{adjustwidth}
  \label{tab:3MAE}%
\end{table}%

\begin{table}[htbp]
  \centering
  \caption{Mean Squared Error (MSE) for the different datasets and solution methods}
      \begin{adjustwidth}{-3.5cm}{0cm}
  \scalebox{0.75}{
    \begin{tabular}{|c|l|lr|lr|lr|lr|lr|}
        \hline
    \multirow{2}[2]{*}{Data} & \multicolumn{1}{c|}{\multirow{2}[2]{*}{Weight}} & \multicolumn{2}{c|}{HEUR1} & \multicolumn{2}{c|}{HEUR2} & \multicolumn{2}{c|}{OWA-SVR} & \multicolumn{2}{c|}{HEUR$_K$} & \multicolumn{2}{c|}{OWA-SVR$_K$} \\
          &       & $(C,\varepsilon,t)$ & \multicolumn{1}{l|}{MSE} & $(C,\varepsilon,t)$ & \multicolumn{1}{l|}{MSE} & $(C,\varepsilon,t)$ & \multicolumn{1}{l|}{MSE} & $(C,\varepsilon,\sigma,t)$ & \multicolumn{1}{l|}{MSE} & $(C,\varepsilon,\sigma,t)$ & \multicolumn{1}{l|}{MSE} \\
        \hline
    \multirow{6}[2]{*}{\begin{sideways}Quandt\end{sideways}} & SUM   & \cellcolor[rgb]{  .859,  .859,  .859}\textbf{(0.1,1,0)} & \cellcolor[rgb]{  .859,  .859,  .859}\textbf{0.72} & \cellcolor[rgb]{  .859,  .859,  .859}\textbf{(0.1,1,0)} & \cellcolor[rgb]{ .851,  .851,  .851}\textbf{0.72} & \cellcolor[rgb]{ .851,  .851,  .851}\textbf{(0.1,1,0)} & \cellcolor[rgb]{ .851,  .851,  .851}\textbf{0.72} & (10,1,10,0) & 2.28  & (100,1,10,0) & 2.18 \\
          & kC    & (0.01,0.1,0) & 14.78 & (0.1,0.1,0) & 3.28  & \cellcolor[rgb]{ .851,  .851,  .851}\textbf{(0.1,0.01,0)} & \cellcolor[rgb]{ .851,  .851,  .851}\textbf{0.70} & (100,1,10,0) & 2.69  & (100,1,10,0) & 2.18 \\
          & MAX   & (0.01,0.1,0) & 6.96  & (1,0.1,0) & 15.34 & \cellcolor[rgb]{ .851,  .851,  .851}\textbf{(10,0.1,0)} & \cellcolor[rgb]{ .851,  .851,  .851}\textbf{0.58} & (100,0.1,10,0) & 17.72 & (100,0,1,0) & 23.23 \\
          & akC   & (0.01,0.1,0) & 51.71 & \cellcolor[rgb]{ .851,  .851,  .851}\textbf{(1,1,0)} & \cellcolor[rgb]{ .851,  .851,  .851}\textbf{0.72} & (100,0.1,0.6) & 1.54  & (10,0.1,10,0) & 5.84  & (100,0,10,0.8) & 4.29 \\
          & MED   & (0.1,0.1,0) & 156.17 & (0.1,0.1,0) & 40.19 & \cellcolor[rgb]{ .851,  .851,  .851}\textbf{(10,0.1,0.5)} & \cellcolor[rgb]{ .851,  .851,  .851}\textbf{1.02} & (100,0,10,0) & 64.61 & (100,0,10,0.1) & 14.77 \\
          & TM    & (0.01,0.1,0) & 18.75 & (100,0.1,0) & 2.02  & \cellcolor[rgb]{ .851,  .851,  .851}\textbf{(1,0.1,4.1)} & \cellcolor[rgb]{ .851,  .851,  .851}\textbf{0.90} & (100,1,10,0) & 1.15  & (100,0,10,4) & 23.85 \\
        \hline
    \multirow{6}[2]{*}{\begin{sideways}Boston\end{sideways}} & SUM   & \cellcolor[rgb]{ .851,  .851,  .851}\textbf{(100,0.1,0)} & \cellcolor[rgb]{ .851,  .851,  .851}\textbf{0.05} & \cellcolor[rgb]{ .851,  .851,  .851}\textbf{(100,0.1,0)} & \cellcolor[rgb]{ .851,  .851,  .851}\textbf{0.05} & \cellcolor[rgb]{ .851,  .851,  .851}\textbf{(100,0.1,0)} & \cellcolor[rgb]{ .851,  .851,  .851}\textbf{0.05} & (10,0.01,0.1,0) & 0.14  & (10,0.01,0.1,0) & 0.14 \\
          & kC    & (0.1,0,0) & 0.07  & \cellcolor[rgb]{ .851,  .851,  .851}\textbf{(100,0.1,0)} & \cellcolor[rgb]{ .851,  .851,  .851}\textbf{0.05} & \cellcolor[rgb]{ .851,  .851,  .851}\textbf{(100,0.1,0)} & \cellcolor[rgb]{ .851,  .851,  .851}\textbf{0.05} & (10,0.1,1,0) & 0.05  & (10,0.01,0.1,0) & 0.13 \\
          & MAX   & \cellcolor[rgb]{ .851,  .851,  .851}\textbf{(1,0.1,0)} & \cellcolor[rgb]{ .851,  .851,  .851}\textbf{0.06} & (0.01,0,0) & 0.43  & (100,0,0) & 0.08  & \cellcolor[rgb]{ .851,  .882,  .949}\textbf{(1,0.01,0.01,0)} & \cellcolor[rgb]{ .851,  .882,  .949}\textbf{0.04} & (1,0.1,0.01,0) & 0.06 \\
          & akC   & (0.01,0,0) & 0.07  & \cellcolor[rgb]{ .851,  .851,  .851}\textbf{(1,0,0)} & \cellcolor[rgb]{ .851,  .851,  .851}\textbf{0.06} & (100,0,145.7) & 0.08  & (10,0.01,0.1,0) & 0.12  & (100,0.01,1,38.5) & 0.07 \\
          & MED   & (10,0.1,0) & 0.08  & (0.1,0.1,0) & 0.08  & \cellcolor[rgb]{ .851,  .851,  .851}\textbf{(100,0.01,5.8)} & \cellcolor[rgb]{ .851,  .851,  .851}\textbf{0.07} & (0.01,0,0.1,0) & 0.11  & (10,0,0.1,3.5) & 0.08 \\
          & TM    & (10,0.1,0) & 0.08  & \cellcolor[rgb]{ .851,  .851,  .851}\textbf{(100,0.1,0)} & \cellcolor[rgb]{ .851,  .851,  .851}\textbf{0.05} & (10,0.01,1307.2) & 0.06  & (10,0,0.1,0) & 0.07  & (1,0.01,1,160.9) & 0.06 \\
        \hline
    \multirow{6}[2]{*}{\begin{sideways}CYGOB1\end{sideways}} & SUM   & \cellcolor[rgb]{ .851,  .851,  .851}\textbf{(0.01,0.1,0)} & \cellcolor[rgb]{ .851,  .851,  .851}\textbf{0.32} & \cellcolor[rgb]{ .851,  .851,  .851}\textbf{(0.01,1,0)} & \cellcolor[rgb]{ .851,  .851,  .851}\textbf{0.32} & \cellcolor[rgb]{ .851,  .851,  .851}\textbf{(0.01,1,0.1)} & \cellcolor[rgb]{ .851,  .851,  .851}\textbf{0.32} & \cellcolor[rgb]{ .851,  .882,  .949}\textbf{(100,0.1,1,0)} & \cellcolor[rgb]{ .851,  .882,  .949}\textbf{0.08} & (10,0.1,0.01,0) & 0.09 \\
          & kC    & (0.01,0.1,0) & 0.46  & \cellcolor[rgb]{ .851,  .851,  .851}\textbf{(0.01,1,0)} & \cellcolor[rgb]{ .851,  .851,  .851}\textbf{0.32} & \cellcolor[rgb]{ .851,  .851,  .851}\textbf{(0.01,1,0.1)} & \cellcolor[rgb]{ .851,  .851,  .851}\textbf{0.32} & (1,0.01,0.01,0) & 0.10  & \cellcolor[rgb]{ .851,  .882,  .949}\textbf{(100,0.1,0.01,0)} & \cellcolor[rgb]{ .851,  .882,  .949}\textbf{0.09} \\
          & MAX   & (0.01,0.1,0) & 0.33  & (0.1,0.1,0) & 1.29  & \cellcolor[rgb]{ .851,  .851,  .851}\textbf{(0.1,0.01,0.1)} & \cellcolor[rgb]{ .851,  .851,  .851}\textbf{0.32} & (10,0,0.01,0) & 0.62  & \cellcolor[rgb]{ .851,  .882,  .949}\textbf{(100,0.1,0.01,0)} & \cellcolor[rgb]{ .851,  .882,  .949}\textbf{0.11} \\
          & akC   & \cellcolor[rgb]{ .851,  .851,  .851}\textbf{(0.01,0.01,0)} & \cellcolor[rgb]{ .851,  .851,  .851}\textbf{0.28} & (0.1,1,0) & 0.32  & (0.01,1,0.1) & 0.33  & \cellcolor[rgb]{ .851,  .882,  .949}\textbf{(100,0.1,1,0)} & \cellcolor[rgb]{ .851,  .882,  .949}\textbf{0.08} & (10,0,0.1,T) & 0.31 \\
          & MED   & \cellcolor[rgb]{ .851,  .851,  .851}\textbf{(0.1,0.1,0)} & \cellcolor[rgb]{ .851,  .851,  .851}\textbf{0.32} & (10,0,0) & 0.62  & \cellcolor[rgb]{ .851,  .851,  .851}\textbf{(0.1,1,0.1)} & \cellcolor[rgb]{ .851,  .851,  .851}\textbf{0.32} & (100,0.1,0.1,0) & 0.32  & \cellcolor[rgb]{ .851,  .882,  .949}\textbf{(100,0,0.1,26.5)} & \cellcolor[rgb]{ .851,  .882,  .949}\textbf{0.28} \\
          & TM    & \cellcolor[rgb]{ .851,  .851,  .851}\textbf{(0.01,0.1,0)} & \cellcolor[rgb]{ .851,  .851,  .851}\textbf{0.31} & (0.01,1,0) & 0.32  & (1,1,0.2) & 0.33  & (100,0.01,0.01,0) & 0.25  & \cellcolor[rgb]{ .851,  .882,  .949}\textbf{(10,0.01,0.01,T)} & \cellcolor[rgb]{ .851,  .882,  .949}\textbf{0.09} \\
        \hline
    \multirow{6}[1]{*}{\begin{sideways}Eilon\end{sideways}} & SUM   & \cellcolor[rgb]{ .851,  .851,  .851}\textbf{(0.01,0.01,0)} & \cellcolor[rgb]{ .851,  .851,  .851}\textbf{0.06} & (1,1,0) & 0.07  & (10,1,0.1) & 0.07  & (100,1,0.01,0) & 0.07  & \cellcolor[rgb]{ .851,  .882,  .949}\textbf{(0.01,0,100,0)} & \cellcolor[rgb]{ .851,  .882,  .949}\textbf{0.06} \\
          & kC    & \cellcolor[rgb]{ .851,  .851,  .851}\textbf{(0.01,0.01,0)} & \cellcolor[rgb]{ .851,  .851,  .851}\textbf{0.06} & (1,1,0) & 0.07  & (100,1,0.1) & 0.07  & (100,1,0.01,0) & 0.07  & (1,0.1,0.1,T) & 0.14 \\
          & MAX   & \cellcolor[rgb]{ .851,  .851,  .851}\textbf{(0.01,0.01,0)} & \cellcolor[rgb]{ .851,  .851,  .851}\textbf{0.06} & (10,0,0) & 0.25  & \cellcolor[rgb]{ .851,  .851,  .851}\textbf{(1,0.1,0.1)} & \cellcolor[rgb]{ .851,  .851,  .851}\textbf{0.06} & \cellcolor[rgb]{ .851,  .882,  .949}\textbf{(100,0.01,0.01,0)} & \cellcolor[rgb]{ .851,  .882,  .949}\textbf{0.06} & \cellcolor[rgb]{ .851,  .882,  .949}\textbf{(100,0.1,0.01,0)} & \cellcolor[rgb]{ .851,  .882,  .949}\textbf{0.06} \\
          & akC   & (0.01,0,0) & 0.07  & \cellcolor[rgb]{ .851,  .851,  .851}\textbf{(100,10,0)} & \cellcolor[rgb]{ .851,  .851,  .851}\textbf{0.06} & \cellcolor[rgb]{ .851,  .851,  .851}\textbf{(100,10,0.1)} & \cellcolor[rgb]{ .851,  .851,  .851}\textbf{0.06} & (100,1,0.1,T) & 0.73  & \cellcolor[rgb]{ .851,  .882,  .949}\textbf{(10,1,0.01,0.1)} & \cellcolor[rgb]{ .851,  .882,  .949}\textbf{0.06} \\
          & MED   & \cellcolor[rgb]{ .851,  .851,  .851}\textbf{(1,0,0)} & \cellcolor[rgb]{ .851,  .851,  .851}\textbf{0.07} & (0.1,0.1,0) & 0.20  & (0.01,0.1,13) & 0.08  & (1,0.1,0.1,T) & 0.20  & \cellcolor[rgb]{ .851,  .882,  .949}\textbf{(1,0.1,0.01,2.7)} & \cellcolor[rgb]{ .851,  .882,  .949}\textbf{0.06} \\
          & TM    & \cellcolor[rgb]{ .851,  .851,  .851}\textbf{(0.01,0.01,0)} & \cellcolor[rgb]{ .851,  .851,  .851}\textbf{0.06} & \cellcolor[rgb]{ .851,  .851,  .851}\textbf{(10,1,0)} & \cellcolor[rgb]{ .851,  .851,  .851}\textbf{0.06} & \cellcolor[rgb]{ .851,  .851,  .851}\textbf{(100,1,0.1)} & \cellcolor[rgb]{ .851,  .851,  .851}\textbf{0.06} & \cellcolor[rgb]{ .851,  .882,  .949}\textbf{(100,1,0.01,0)} & \cellcolor[rgb]{ .851,  .882,  .949}\textbf{0.06} & (100,10,100,0.2) & 0.38 \\
          \hline
    \multirow{6}[1]{*}{\begin{sideways}DW2\end{sideways}} & SUM   & \cellcolor[rgb]{ .851,  .851,  .851}\textbf{(100,0.1,0)} & \cellcolor[rgb]{ .851,  .851,  .851}\textbf{0.25} & \cellcolor[rgb]{ .851,  .851,  .851}\textbf{(100,0.1,0)} & \cellcolor[rgb]{ .851,  .851,  .851}\textbf{0.25} & \cellcolor[rgb]{ .851,  .851,  .851}\textbf{(100,0.1,0.4)} & \cellcolor[rgb]{ .851,  .851,  .851}\textbf{0.25} & \cellcolor[rgb]{ .851,  .882,  .949}\textbf{(1,0.1,0.1,T)} & \cellcolor[rgb]{ .851,  .882,  .949}\textbf{0.02} & \cellcolor[rgb]{ .851,  .882,  .949}\textbf{(1,0.01,0.1,T)} & \cellcolor[rgb]{ .851,  .882,  .949}\textbf{0.02} \\
          & kC    & \cellcolor[rgb]{ .851,  .851,  .851}\textbf{(0.01,0.1,0)} & \cellcolor[rgb]{ .851,  .851,  .851}\textbf{0.11} & (10,0.1,0) & 0.41  & (100,0.1,0.4) & 0.27  & \cellcolor[rgb]{ .851,  .882,  .949}\textbf{(1,0.1,0.1,T)} & \cellcolor[rgb]{ .851,  .882,  .949}\textbf{0.02} & (1,0,0.1,T) & 0.03 \\
          & MAX   & \cellcolor[rgb]{ .851,  .851,  .851}\textbf{(0.01,0.1,0)} & \cellcolor[rgb]{ .851,  .851,  .851}\textbf{0.11} & (0.01,1,0) & 0.45  & (100,0.01,0.3) & 0.31  & (100,0.1,0.01,0) & 0.10  & \cellcolor[rgb]{ .851,  .882,  .949}\textbf{(100,0.1,0.1,T)} & \cellcolor[rgb]{ .851,  .882,  .949}\textbf{0.03} \\
          & akC   & (0.1,0.1,0) & 0.47  & (100,0,0) & 0.19  & \cellcolor[rgb]{ .851,  .851,  .851}\textbf{(0.01,1,0.2)} & \cellcolor[rgb]{ .851,  .851,  .851}\textbf{0.10} & \cellcolor[rgb]{ .851,  .882,  .949}\textbf{(10,0,0.1,T)} & \cellcolor[rgb]{ .851,  .882,  .949}\textbf{0.02} & \cellcolor[rgb]{ .851,  .882,  .949}\textbf{(100,0,0.1,T)} & \cellcolor[rgb]{ .851,  .882,  .949}\textbf{0.02} \\
          & MED   & \cellcolor[rgb]{ .851,  .851,  .851}\textbf{(0.1,0.1,0)} & \cellcolor[rgb]{ .851,  .851,  .851}\textbf{0.11} & (100,0.01,0) & 0.15  & (100,0.01,125) & 0.20  & (10,0.01,1,0) & 0.15  & (100,0.01,0.1,T) & 0.20 \\
          & TM    & (1,0.1,0) & 0.30  & \cellcolor[rgb]{ .851,  .851,  .851}\textbf{(100,0.1,0)} & \cellcolor[rgb]{ .851,  .851,  .851}\textbf{0.18} & \cellcolor[rgb]{ .851,  .851,  .851}\textbf{(100,0.1,239.5)} & \cellcolor[rgb]{ .851,  .851,  .851}\textbf{0.18} & \cellcolor[rgb]{ .851,  .882,  .949}\textbf{(1,0.1,0.1,T)} & \cellcolor[rgb]{ .851,  .882,  .949}\textbf{0.02} & \cellcolor[rgb]{ .851,  .882,  .949}\textbf{(10,0,0.1,T)} & \cellcolor[rgb]{ .851,  .882,  .949}\textbf{0.02} \\
        \hline
    \multirow{6}[2]{*}{\begin{sideways}alcohol\end{sideways}} & SUM   & \cellcolor[rgb]{ .851,  .851,  .851}\textbf{(10,0.1,0)} & \cellcolor[rgb]{ .851,  .851,  .851}\textbf{1.26} & \cellcolor[rgb]{ .851,  .851,  .851}\textbf{(10,0.1,0)} & \cellcolor[rgb]{ .851,  .851,  .851}\textbf{1.26} & \cellcolor[rgb]{ .851,  .851,  .851}\textbf{(10,0.1,0.1)} & \cellcolor[rgb]{ .851,  .851,  .851}\textbf{1.26} & \cellcolor[rgb]{ .851,  .882,  .949}\textbf{(100,0,100,0)} & \cellcolor[rgb]{ .851,  .882,  .949}\textbf{0.85} & (100,0.1,100,0) & 0.94 \\
          & kC    & (100,10,0) & 145.05 & (100,0.01,0) & 1.29  & \cellcolor[rgb]{ .851,  .851,  .851}\textbf{(10,0.01,0.1)} & \cellcolor[rgb]{ .851,  .851,  .851}\textbf{1.20} & (10,1,100,0) & 1.85  & (100,0.1,100,0) & 1.21 \\
          & MAX   & (100,10,0) & 145.05 & (1,0,0) & 19.61 & \cellcolor[rgb]{ .851,  .851,  .851}\textbf{(0.1,0.1,0.1)} & \cellcolor[rgb]{ .851,  .851,  .851}\textbf{0.19} & (100,0.1,100,0) & 19.19 & (100,1,100,0) & 1.85 \\
          & akC   & (10,10,0) & 145.05 & (10,0.01,0) & 0.88  & \cellcolor[rgb]{ .851,  .851,  .851}\textbf{(100,0,T)} & \cellcolor[rgb]{ .851,  .851,  .851}\textbf{0.80} & (100,0.1,100,0) & 1.10  & (100,0.01,100,T) & 15.23 \\
          & MED   & (0.01,1,0) & 21.36 & (0.01,0.01,0) & 17.33 & \cellcolor[rgb]{ .851,  .851,  .851}\textbf{(0.1,0.01,709.6)} & \cellcolor[rgb]{ .851,  .851,  .851}\textbf{0.28} & (1,0.1,100,0) & 19.48 & (100,0.01,100,5.5) & 13.99 \\
          & TM    & (0.01,1,0) & 21.36 & \cellcolor[rgb]{ .851,  .851,  .851}\textbf{(100,0.1,0)} & \cellcolor[rgb]{ .851,  .851,  .851}\textbf{2.07} & (100,0.01,T) & 2.52  & \cellcolor[rgb]{ .851,  .882,  .949}\textbf{(10,0.1,100,0)} & \cellcolor[rgb]{ .851,  .882,  .949}\textbf{1.12} & (100,0.1,100,T) & 1.24 \\
       \hline
    \multirow{6}[2]{*}{\begin{sideways}DW3\end{sideways}} & SUM   & \cellcolor[rgb]{ .851,  .851,  .851}\textbf{(1,0.01,0)} & \cellcolor[rgb]{ .851,  .851,  .851}\textbf{0.01} & \cellcolor[rgb]{ .851,  .851,  .851}\textbf{(1,0.01,0)} & \cellcolor[rgb]{ .851,  .851,  .851}\textbf{0.01} & \cellcolor[rgb]{ .851,  .851,  .851}\textbf{(1,0.01,0.4)} & \cellcolor[rgb]{ .851,  .851,  .851}\textbf{0.01} & \cellcolor[rgb]{ .851,  .882,  .949}\textbf{(1,0.01,0.1,0.1)} & \cellcolor[rgb]{ .851,  .882,  .949}\textbf{0.01} & \cellcolor[rgb]{ .851,  .882,  .949}\textbf{(1,0.01,0.1,0)} & \cellcolor[rgb]{ .851,  .882,  .949}\textbf{0.01} \\
          & kC    & (0.01,0.1,0) & 0.17  & (10,0.1,0) & 0.01  & \cellcolor[rgb]{ .851,  .851,  .851}\textbf{(10,0,0.4)} & \cellcolor[rgb]{ .851,  .851,  .851}\textbf{0.01} & \cellcolor[rgb]{ .851,  .882,  .949}\textbf{(100,0.1,1,0.1)} & \cellcolor[rgb]{ .851,  .882,  .949}\textbf{0.01} & (1,0.01,0.1,0) & 0.01 \\
          & MAX   & (0.01,0.1,0) & 0.11  & (100,0.01,0) & 0.11  & \cellcolor[rgb]{ .851,  .851,  .851}\textbf{(10,0.1,0.3)} & \cellcolor[rgb]{ .851,  .851,  .851}\textbf{0.01} & (10,0.01,1,0) & 0.11  & \cellcolor[rgb]{ .851,  .882,  .949}\textbf{(10,0.1,1,0)} & \cellcolor[rgb]{ .851,  .882,  .949}\textbf{0.01} \\
          & akC   & (0.01,0,0) & 0.08  & \cellcolor[rgb]{ .851,  .851,  .851}\textbf{(100,0,0)} & \cellcolor[rgb]{ .851,  .851,  .851}\textbf{0.01} & \cellcolor[rgb]{ .851,  .851,  .851}\textbf{(100,0,T)} & \cellcolor[rgb]{ .851,  .851,  .851}\textbf{0.01} & (10,0,0.1,0.1) & 0.01  & \cellcolor[rgb]{ .851,  .882,  .949}\textbf{(100,0,1,T)} & \cellcolor[rgb]{ .851,  .882,  .949}\textbf{0.01} \\
          & MED   & (0.01,0.1,0) & 0.11  & (0.01,0,0) & 0.21  & \cellcolor[rgb]{ .851,  .851,  .851}\textbf{(100,0,T)} & \cellcolor[rgb]{ .851,  .851,  .851}\textbf{0.01} & (0.1,0.1,0.1,0) & 0.09  & \cellcolor[rgb]{ .851,  .882,  .949}\textbf{(100,0,1,66.2)} & \cellcolor[rgb]{ .851,  .882,  .949}\textbf{0.01} \\
          & TM    & (0.1,0,0) & 0.11  & (10,0,0) & 0.01  & \cellcolor[rgb]{ .851,  .851,  .851}\textbf{(10,0,T)} & \cellcolor[rgb]{ .851,  .851,  .851}\textbf{0.01} & (100,0.01,0.1,0.1) & 0.01  & (10,0.01,0.1,T) & 0.01 \\
        \hline
    \multirow{6}[2]{*}{\begin{sideways}HBK\end{sideways}} & SUM   & \cellcolor[rgb]{ .851,  .851,  .851}\textbf{(1,1,0)} & \cellcolor[rgb]{ .851,  .851,  .851}\textbf{0.72} & \cellcolor[rgb]{ .851,  .851,  .851}\textbf{(1,1,0)} & \cellcolor[rgb]{ .851,  .851,  .851}\textbf{0.72} & \cellcolor[rgb]{ .851,  .851,  .851}\textbf{(1,1,0.5)} & \cellcolor[rgb]{ .851,  .851,  .851}\textbf{0.72} & (10,0,1,0) & 1.82  & (10,0,1,0) & 1.82 \\
          & kC    & (10,10,0) & 43.45 & (0.01,0.1,0) & 8.83  & \cellcolor[rgb]{ .851,  .851,  .851}\textbf{(1,1,0.7)} & \cellcolor[rgb]{ .851,  .851,  .851}\textbf{0.72} & \cellcolor[rgb]{ .851,  .882,  .949}\textbf{(100,1,10,0)} & \cellcolor[rgb]{ .851,  .882,  .949}\textbf{0.56} & (10,0,1,0) & 1.82 \\
          & MAX   & \cellcolor[rgb]{ .851,  .851,  .851}\textbf{(0.01,0.1,0)} & \cellcolor[rgb]{ .851,  .851,  .851}\textbf{2.93} & (0.01,0,0) & 11.01 & (0.01,10,0.5) & 9.67  & (100,1,1,0) & 9.66  & \cellcolor[rgb]{ .851,  .882,  .949}\textbf{(100,1,10,0)} & \cellcolor[rgb]{ .851,  .882,  .949}\textbf{0.80} \\
          & akC   & (100,10,0) & 43.45 & \cellcolor[rgb]{ .851,  .851,  .851}\textbf{(0.01,10,0)} & \cellcolor[rgb]{ .851,  .851,  .851}\textbf{12.04} & \cellcolor[rgb]{ .851,  .851,  .851}\textbf{(0.01,10,0.3)} & \cellcolor[rgb]{ .851,  .851,  .851}\textbf{12.04} & \cellcolor[rgb]{ .851,  .882,  .949}\textbf{(100,0.01,10,0.1)} & \cellcolor[rgb]{ .851,  .882,  .949}\textbf{0.69} & (100,10,1,0.5) & 13.73 \\
          & MED   & \cellcolor[rgb]{ .851,  .851,  .851}\textbf{(1,0,0)} & \cellcolor[rgb]{ .851,  .851,  .851}\textbf{3.80} & (0.1,1,0) & 10.54 & (0.01,1,0.3) & 9.66  & (10,1,10,0) & 10.54 & (0.01,1,0.01,0.3) & 9.66 \\
          & TM    & (1,1,0) & 82.54 & \cellcolor[rgb]{ .851,  .851,  .851}\textbf{(0.01,1,0)} & \cellcolor[rgb]{ .851,  .851,  .851}\textbf{0.66} & (10,10,0.6) & 9.66  & (10,0.01,1,0) & 2.19  & (100,0.01,10,T) & 0.71 \\
    \hline
    \end{tabular}%
    }
    \end{adjustwidth}
  \label{tab:4MSE}%
\end{table}%

\newpage
\subsection{Graphical insights for evaluation metrics}
\label{graphical}

\begin{figure}[h!]
\caption{ Average evaluation metrics obtained in the training sets for\\ each hyperparameter combination in Quandt dataset.}
\hspace*{-3cm}
    \centering
    \begin{subfigure}{.55\textwidth}
    \caption{Raincloud for average MAE results}
        \centering
        \includegraphics[width=.85\linewidth]{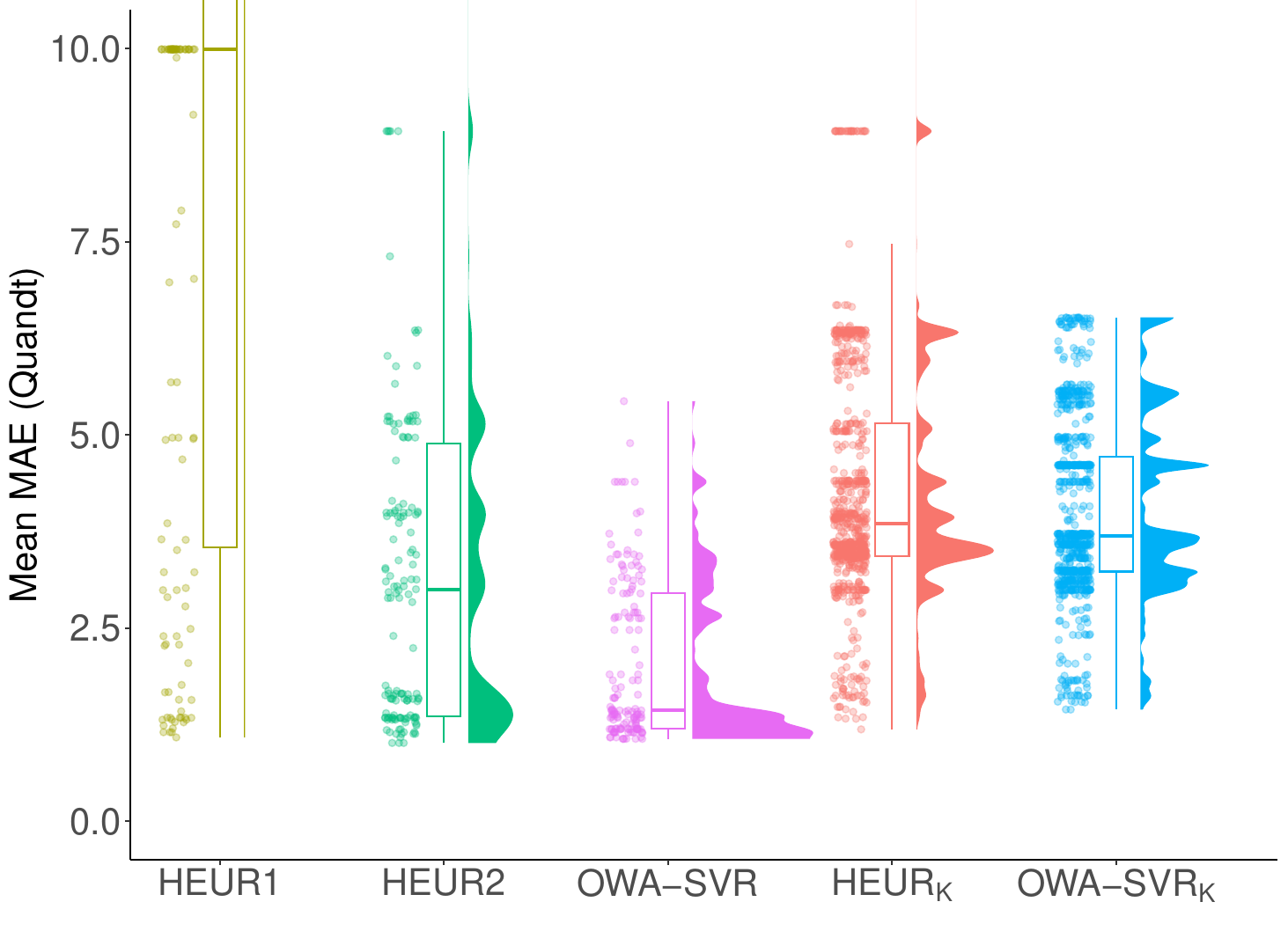}
    \end{subfigure}
    \begin{subfigure}{.55\textwidth}
    \caption{Raincloud for average MSE results}
        \centering
        \includegraphics[width=0.85\linewidth]{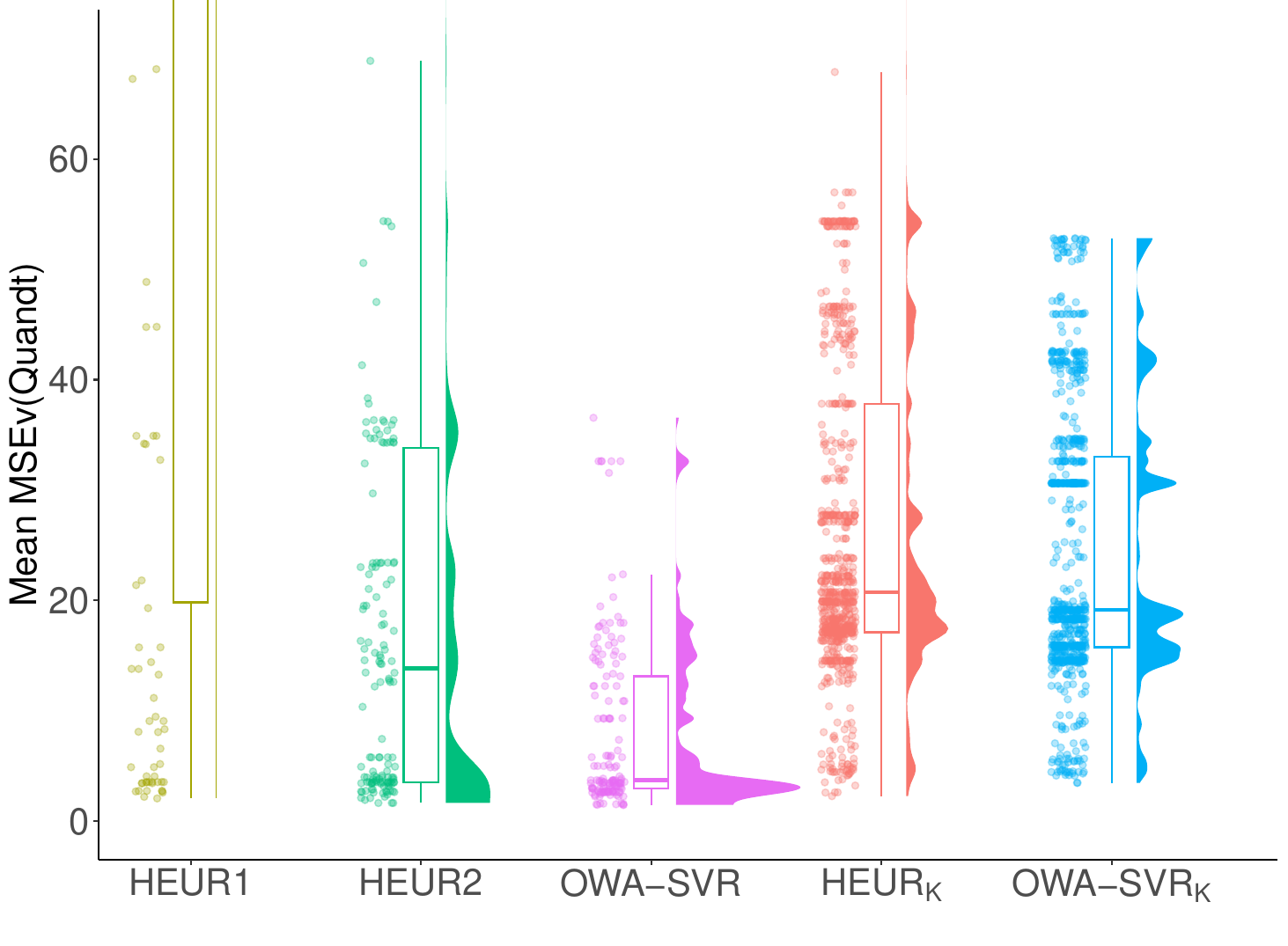}        
    \end{subfigure}
    \label{fig:Quandt}
\end{figure}

\begin{figure}[h!]
\caption{ Average evaluation metrics obtained in the training sets for\\ each hyperparameter combination in Boston dataset.}
\hspace*{-3cm}
    \centering
    \begin{subfigure}{.55\textwidth}
    \caption{Raincloud for average MAE results}
        \centering
        \includegraphics[width=0.85\linewidth]{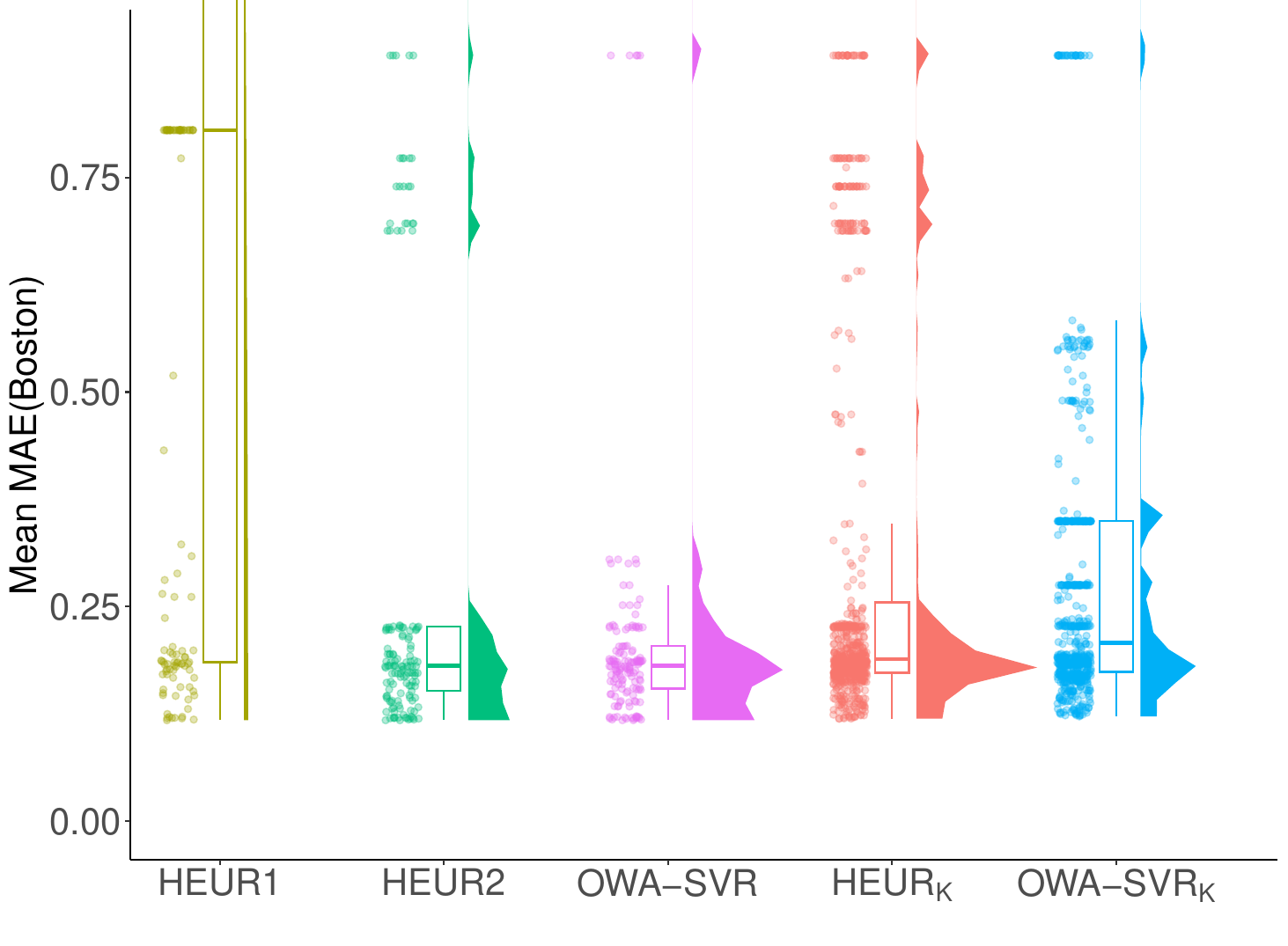}
    \end{subfigure}%
    \begin{subfigure}{0.55\textwidth}
    \caption{Raincloud for average MSE results}
        \centering
        \includegraphics[width=0.85\linewidth]{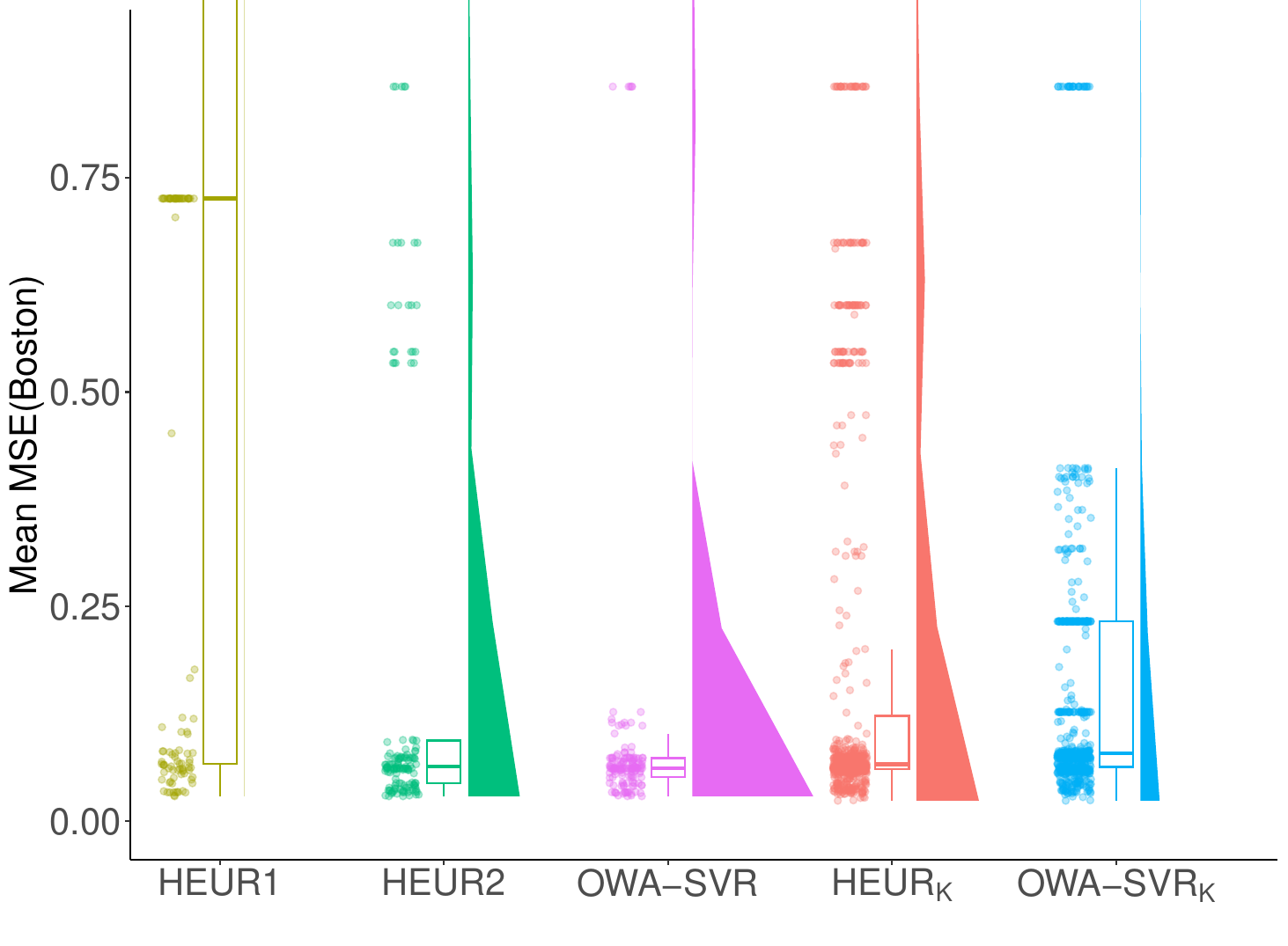}        
    \end{subfigure}
    \label{fig:Boston}
\end{figure}

\begin{figure}[h!]
\caption{ Average evaluation metrics obtained in the training sets for\\ each hyperparameter combination in CYGOB1 dataset.}
\hspace*{-3cm}
    \centering
    \begin{subfigure}{.55\textwidth}
    \caption{Raincloud for average MAE results}
        \centering
        \includegraphics[width=0.85\linewidth]{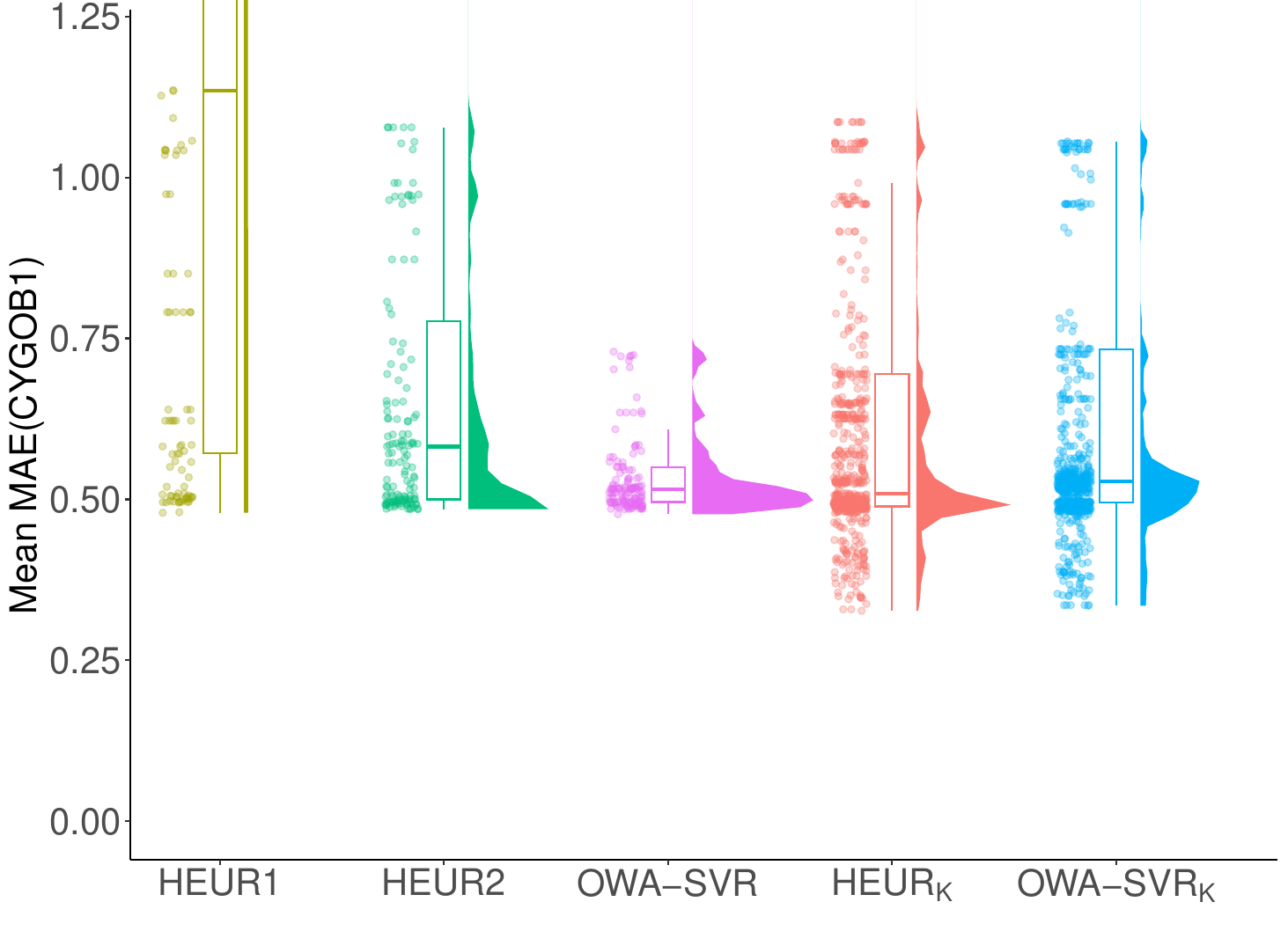}
    \end{subfigure}%
    \begin{subfigure}{0.55\textwidth}
    \caption{Raincloud for average MSE results}
        \centering
        \includegraphics[width=0.85\linewidth]{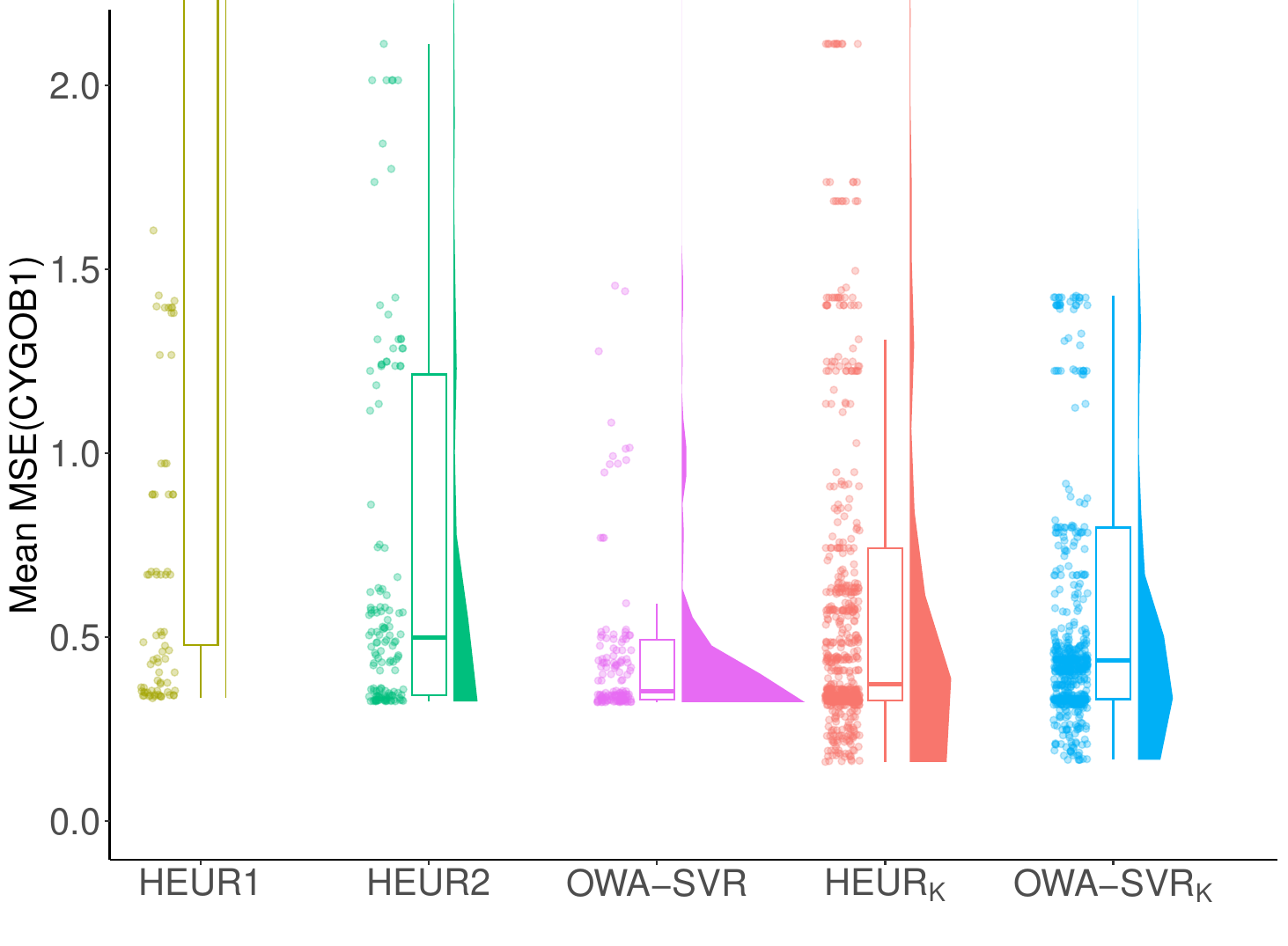}        
    \end{subfigure}
    \label{fig:CYGOB1}
\end{figure}

\begin{figure}[h!]
\caption{ Average evaluation metrics obtained in the training sets for\\ each hyperparameter combination in Eilon dataset.}
\hspace*{-3cm}
    \centering
    \begin{subfigure}{0.55\textwidth}
    \caption{Raincloud for average MAE results}
        \centering
        \includegraphics[width=0.85\linewidth]{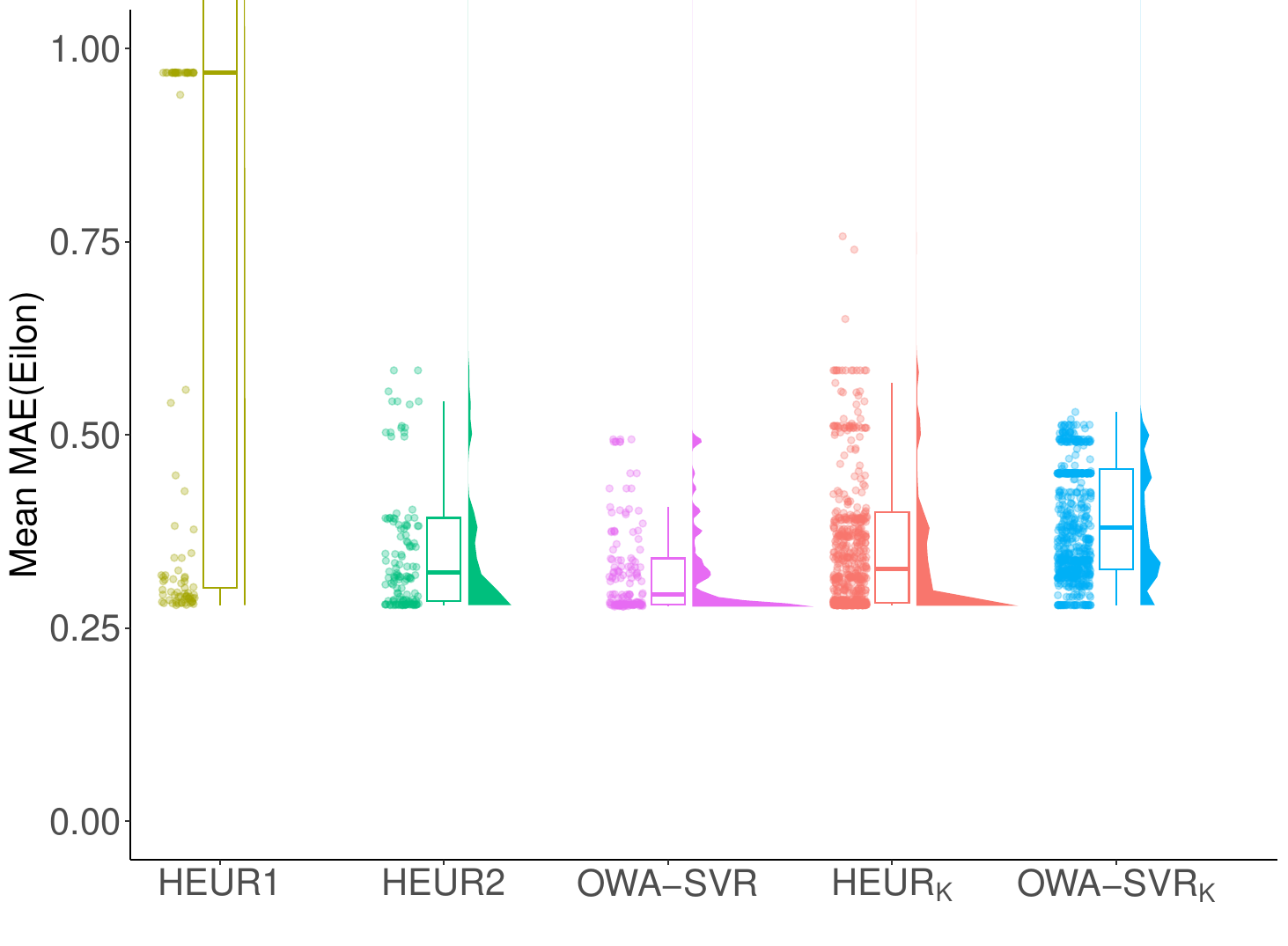}
    \end{subfigure}  
    \begin{subfigure}{0.55\textwidth}
    \caption{Raincloud for average MSE results}
        \centering
        \includegraphics[width=0.85\linewidth]{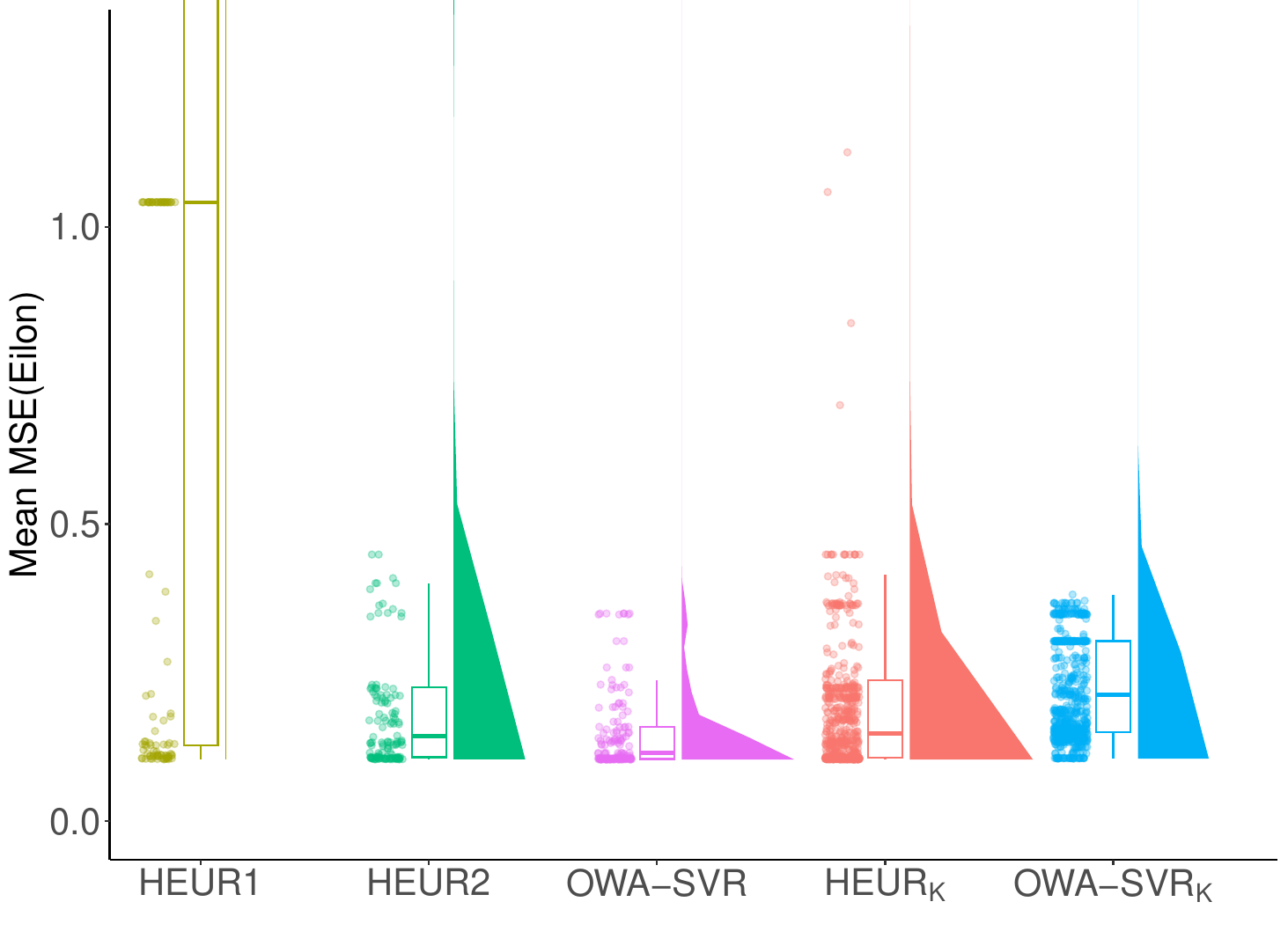}        
    \end{subfigure}
    \label{fig:Eilon}
\end{figure}

\begin{figure}[h!]
\caption{  Average evaluation metrics obtained in the training sets for\\ each hyperparameter combination in DW2 dataset.}
\hspace*{-3cm}
    \centering
    \begin{subfigure}{0.55\textwidth}
    \caption{Raincloud for average MAE results}
        \centering
        \includegraphics[width=0.85\linewidth]{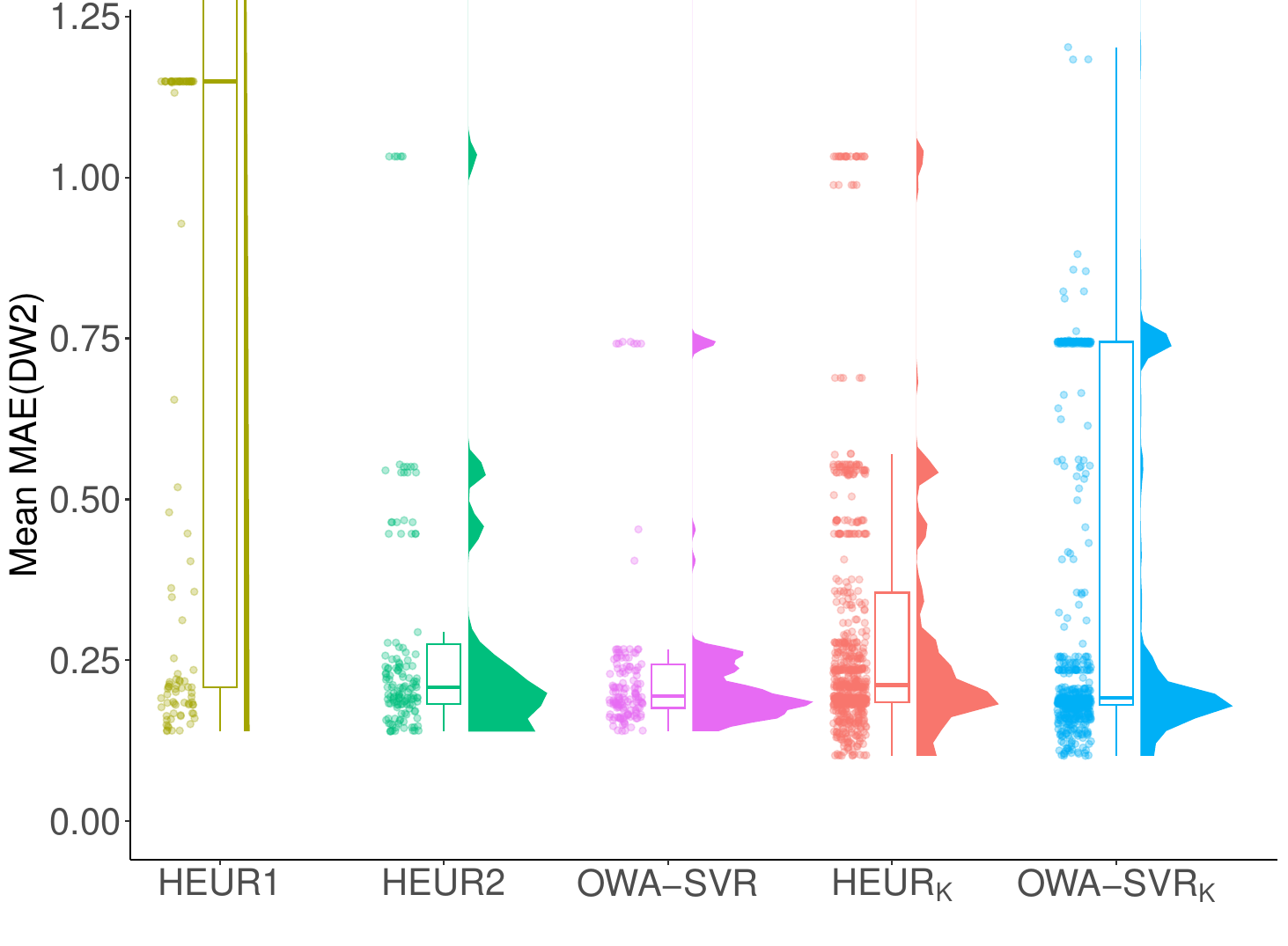}
    \end{subfigure} 
    \begin{subfigure}{0.55\textwidth}
    \caption{Raincloud for average MSE results}
        \centering
        \includegraphics[width=0.85\linewidth]{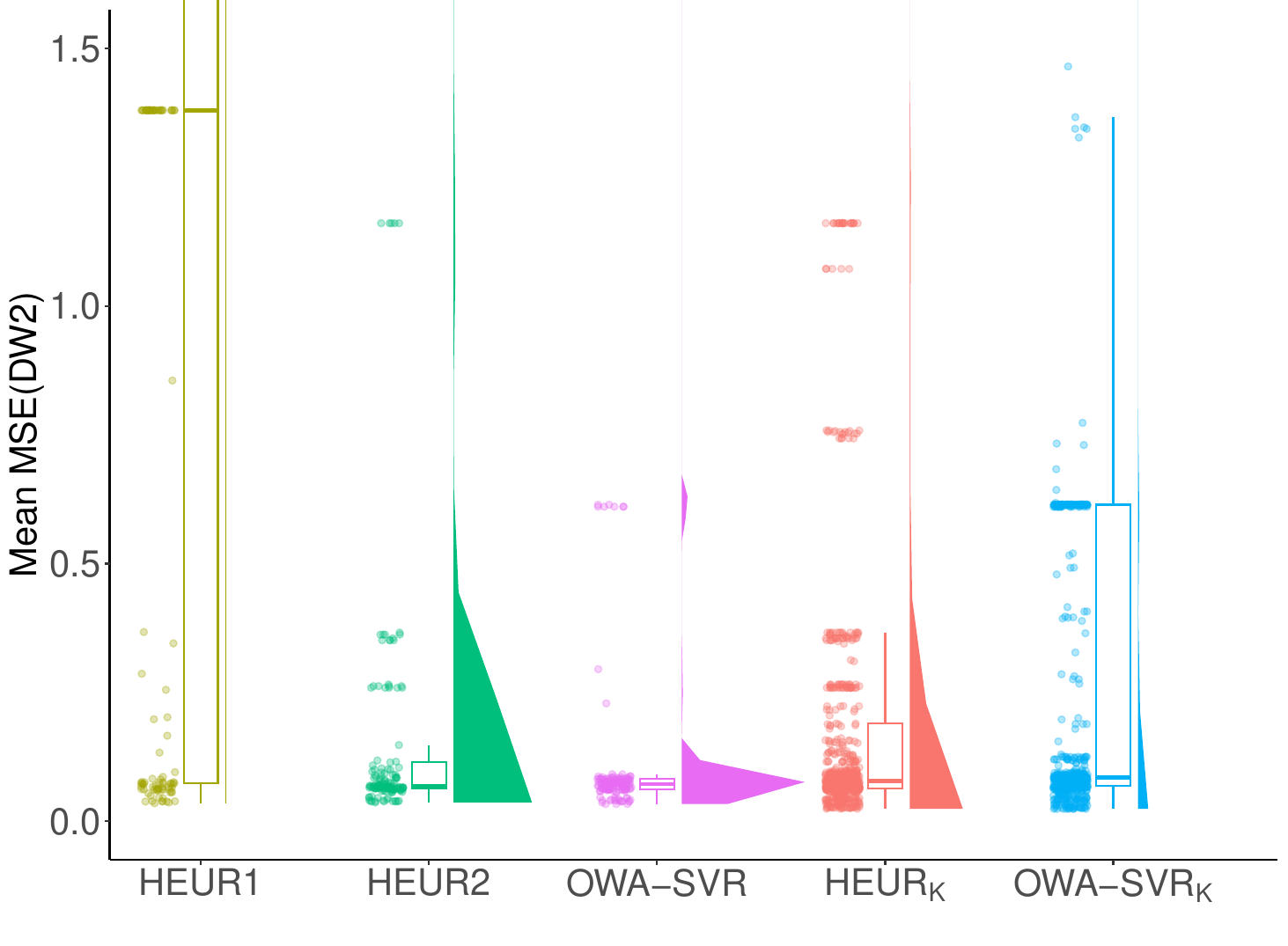}        
    \end{subfigure}
    \label{fig:DW2}
\end{figure}

\begin{figure}[h!]
\caption{ Average evaluation metrics obtained in the training sets for\\ each hyperparameter combination in Alcohol dataset.}
\hspace*{-3cm}
    \centering
    \begin{subfigure}{0.55\textwidth}
    \caption{Raincloud for average MAE results}
        \centering
        \includegraphics[width=0.85\linewidth]{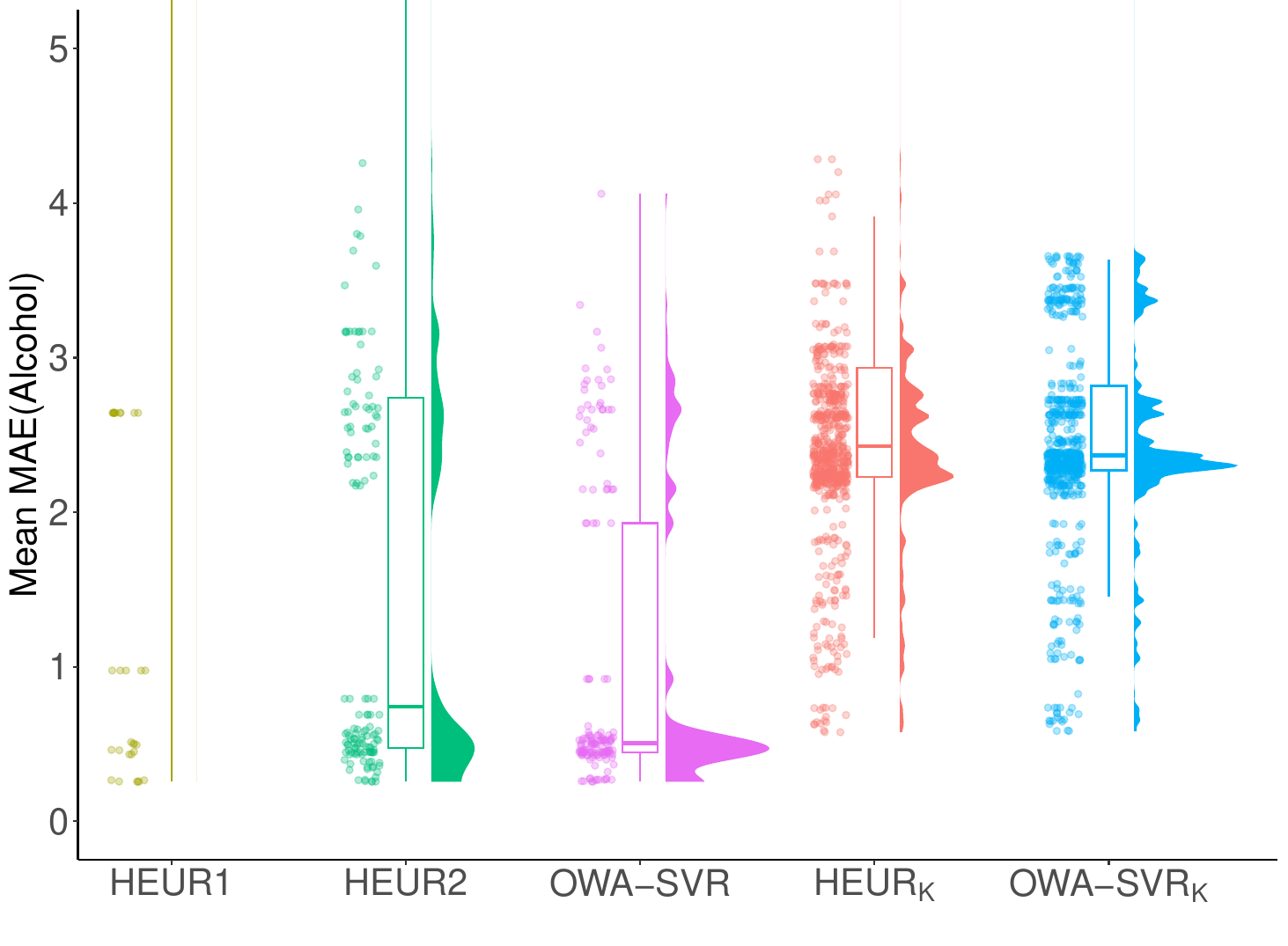}
    \end{subfigure}   
    \begin{subfigure}{0.55\textwidth}
    \caption{Raincloud for average MSE results}
        \centering
        \includegraphics[width=0.85\linewidth]{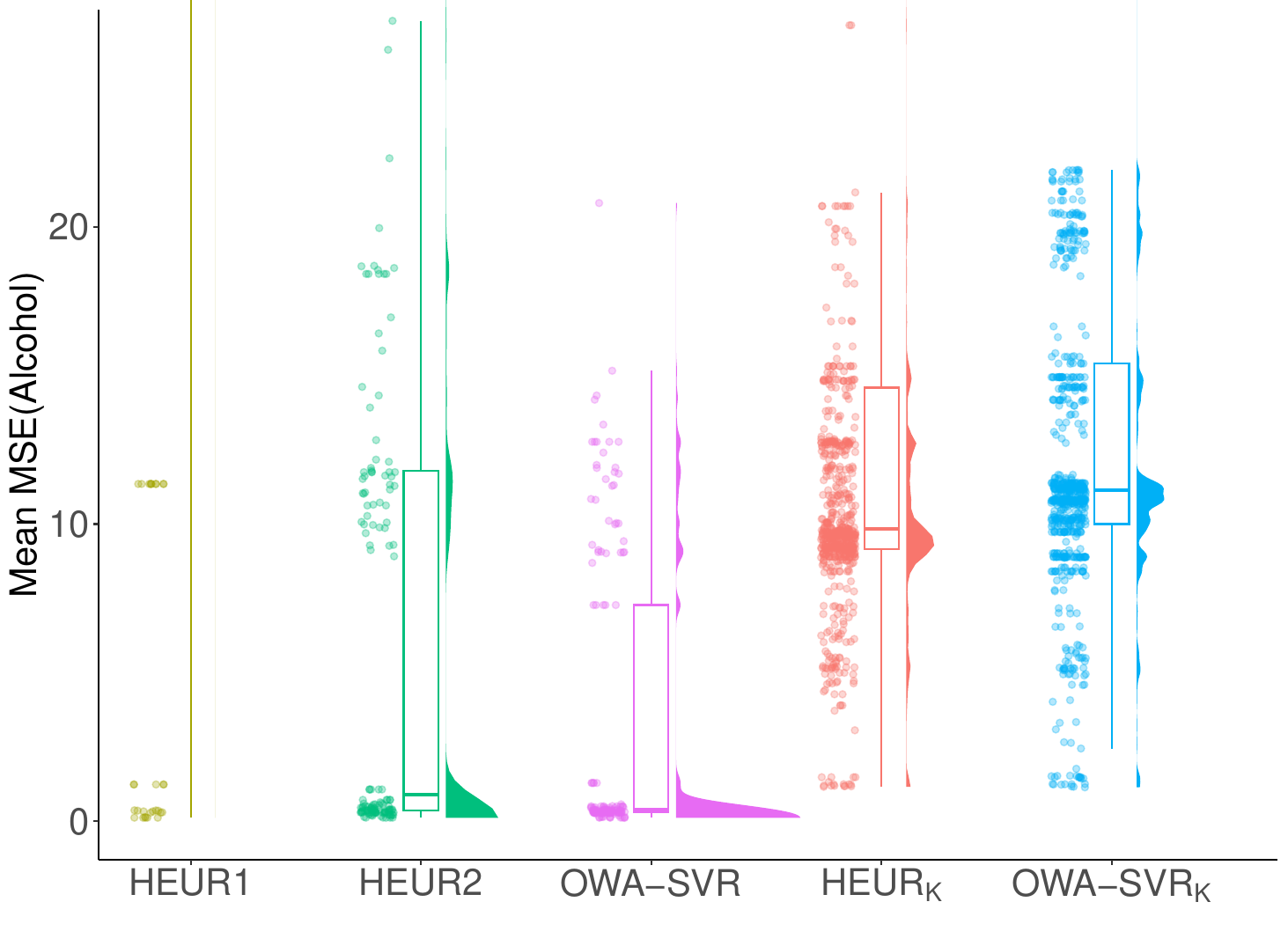}        
    \end{subfigure}
    \label{fig:Alochol}
\end{figure}

\begin{figure}[h!]
\caption{ Average evaluation metrics obtained in the training sets for\\ each hyperparameter combination in DW3 dataset.}
\hspace*{-3cm}
    \centering
    \begin{subfigure}{0.55\textwidth}
    \caption{Raincloud for average MAE results}
        \centering
        \includegraphics[width=0.85\linewidth]{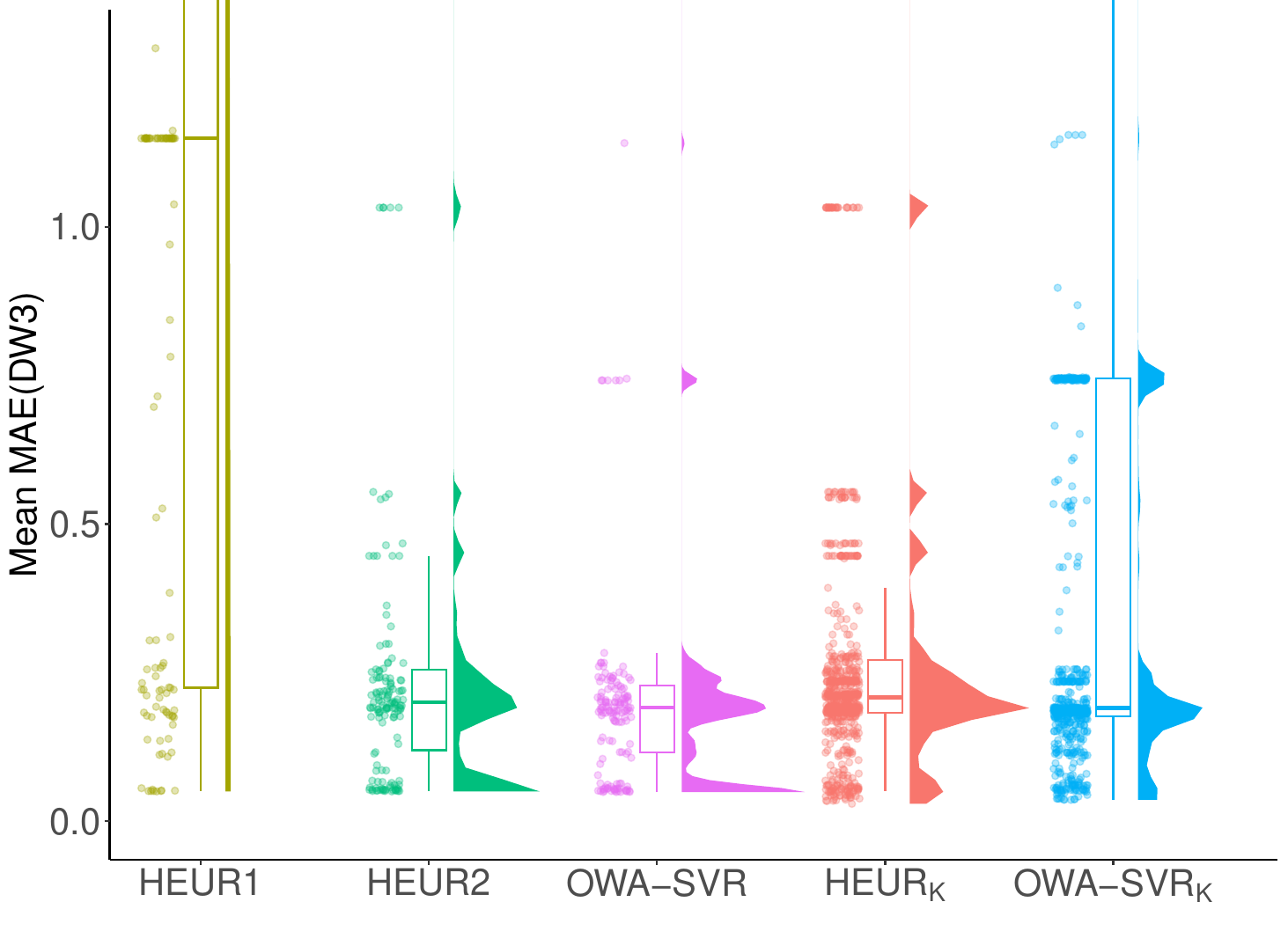}
    \end{subfigure}      
    \begin{subfigure}{0.55\textwidth}
      \caption{Raincloud for average MSE results}
        \centering
        \includegraphics[width=0.85\linewidth]{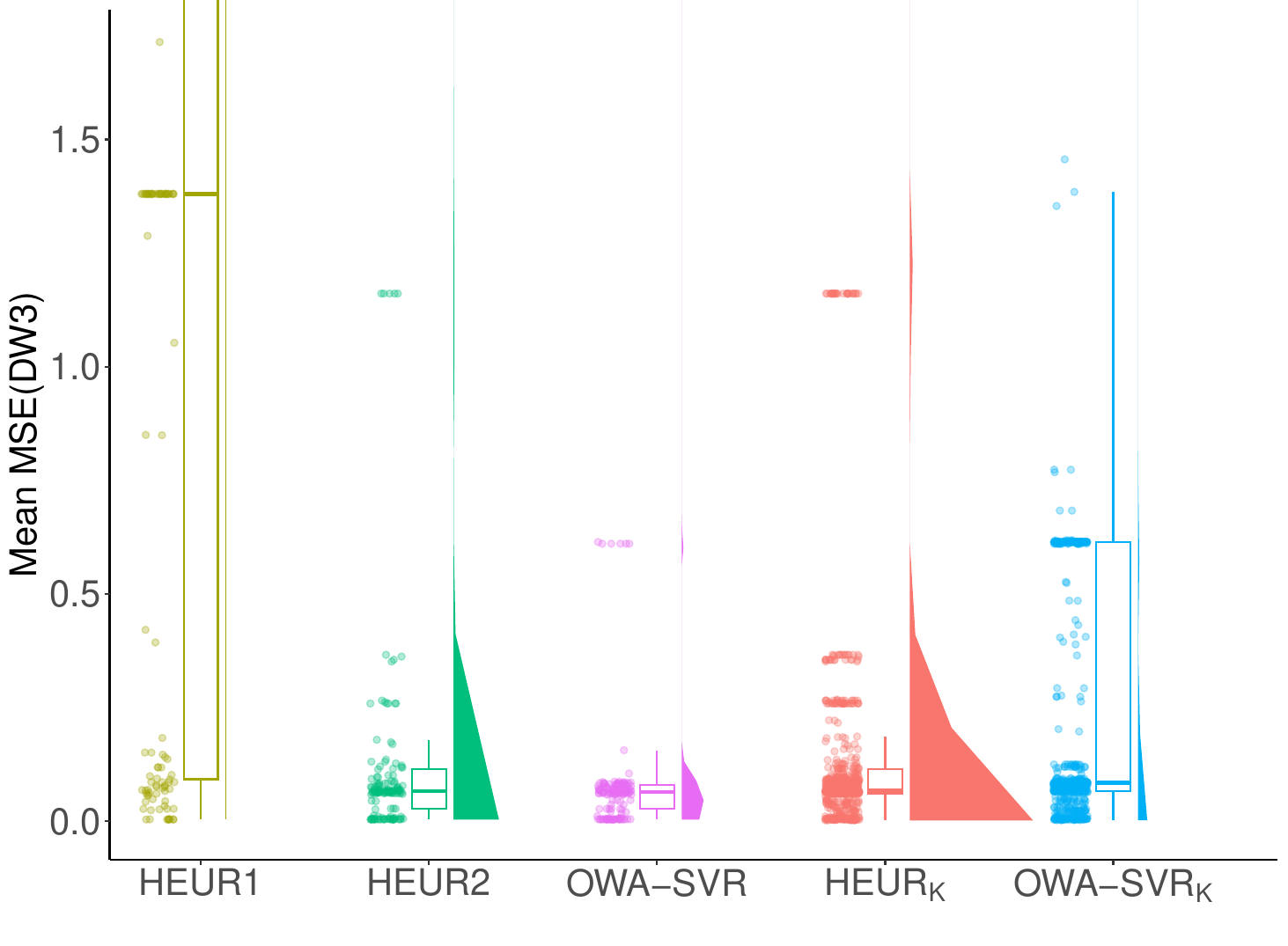}        
    \end{subfigure}   
    \label{fig:DW3}
\end{figure}

\begin{figure}[h!]
\caption{ Average evaluation metrics obtained in the training sets for\\ each hyperparameter combination in HBK dataset.}
\hspace*{-3cm}
    \centering
     \begin{subfigure}{0.55\textwidth}
     \caption{Raincloud for average MAE results}
        \centering
        \includegraphics[width=0.85\linewidth]{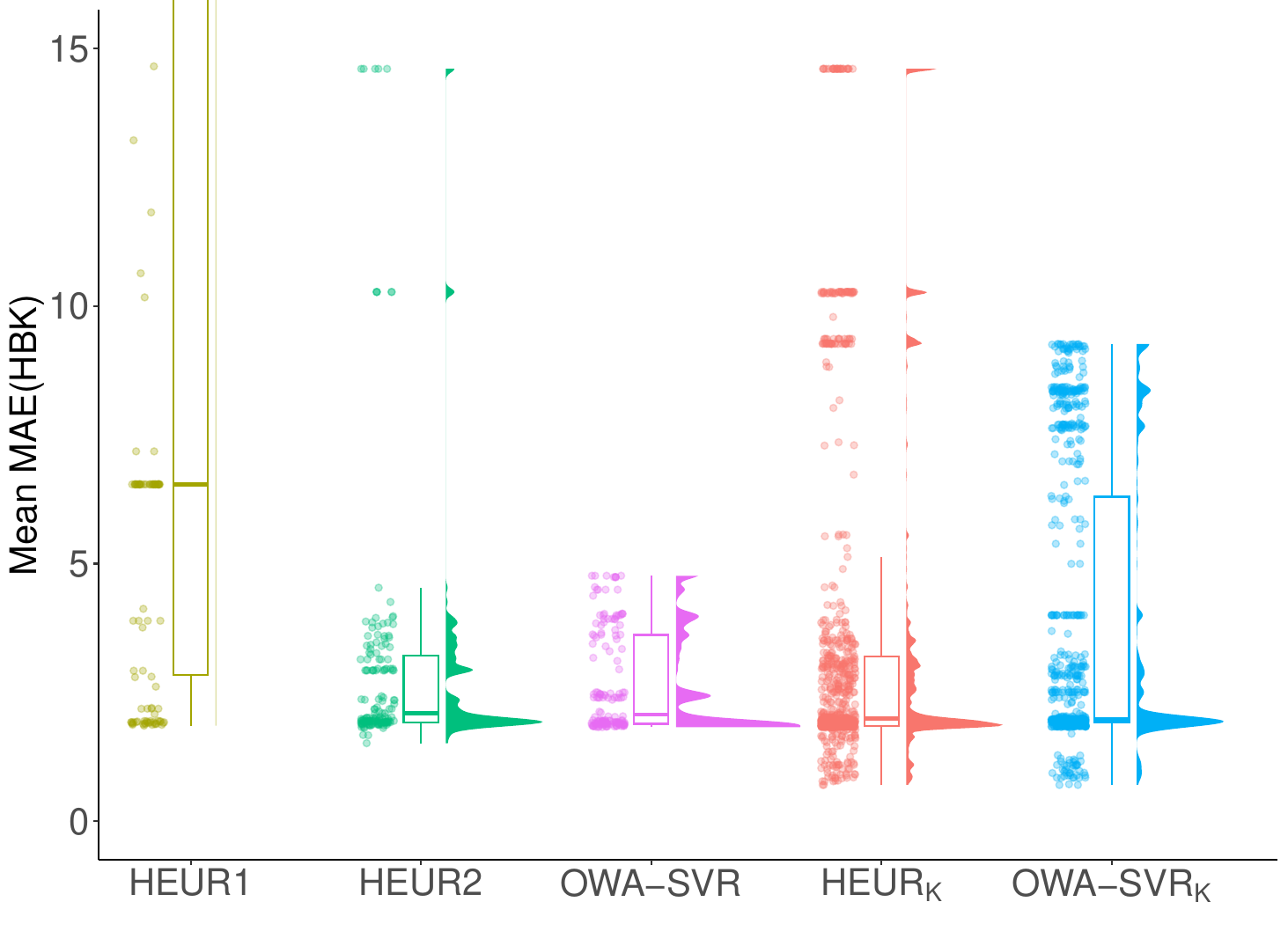}
    \end{subfigure}  
     \begin{subfigure}{0.55\textwidth}
     \caption{Raincloud for average MSE results}
        \centering
        \includegraphics[width=0.85\linewidth]{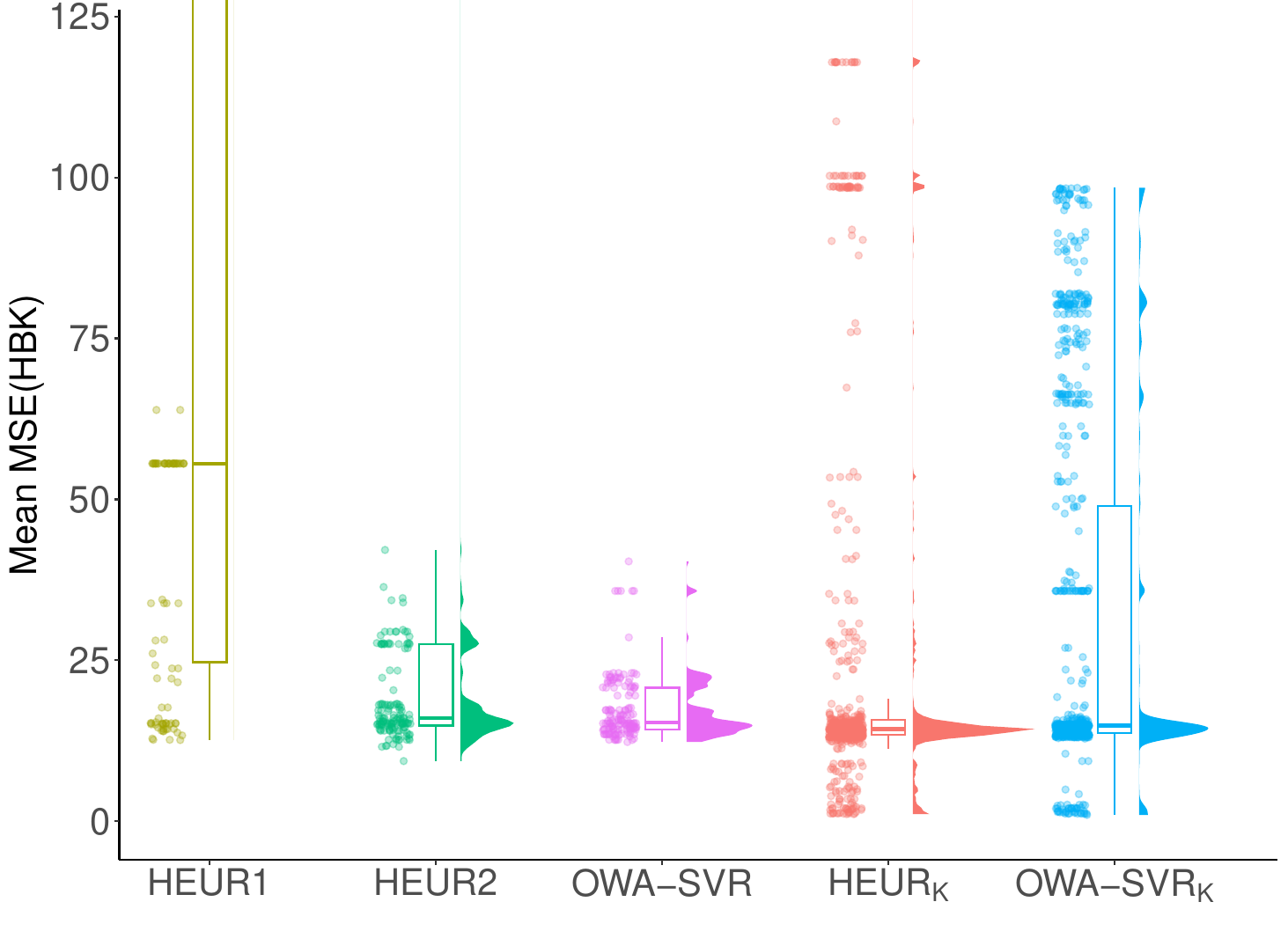}        
    \end{subfigure}
    \label{fig:HBK}
\end{figure}

This section seeks to offer a global insight on the performance of the models in terms of MAE and MSE. For this purpose, we report  two groups of raincloud plots related to each dataset (Figures \ref{fig:Quandt} to \ref{fig:HBK}). Each raincloud plot includes a  jitter plot, a boxplot (the horizontal bar represents the median), and an
illustration of the data distribution. For each dataset, the subfigure (a) illustrates the average MAE values obtained considering each studied weight and hyperparameter combination. Similarly, the subfigure (b) shows the raincloud plots considering the average MSE values.

Note that in some of the figures \ref{fig:Quandt} to \ref{fig:HBK} the median value in HEUR1 raincloud plots do not appear. This is due to the fact that this median is greater than the ones associated with the remaining approaches. Focusing on HEUR1, HEUR2 and OWA-SVR, we can conclude that, in most of the cases, the exact approach (OWA-SVR) provides the best median values. Besides, HEUR2 always present a better median value than HEUR1. Then, we could state that, for the linear case, the best heuristic approach seems to be HEUR2.

Regarding HEUR$_K$ and OWA-SVR$_K$, observe that they report a similar performance. In some cases, it is better to use the heuristic (Boston and CYGOB1), but in other instances, the exact approach provides better results (Quandt and Eilon).

In order to provide a deeper insight about the performance of MAE and MSE measures, some statistical test results are reported in Appendix \ref{app:test}. Particularly, we apply Kruskal-Wallis test using the mean values of MAE  obtained in the training. Note that, all combination of hypeparameter values and weights are considered. This test is performed for each dataset and it allows us to decide whether there is a significant difference between the mean MAE values of the different approaches (HEUR1, HEUR2, OWA-SVR, HEUR$_K$, OWA-SVR$_K$) or not. The same test is carried out using the mean MSE results for each dataset. Regarding the results, we observe that all $p$-values are smaller than $0.05$. Consequently, we can state that there is a significant difference between the approaches for mean MAE and MSE. Besides, Appendix \ref{app:test} reports pairwise Wilcoxon tests for each pair of approaches. We can appreciate that, in most of the cases, there exists a significant difference in the mean MAE and MSE values.


\clearpage

\subsection{Time performance analysis}\label{time}

In this section, we report the average results in terms of training sets with 75\% of the sample using the best hyperparameter values obtained in the  cross-validation process. To analyze this aspect we have generated a figure for each weight considered $\lambda$. In each figure, eight datasets appear: Quandt (size: 5), Boston (size: 5/6), CYGOB1 (size: 7),  Eilon (size: 7/8), DW2  (size: 10/11), Alcohol  (size: 6/7), DW3  (size: 10/11) and HBK  (size: 11/12). The goal of this section is to analyze how the performance time varies with the sample size  or with  multivariate data. 

\begin{figure}[!htb]
\hspace*{-3cm}
    \centering
    \begin{minipage}{.65\textwidth}
    \caption{Average computing time for SUM weight}
        \centering
        \includegraphics[width=0.85\linewidth]{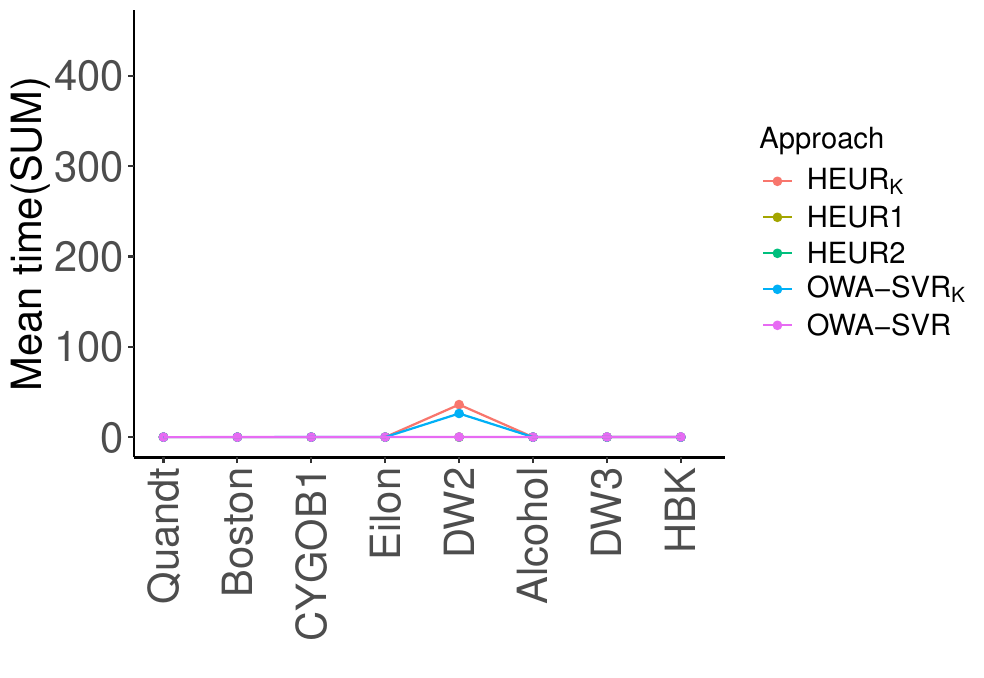}
        \label{SUM}
    \end{minipage}%
    \begin{minipage}{0.65\textwidth}
    \caption{Average computing time for kC weight}
        \centering
        \includegraphics[width=0.85\linewidth]{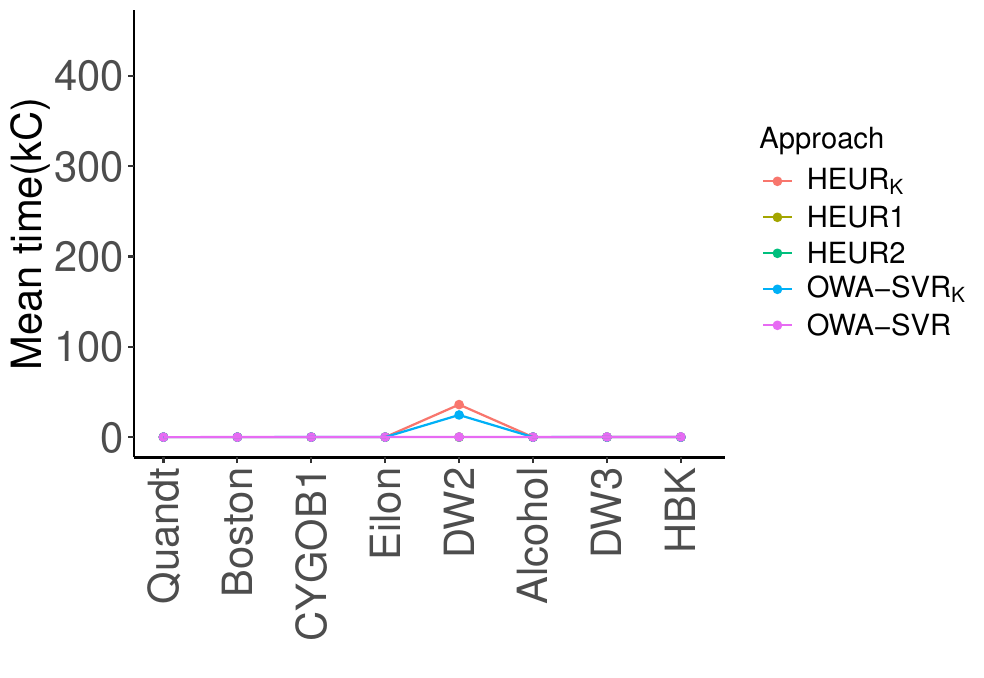}      
        \label{kC}
    \end{minipage}
\end{figure}

Regarding Figure \ref{SUM}, we can observe that the average computing times are very small in all the cases. Note that the largest times are associated with DW2 dataset, this may be due to the intrinsic characteristics of this dataset.

\begin{figure}[!htb]
\hspace*{-3cm}
    \centering
    \begin{minipage}{.65\textwidth}
    \caption{Average computing time for MAX weight}
        \centering
        \includegraphics[width=0.85\linewidth]{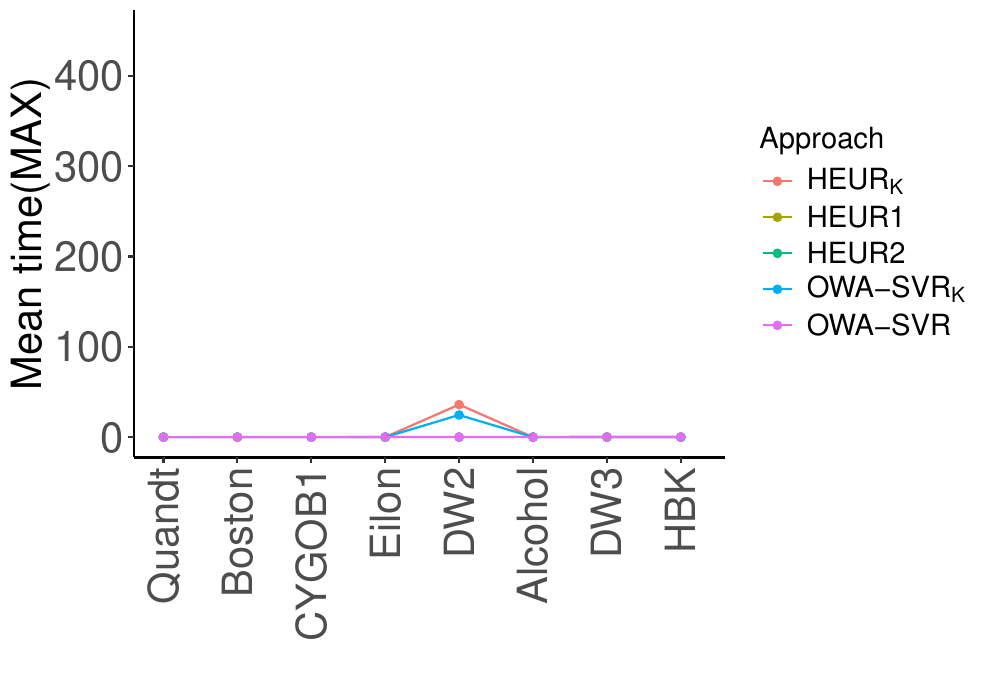}
        \label{MAX}
    \end{minipage}%
    \begin{minipage}{0.65\textwidth}
    \caption{Average computing time for akC weight}
        \centering
        \includegraphics[width=0.85\linewidth]{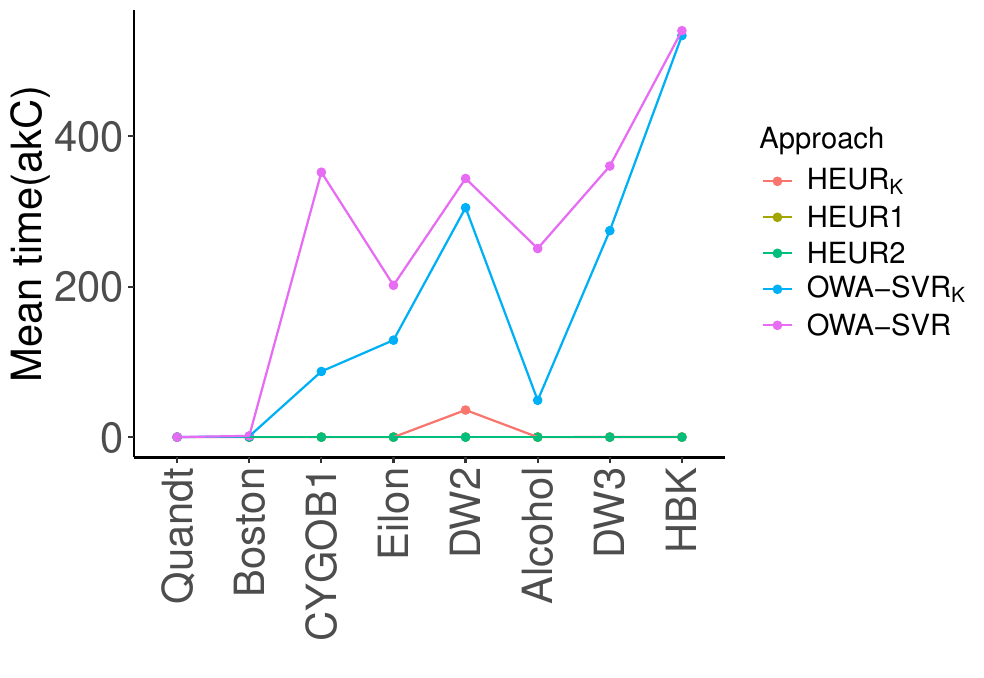} 
        \label{akC}
    \end{minipage}
\end{figure}

\begin{figure}[!htb]
\hspace*{-3cm}
    \centering
    \begin{minipage}{.65\textwidth}
    \caption{Average computing time for MED weight}
        \centering
        \includegraphics[width=0.85\linewidth]{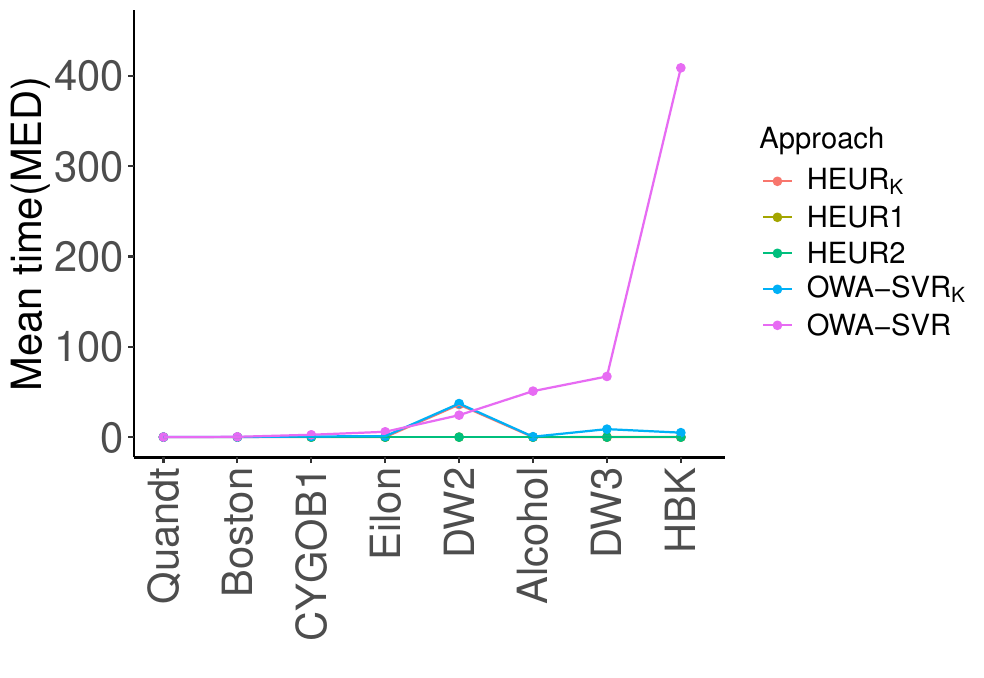}
        \label{MED}
    \end{minipage}%
    \begin{minipage}{0.65\textwidth}
    \caption{Average computing time for TM weight}
        \centering
        \includegraphics[width=0.85\linewidth]{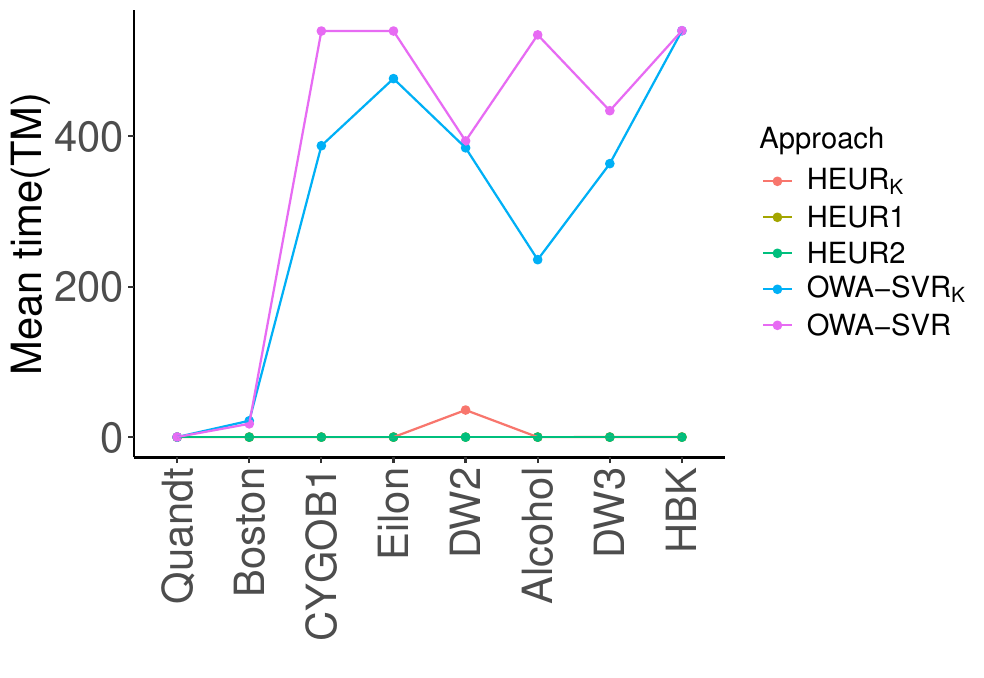}  
        \label{TM}
    \end{minipage}
\end{figure}

 In figures \ref{kC} and \ref{MAX} (akC and MAX weights), we can see a performance similar to the one of SUM weight. Recall that the considered weights in these figures are convex and the resulting exact approaches are easier to solve.

 Figure \ref{akC} reports a different behavior with respect to the previous figures. In this case, the exact approaches (OWA-SVR) and (OWA-SVR$_K$) require greater computing times. Figure \ref{MED} (MED weight) is also related to a non-convex operator but the computing times seem to be better than in the case of akC weight.

 Finally, Figure \ref{TM} reports the results of TM weight. We observe that the computing time is greater when the sample size and number of variables increase. From this study, we can conclude that largest times are obtained when using non-convex operators and exact approaches.

Next section is devoted to Friedman datasets which  contain larger sample size and number of variables. Then, we will see how the results obtained in these bigger datasets differ from the ones obtained in Section \ref{gof}, \ref{graphical} and \ref{time}.

\subsection{Analysis of Friedman datasets}
\label{friedman}

In this Section we focus on the results obtained from Friedman (150), Friedman (300) and Friedman (500). In this case, we only consider the heuristic approaches since previous experimental studies show that the time limit will be reached  in almost all the cases. 

First, we analyze the GoF results presented in Table \ref{tab:gof_friedman}. In this table we observe that, in general, best GoF values are obtained when using HEUR2 or HEUR$_K$. However, note that HEUR1 seems to obtain the best GoF for  MAX and MED weights. To explain this, recall that in these  weights, the OWA-SVR approaches try to optimize the deviation associated to just one position (maximum or median). Then, it makes sense that GoF value can be smaller than desired.

\begin{table}[htbp]
  \centering
  \caption{GoF for the different Friedman datasets and solution approaches}
  \scalebox{0.75}{
    \begin{tabular}{|c|l|lr|lr|lr|}
    \hline
    \multirow{2}[2]{*}{Data} & \multicolumn{1}{c|}{\multirow{2}[2]{*}{Weight}} & \multicolumn{2}{c|}{HEUR1} & \multicolumn{2}{c|}{HEUR2} & \multicolumn{2}{c|}{HEUR$_K$} \\
          &       & $(C,\varepsilon, t)$ & \multicolumn{1}{l|}{GoF} & $(C,\varepsilon, t)$ & \multicolumn{1}{l|}{GoF} & $(C,\varepsilon,\sigma, t)$ & \multicolumn{1}{l|}{GoF} \\
    \hline
    \multirow{6}[1]{*}{\begin{sideways}Friedman (150)\end{sideways}} & SUM   & \cellcolor[rgb]{ .859,  .859,  .859}\textbf{(10,10,0.1)} & \cellcolor[rgb]{ .859,  .859,  .859}\textbf{1.00} & \cellcolor[rgb]{ .859,  .859,  .859}\textbf{(10,10,0)} & \cellcolor[rgb]{ .859,  .859,  .859}\textbf{1.00} & \cellcolor[rgb]{ .851,  .882,  .949}\textbf{(10,10,0.01,0.1)} & \cellcolor[rgb]{ .851,  .882,  .949}\textbf{1.00} \\
          & kC    & (10,10,0.1) & 0.73  & \cellcolor[rgb]{ .859,  .859,  .859}\textbf{(10,10,0)} & \cellcolor[rgb]{ .859,  .859,  .859}\textbf{1.00} & \cellcolor[rgb]{ .851,  .882,  .949}\textbf{(10,10,0.01,0)} & \cellcolor[rgb]{ .851,  .882,  .949}\textbf{1.00} \\
          & MAX   & \cellcolor[rgb]{ .859,  .859,  .859}\textbf{(10,10,0.1)} & \cellcolor[rgb]{ .859,  .859,  .859}\textbf{0.66}  & (10,1,0) & -0.08 & (100,0.1,0.01,0) & 0.07 \\
          & akC   & (0.01,1,0.1) & 0.08  & \cellcolor[rgb]{ .859,  .859,  .859}\textbf{(100,1,0)} & \cellcolor[rgb]{ .859,  .859,  .859}\textbf{0.87} & \cellcolor[rgb]{ .851,  .882,  .949}\textbf{(10,0,0.01,0)} & \cellcolor[rgb]{ .851,  .882,  .949}\textbf{1.00} \\
          & MED   & \cellcolor[rgb]{ .859,  .859,  .859}\textbf{(0.01,1,0.1)} & \cellcolor[rgb]{ .859,  .859,  .859}\textbf{0.15} & \cellcolor[rgb]{ .859,  .859,  .859}\textbf{(100,0,0)} & \cellcolor[rgb]{ .859,  .859,  .859}\textbf{0.15} & (10,0.01,1,0.3) & 0.12 \\
          & TM    & (1,0.01,0.1) & 0.18  & \cellcolor[rgb]{ .859,  .859,  .859}\textbf{(10,1,0)} & \cellcolor[rgb]{ .859,  .859,  .859}\textbf{0.71} & \cellcolor[rgb]{ .851,  .882,  .949}\textbf{(100,1,1,0.4)} & \cellcolor[rgb]{ .851,  .882,  .949}\textbf{1.00} \\
          \hline
    \multirow{6}[0]{*}{\begin{sideways}Friedman (300)\end{sideways}} & SUM   & \cellcolor[rgb]{ .859,  .859,  .859}\textbf{(10,10,0.4)} & \cellcolor[rgb]{ .859,  .859,  .859}\textbf{1.00} & \cellcolor[rgb]{ .859,  .859,  .859}\textbf{(10,10,0.1)} & \cellcolor[rgb]{ .859,  .859,  .859}\textbf{1.00} & \cellcolor[rgb]{ .851,  .882,  .949}\textbf{(10,10,0.01,0.2)} & \cellcolor[rgb]{ .851,  .882,  .949}\textbf{1.00} \\
          & kC    & (1,10,0.4) & 0.81  & \cellcolor[rgb]{ .859,  .859,  .859}\textbf{(10,10,0)} & \cellcolor[rgb]{ .859,  .859,  .859}\textbf{1.00} & \cellcolor[rgb]{ .851,  .882,  .949}\textbf{(10,10,0.01,0.2)} & \cellcolor[rgb]{ .851,  .882,  .949}\textbf{1.00} \\
          & MAX   & \cellcolor[rgb]{ .859,  .859,  .859}\textbf{(1,10,0.4)} & \cellcolor[rgb]{ .859,  .859,  .859}\textbf{0.41} & (1,1,0) & -0.12 & (100,0,0.1,0.4) & -0.02 \\
          & akC   & (0.01,0,0.4) & -0.07 & \cellcolor[rgb]{ .859,  .859,  .859}\textbf{(100,1,0)} & \cellcolor[rgb]{ .859,  .859,  .859}\textbf{0.94} & \cellcolor[rgb]{ .851,  .882,  .949}\textbf{(10,0,0.1,0.5)} & \cellcolor[rgb]{ .851,  .882,  .949}\textbf{1.00} \\
          & MED   & \cellcolor[rgb]{ .859,  .859,  .859}\textbf{(100,0,0.4)} & \cellcolor[rgb]{ .859,  .859,  .859}\textbf{0.08} & (100,0.1,0) & 0.01  & (100,0,0.1,0.4) & 0.01 \\
          & TM    & (1,0,0.4) & 0.20  & \cellcolor[rgb]{ .859,  .859,  .859}\textbf{(100,1,0)} & \cellcolor[rgb]{ .859,  .859,  .859}\textbf{0.73} & \cellcolor[rgb]{ .851,  .882,  .949}\textbf{(100,1,1,1.7)} & \cellcolor[rgb]{ .851,  .882,  .949}\textbf{1.00} \\
          \hline
    \multirow{6}[1]{*}{\begin{sideways}Friedman (500)\end{sideways}} & SUM   & \cellcolor[rgb]{ .859,  .859,  .859}\textbf{(10,10,1.1)} & \cellcolor[rgb]{ .859,  .859,  .859}\textbf{1.00} & \cellcolor[rgb]{ .859,  .859,  .859}\textbf{(10,10,0)} & \cellcolor[rgb]{ .859,  .859,  .859}\textbf{1.00} & \cellcolor[rgb]{ .851,  .882,  .949}\textbf{(10,10,0.01,0.5)} & \cellcolor[rgb]{ .851,  .882,  .949}\textbf{1.00} \\
          & kC    & (1,10,1.1) & 0.84  & \cellcolor[rgb]{ .859,  .859,  .859}\textbf{(10,10,0)} & \cellcolor[rgb]{ .859,  .859,  .859}\textbf{1.00} & \cellcolor[rgb]{ .851,  .882,  .949}\textbf{(10,10,0.01,0.5)} & \cellcolor[rgb]{ .851,  .882,  .949}\textbf{1.00} \\
          & MAX   &\cellcolor[rgb]{ .859,  .859,  .859}\textbf{ (1,10,1.2) }& \cellcolor[rgb]{ .859,  .859,  .859}\textbf{0.24}  & (1,1,0) & -0.60 & (100,0.1,0.01,0.5) & -0.02 \\
          & akC   & (0.01,0,1.1) & -0.07 & \cellcolor[rgb]{ .859,  .859,  .859}\textbf{(100,1,0)} & \cellcolor[rgb]{ .859,  .859,  .859}\textbf{0.92} & \cellcolor[rgb]{ .851,  .882,  .949}\textbf{(10,0,0.01,0.5)} & \cellcolor[rgb]{ .851,  .882,  .949}\textbf{1.00} \\
          & MED   & \cellcolor[rgb]{ .859,  .859,  .859}\textbf{(0.01,0.01,1.1)} & \cellcolor[rgb]{ .859,  .859,  .859}\textbf{0.01} & \cellcolor[rgb]{ .859,  .859,  .859}\textbf{(1,0.1,0)} & \cellcolor[rgb]{ .859,  .859,  .859}\textbf{0.01} & \cellcolor[rgb]{ .851,  .882,  .949}\textbf{(100,0.1,10,4)} & \cellcolor[rgb]{ .851,  .882,  .949}\textbf{0.02} \\
          & TM    & (0.1,0.01,1.1) & 0.23  & \cellcolor[rgb]{ .859,  .859,  .859}\textbf{(10,1,0)} & \cellcolor[rgb]{ .859,  .859,  .859}\textbf{0.74} & \cellcolor[rgb]{ .851,  .882,  .949}\textbf{(100,1,1,5.1)} & \cellcolor[rgb]{ .851,  .882,  .949}\textbf{1.00} \\
    \hline
    \end{tabular}%
    }
  \label{tab:gof_friedman}%
\end{table}%

As in the previous datasets, to carry out this analysis, each dataset was divided into two parts: training set (75\% of the individuals) and test set (25\%). Then, a 5-fold cross-validation (5-FCV) was performed in each training set to determine the most appropriate parameter values in terms of different performance metrics. Table \ref{friedmanc} shows the characteristics of theses datasets in the fold-cross validation.

Tables \ref{tab:mae_friedman} and \ref{tab:mse_friedman} report the results of MAE and MSE respectively, for the best parameters values obtained in the training phase. We can derive the same conclusions as mentioned before. In general, the two-step heuristics (HEUR2 and HEUR$_K$) provide the best results in terms of MAE and MSE.

\begin{table}[htbp]
  \centering
        \caption{Characteristics of Friedman datasets}
    \begin{tabular}{lrrrr}
    \hline
    \hline
    Data  & \multicolumn{1}{l}{Number of } & \multicolumn{1}{l}{Explanatory} & \multicolumn{1}{l}{Folds size } & \multicolumn{1}{l}{Test set} \\
    name  & \multicolumn{1}{l}{observations} & \multicolumn{1}{l}{variables} & \multicolumn{1}{l}{in FCV} & \multicolumn{1}{l}{size}\\
    \hline
    \hline
    Friedman (150) & 150   & 8     & 22/23&38 \\
    Friedman (300) & 300    & 8     & 45&75 \\
    Friedman (500) & 500   &8     & 75 &125\\
    \hline
    \hline
    \end{tabular}%
  \label{friedmanc}
\end{table}%

\begin{table}[htbp]
  \centering
  \caption{MAE for the different Friedman datasets and solution approaches}
  \scalebox{0.75}{
    \begin{tabular}{|c|l|lr|lr|lr|}
    \hline
    \multirow{2}[2]{*}{Data} & \multicolumn{1}{c|}{\multirow{2}[2]{*}{Weight}} & \multicolumn{2}{c|}{HEUR1} & \multicolumn{2}{c|}{HEUR2} & \multicolumn{2}{c|}{HEUR$_K$} \\
          &       & Param & \multicolumn{1}{l|}{MAE} & Param & \multicolumn{1}{l|}{MAE} & Param & \multicolumn{1}{l|}{MAE} \\
    \hline
    \multirow{6}[1]{*}{\begin{sideways}Friedman (150)\end{sideways}} & SUM   & \cellcolor[rgb]{ .859,  .859,  .859}\textbf{(100,1,0)} & \cellcolor[rgb]{ .859,  .859,  .859}\textbf{2.36} & \cellcolor[rgb]{ .859,  .859,  .859}\textbf{(100,1,0)} & \cellcolor[rgb]{ .859,  .859,  .859}\textbf{2.36} & \cellcolor[rgb]{ .851,  .882,  .949}\textbf{(100,0.1,1,0.4)} & \cellcolor[rgb]{ .851,  .882,  .949}\textbf{1.31} \\
          & kC    & (0.1,0.1,0) & 4.32  & \cellcolor[rgb]{ .859,  .859,  .859}\textbf{(100,0.1,0)} & \cellcolor[rgb]{ .859,  .859,  .859}\textbf{2.48} & \cellcolor[rgb]{ .851,  .882,  .949}\textbf{(100,1,1,0.2)} & \cellcolor[rgb]{ .851,  .882,  .949}\textbf{1.63} \\
          & MAX   & \cellcolor[rgb]{ .859,  .859,  .859}\textbf{(10,10,0.1)} & \cellcolor[rgb]{ .859,  .859,  .859}\textbf{3.64} & (100,1,0) & 3.97  & (100,1,1,0.2) & 5.30 \\
          & akC   & (0.01,0.01,0) & 3.72  & \cellcolor[rgb]{ .859,  .859,  .859}\textbf{(10,1,0)} & \cellcolor[rgb]{ .859,  .859,  .859}\textbf{2.36} & \cellcolor[rgb]{ .851,  .882,  .949}\textbf{(100,0.01,1,0.2)} & \cellcolor[rgb]{ .851,  .882,  .949}\textbf{1.78} \\
          & MED   & \cellcolor[rgb]{ .859,  .859,  .859}\textbf{(0.1,0,0.1)} & \cellcolor[rgb]{ .859,  .859,  .859}\textbf{3.63} & (1,1,0) & 8.66  & (10,0.1,1,0.2) & 4.98 \\
          & TM    & (1,0,0) & 2.89  & \cellcolor[rgb]{ .859,  .859,  .859}\textbf{(10,1,0)} & \cellcolor[rgb]{ .859,  .859,  .859}\textbf{2.49} & \cellcolor[rgb]{ .851,  .882,  .949}\textbf{(100,0,1,0.2)} & \cellcolor[rgb]{ .851,  .882,  .949}\textbf{1.54} \\
          \hline
    \multirow{6}[0]{*}{\begin{sideways}Friedman (300)\end{sideways}} & SUM   & \cellcolor[rgb]{ .859,  .859,  .859}\textbf{(100,0,0.2)} & \cellcolor[rgb]{ .859,  .859,  .859}\textbf{2.31} & \cellcolor[rgb]{ .859,  .859,  .859}\textbf{(100,0,0)} & \cellcolor[rgb]{ .859,  .859,  .859}\textbf{2.31} & \cellcolor[rgb]{ .851,  .882,  .949}\textbf{(100,0.01,1,1.5)} & \cellcolor[rgb]{ .851,  .882,  .949}\textbf{0.82} \\
          & kC    & (100,10,0.2) & 3.36  & \cellcolor[rgb]{ .859,  .859,  .859}\textbf{(100,1,0)} & \cellcolor[rgb]{ .859,  .859,  .859}\textbf{2.42} & \cellcolor[rgb]{ .851,  .882,  .949}\textbf{(100,1,1,0.9)} & \cellcolor[rgb]{ .851,  .882,  .949}\textbf{1.20} \\
          & MAX   & \cellcolor[rgb]{ .859,  .859,  .859}\textbf{(10,10,0.2)} & \cellcolor[rgb]{ .859,  .859,  .859}\textbf{3.60} & (10,1,0) & 4.35  & (100,1,1,0.8) & 4.35 \\
          & akC   & (0.01,0.01,0.2) & 4.31  & \cellcolor[rgb]{ .859,  .859,  .859}\textbf{(100,0.01,0)} & \cellcolor[rgb]{ .859,  .859,  .859}\textbf{2.38} & \cellcolor[rgb]{ .851,  .882,  .949}\textbf{(100,0,1,0.9)} & \cellcolor[rgb]{ .851,  .882,  .949}\textbf{1.39} \\
          & MED   & \cellcolor[rgb]{ .859,  .859,  .859}\textbf{(0.1,0.01,0.2)} & \cellcolor[rgb]{ .859,  .859,  .859}\textbf{4.14} & (0.1,1,0) & 4.78  & (0.01,0.01,1,0.7) & 5.05 \\
          & TM    & (1,0.1,0.2) & 3.49  & \cellcolor[rgb]{ .859,  .859,  .859}\textbf{(10,1,0)} & \cellcolor[rgb]{ .859,  .859,  .859}\textbf{2.40} & \cellcolor[rgb]{ .851,  .882,  .949}\textbf{(100,0.01,1,0.9)} & \cellcolor[rgb]{ .851,  .882,  .949}\textbf{1.47} \\
          \hline
    \multirow{6}[1]{*}{\begin{sideways}Friedman (500)\end{sideways}} & SUM   & \cellcolor[rgb]{ .859,  .859,  .859}\textbf{(100,1,0.7)} & \cellcolor[rgb]{ .859,  .859,  .859}\textbf{1.99} & \cellcolor[rgb]{ .859,  .859,  .859}\textbf{(100,1,0.1)} & \cellcolor[rgb]{ .859,  .859,  .859}\textbf{1.99} & \cellcolor[rgb]{ .851,  .882,  .949}\textbf{(100,0,1,5)} & \cellcolor[rgb]{ .851,  .882,  .949}\textbf{0.49} \\
          & kC    & (10,10,0.6) & 3.71  & \cellcolor[rgb]{ .859,  .859,  .859}\textbf{(1,1,0)} & \cellcolor[rgb]{ .859,  .859,  .859}\textbf{2.45} & \cellcolor[rgb]{ .851,  .882,  .949}\textbf{(100,0.1,1,2.6)} & \cellcolor[rgb]{ .851,  .882,  .949}\textbf{0.84} \\
          & MAX   & \cellcolor[rgb]{ .859,  .859,  .859}\textbf{(10,10,0.6)} & \cellcolor[rgb]{ .859,  .859,  .859}\textbf{3.71} & (1,1,0) & 4.69  & (10,1,1,2.2) & 4.69 \\
          & akC   & (1,10,0.6) & 4.89  & \cellcolor[rgb]{ .859,  .859,  .859}\textbf{(100,1,0)} & \cellcolor[rgb]{ .859,  .859,  .859}\textbf{1.98} & \cellcolor[rgb]{ .851,  .882,  .949}\textbf{(100,0.1,1,2.7)} & \cellcolor[rgb]{ .851,  .882,  .949}\textbf{0.85} \\
          & MED   & \cellcolor[rgb]{ .859,  .859,  .859}\textbf{(0.01,0.1,0.6)} & \cellcolor[rgb]{ .859,  .859,  .859}\textbf{4.46} & (10,0.1,0) & 4.67  & (10,0.01,0.01,0.3) & 5.45 \\
          & TM    & (1,1,0.6) & 3.6   & \cellcolor[rgb]{ .859,  .859,  .859}\textbf{(10,1,0)} & \cellcolor[rgb]{ .859,  .859,  .859}\textbf{2.04} & \cellcolor[rgb]{ .851,  .882,  .949}\textbf{(100,0,1,2.7)} & \cellcolor[rgb]{ .851,  .882,  .949}\textbf{0.72} \\
    \hline
    \end{tabular}%
    }
  \label{tab:mae_friedman}%
\end{table}%

\begin{table}[htbp]
  \centering
  \caption{MSE for the different Friedman datasets and solution approaches}
  \scalebox{0.75}{
    \begin{tabular}{|c|l|lr|lr|lr|}
     \hline
    \multirow{2}[2]{*}{Data} & \multicolumn{1}{c|}{\multirow{2}[2]{*}{Weight}} & \multicolumn{2}{c|}{HEUR1} & \multicolumn{2}{c|}{HEUR2} & \multicolumn{2}{c|}{HEUR$_K$} \\
          &       & Param & \multicolumn{1}{l|}{MSE} & Param & \multicolumn{1}{l|}{MSE} & Param & \multicolumn{1}{l|}{MSE} \\
     \hline
    \multirow{6}[1]{*}{\begin{sideways}Friedman (150)\end{sideways}} & SUM   & \cellcolor[rgb]{ .859,  .859,  .859}\textbf{(10,1,0)} & \cellcolor[rgb]{ .859,  .859,  .859}\textbf{8.32} & \cellcolor[rgb]{ .859,  .859,  .859}\textbf{(10,1,0)} & \cellcolor[rgb]{ .859,  .859,  .859}\textbf{8.32} & \cellcolor[rgb]{ .851,  .882,  .949}\textbf{(100,0.1,1,0.4)} & \cellcolor[rgb]{ .851,  .882,  .949}\textbf{3.21} \\
          & kC    & (0.1,0.1,0.1) & 29.35 & \cellcolor[rgb]{ .859,  .859,  .859}\textbf{(10,0.01,0)}
          & \cellcolor[rgb]{ .859,  .859,  .859}\textbf{11.02}
          & \cellcolor[rgb]{ .851,  .882,  .949}\textbf{(100,1,1,0.2)} & \cellcolor[rgb]{ .851,  .882,  .949}\textbf{4.30} \\
          & MAX   & \cellcolor[rgb]{ .859,  .859,  .859}\textbf{(10,10,0)} & \cellcolor[rgb]{ .859,  .859,  .859}\textbf{19.85} & (100,1,0) & 22.63 & (100,1,1,0.2) & 40.41 \\
          & akC   & (0.01,0.01,0) & 19.12 & \cellcolor[rgb]{ .859,  .859,  .859}\textbf{(10,1,0)} & \cellcolor[rgb]{ .859,  .859,  .859}\textbf{8.41} & \cellcolor[rgb]{ .851,  .882,  .949}\textbf{(100,0.01,1,0.2)} & \cellcolor[rgb]{ .851,  .882,  .949}\textbf{5.84} \\
          & MED   & \cellcolor[rgb]{ .859,  .859,  .859}\textbf{(0.1,0,0)} & \cellcolor[rgb]{ .859,  .859,  .859}\textbf{18.40} & (1,1,0) & 93.04 & (10,0.1,1,0.2) & 36.02 \\
          & TM    & (1,0,0.1) & 12.79 & \cellcolor[rgb]{ .859,  .859,  .859}\textbf{(10,1,0)} & \cellcolor[rgb]{ .859,  .859,  .859}\textbf{9.50} & \cellcolor[rgb]{ .851,  .882,  .949}\textbf{(100,0,1,0.2)} & \cellcolor[rgb]{ .851,  .882,  .949}\textbf{5.19} \\
           \hline
    \multirow{6}[0]{*}{\begin{sideways}Friedman (300)\end{sideways}} & SUM   & \cellcolor[rgb]{ .859,  .859,  .859}(100,0,0.2) & \cellcolor[rgb]{ .859,  .859,  .859}9.62 & \cellcolor[rgb]{ .859,  .859,  .859}(100,0,0) & \cellcolor[rgb]{ .859,  .859,  .859}9.62 & \cellcolor[rgb]{ .851,  .882,  .949}\textbf{(100,0.01,1,1.6)} & \cellcolor[rgb]{ .851,  .882,  .949}\textbf{1.75} \\
          & kC    & (100,10,0.2) & 16.14 & \cellcolor[rgb]{ .859,  .859,  .859}\textbf{(100,1,0)} & \cellcolor[rgb]{ .859,  .859,  .859}\textbf{9.23} & \cellcolor[rgb]{ .851,  .882,  .949}\textbf{(100,1,1,0.7)} & \cellcolor[rgb]{ .851,  .882,  .949}\textbf{2.66} \\
          & MAX   & \cellcolor[rgb]{ .859,  .859,  .859}\textbf{(10,10,0.2)} & \cellcolor[rgb]{ .859,  .859,  .859}\textbf{20.18} & (10,1,0) & 26.97 & (10,1,1,0.9) & 26.97 \\
          & akC   & (0.01,0.01,0.2) & 28.34 & \cellcolor[rgb]{ .859,  .859,  .859}\textbf{(10,1,0)} & \cellcolor[rgb]{ .859,  .859,  .859}\textbf{9.27} & \cellcolor[rgb]{ .851,  .882,  .949}\textbf{(100,0.01,1,0.9)} & \cellcolor[rgb]{ .851,  .882,  .949}\textbf{5.03} \\
          & MED   & \cellcolor[rgb]{ .859,  .859,  .859}\textbf{(0.01,0.01,0.2)} & \cellcolor[rgb]{ .859,  .859,  .859}\textbf{26.19} & (0.1,1,0) & 31.74 & (100,0.1,1,0.7) & 26.72 \\
          & TM    & (1,0.1,0.2) & 21.56 & \cellcolor[rgb]{ .859,  .859,  .859}\textbf{(10,1,0)} & \cellcolor[rgb]{ .859,  .859,  .859}\textbf{9.88} & \cellcolor[rgb]{ .851,  .882,  .949}\textbf{(100,0.01,1,0.9)} & \cellcolor[rgb]{ .851,  .882,  .949}\textbf{6.87} \\
           \hline
    \multirow{6}[1]{*}{\begin{sideways}Friedman (500)\end{sideways}} & SUM   & \cellcolor[rgb]{ .859,  .859,  .859}\textbf{(10,1,0.6)} & \cellcolor[rgb]{ .859,  .859,  .859}\textbf{6.66} & \cellcolor[rgb]{ .859,  .859,  .859}\textbf{(10,1,0.1)} & \cellcolor[rgb]{ .859,  .859,  .859}\textbf{6.66} & \cellcolor[rgb]{ .851,  .882,  .949}\textbf{(100,0.01,1,4.9)} & \cellcolor[rgb]{ .851,  .882,  .949}\textbf{0.42} \\
          & kC    & (1,10,0.6) & 21.49 & \cellcolor[rgb]{ .859,  .859,  .859}\textbf{(1,1,0)} & \cellcolor[rgb]{ .859,  .859,  .859}\textbf{9.47} & \cellcolor[rgb]{ .851,  .882,  .949}\textbf{(100,0.1,1,2.2)} & \cellcolor[rgb]{ .851,  .882,  .949}\textbf{1.06} \\
          & MAX   & \cellcolor[rgb]{ .859,  .859,  .859}\textbf{(10,10,0.6)} & \cellcolor[rgb]{ .859,  .859,  .859}\textbf{22.27} & (1,1,0) & 31.37 & (10,1,1,2.6) & 31.37 \\
          & akC   & (0.01,0.1,0.6) & 34.39 & \cellcolor[rgb]{ .859,  .859,  .859}\textbf{(10,1,0)} & \cellcolor[rgb]{ .859,  .859,  .859}\textbf{6.67} & \cellcolor[rgb]{ .851,  .882,  .949}\textbf{(100,0,1,2.7)} & \cellcolor[rgb]{ .851,  .882,  .949}\textbf{1.14} \\
          & MED   & \cellcolor[rgb]{ .859,  .859,  .859}\textbf{(0.01,0.1,0.6)} & \cellcolor[rgb]{ .859,  .859,  .859}\textbf{29.27} & (10,0.1,0) & 35.43 & (10,0.01,0.01,0.3) & 41.91 \\
          & TM    & (1,1,0.6) & 19.94 & \cellcolor[rgb]{ .859,  .859,  .859}\textbf{(10,1,0)} & \cellcolor[rgb]{ .859,  .859,  .859}\textbf{6.85} & \cellcolor[rgb]{ .851,  .882,  .949}\textbf{(100,0.1,1,2.7)} & \cellcolor[rgb]{ .851,  .882,  .949}\textbf{1.03} \\
    \hline
    \end{tabular}%
    }
  \label{tab:mse_friedman}%
\end{table}%

\begin{figure}[!htb]
 \caption{ Average MAE and MSE  training values for Friedman(150)  dataset}
 \hspace*{-3cm}
    \centering
    \begin{minipage}{.55\textwidth}
        \centering
        \includegraphics[width=0.85\linewidth]{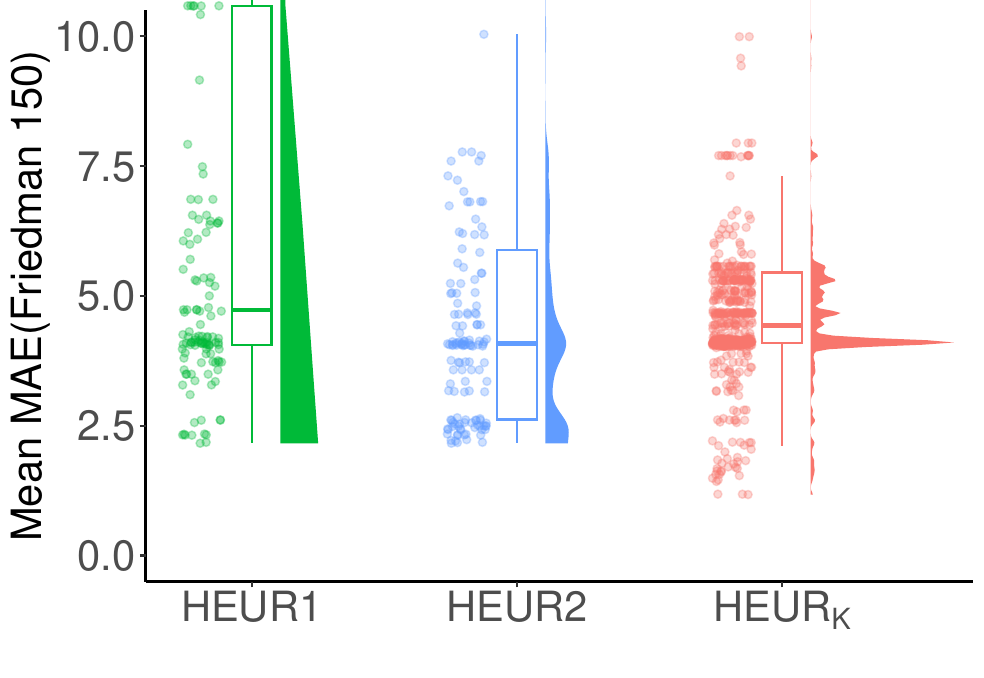}
    \end{minipage}%
    \begin{minipage}{0.55\textwidth}
        \centering
        \includegraphics[width=0.85\linewidth]{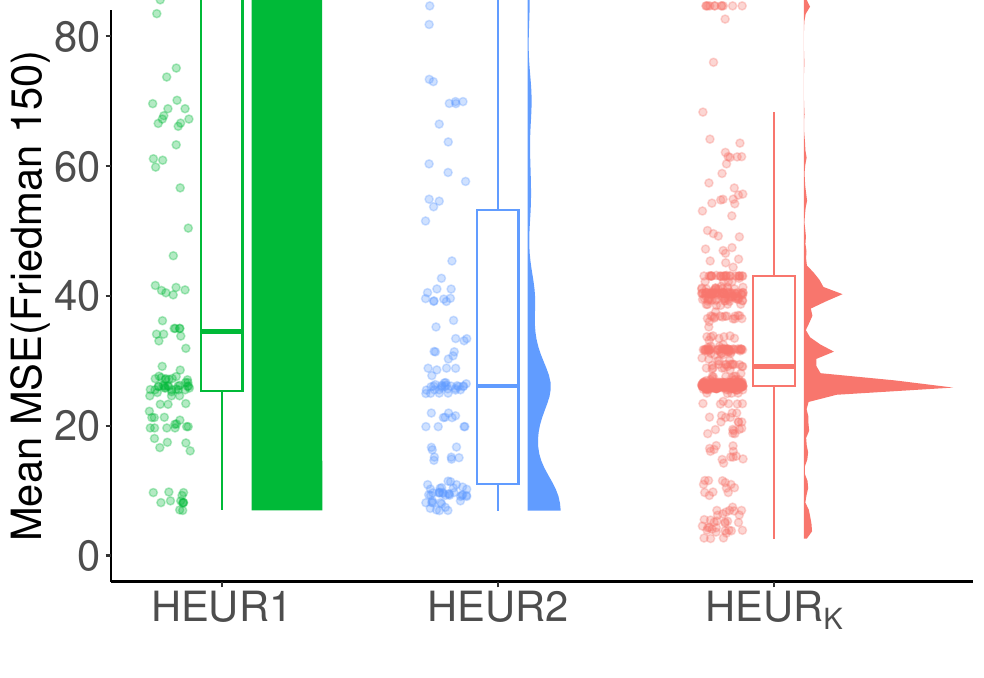}        
    \end{minipage}

    \label{fig:FRIEDMAN_150}
\end{figure}

\begin{figure}[!htb]
     \caption{ Average MAE and MSE  training values for Friedman(300) dataset}
     \hspace*{-3cm}
    \centering
    \begin{minipage}{.55\textwidth}
        \centering
        \includegraphics[width=0.85\linewidth]{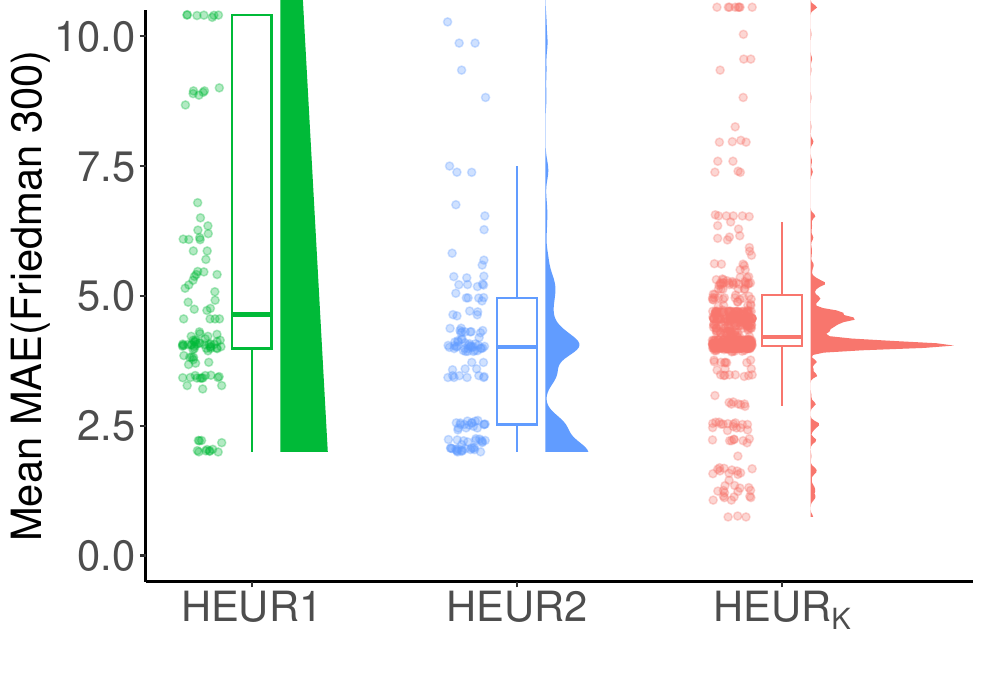}
    \end{minipage}%
    \begin{minipage}{0.55\textwidth}
        \centering
        \includegraphics[width=0.85\linewidth]{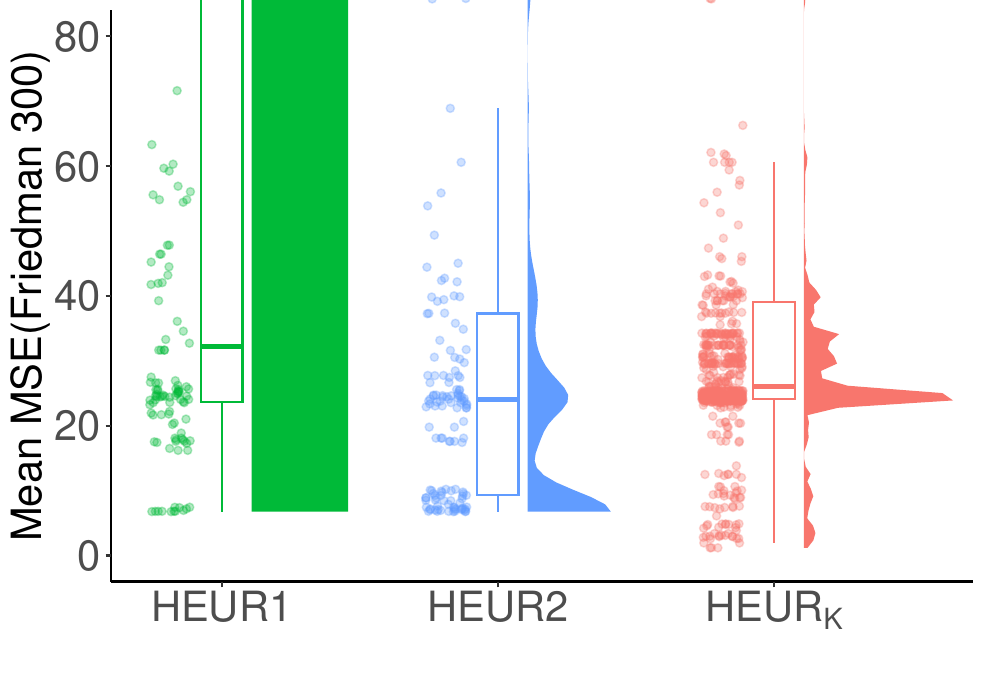}        
    \end{minipage}
    \label{fig:FRIEDMAN_300}
\end{figure}

\begin{figure}[!htb]
  \caption{ Average MAE and MSE  training values for Friedman(500) dataset}
  \hspace*{-3cm}
    \centering
    \begin{minipage}{.55\textwidth}
        \centering
        \includegraphics[width=0.85\linewidth]{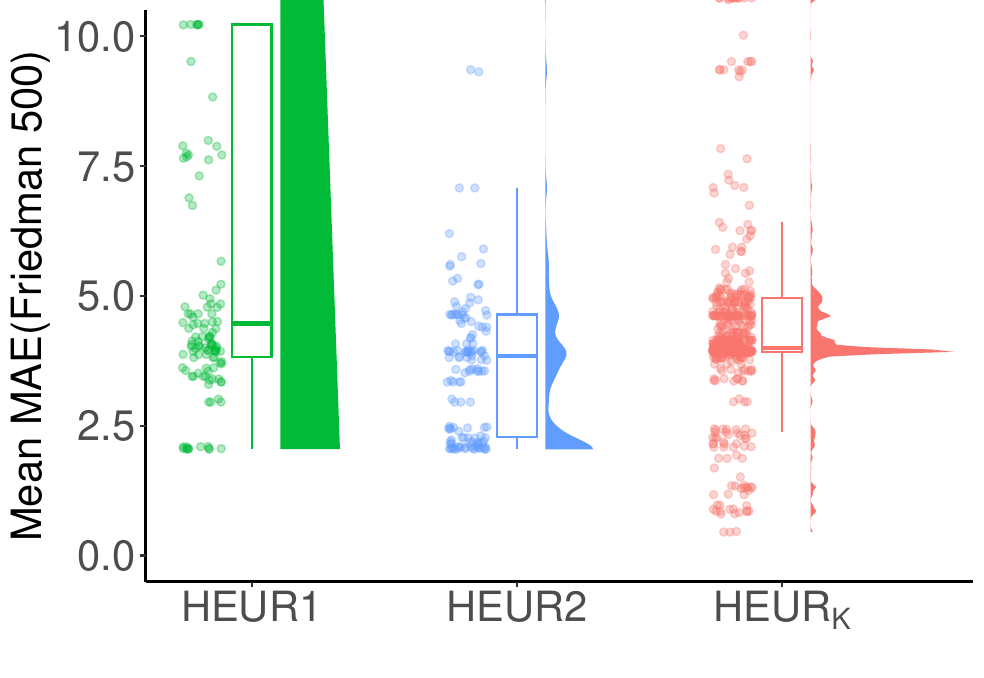}
    \end{minipage}%
    \begin{minipage}{0.55\textwidth}
        \centering
        \includegraphics[width=0.85\linewidth]{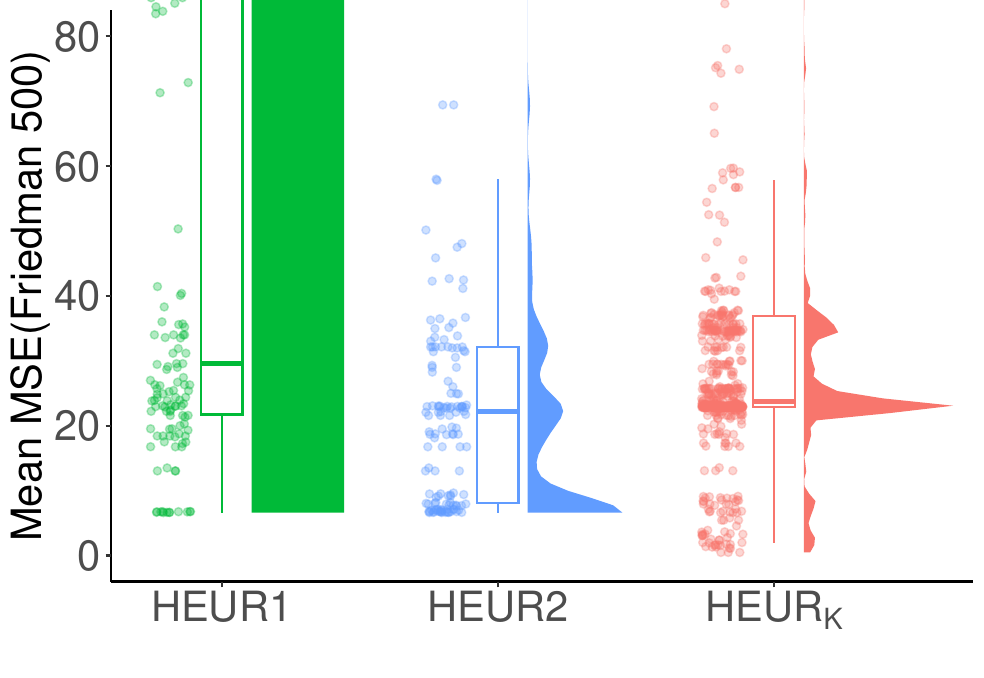}        
    \end{minipage}
    \label{fig:FRIEDMAN_500}
\end{figure}

Figures \ref{fig:FRIEDMAN_150} to \ref{fig:FRIEDMAN_500}
show the raincloud plots associated with the average MAE and MSE in the training sample. We observe that the best MAE and MSE results are derived from HEUR2 and HEUR$_K$. 

\begin{figure}[h!]
\hspace*{-3cm}
    \centering
    \begin{minipage}{.65\textwidth}
     \caption{Average computing time for akC weight}
        \centering
        \includegraphics[width=0.85\linewidth]{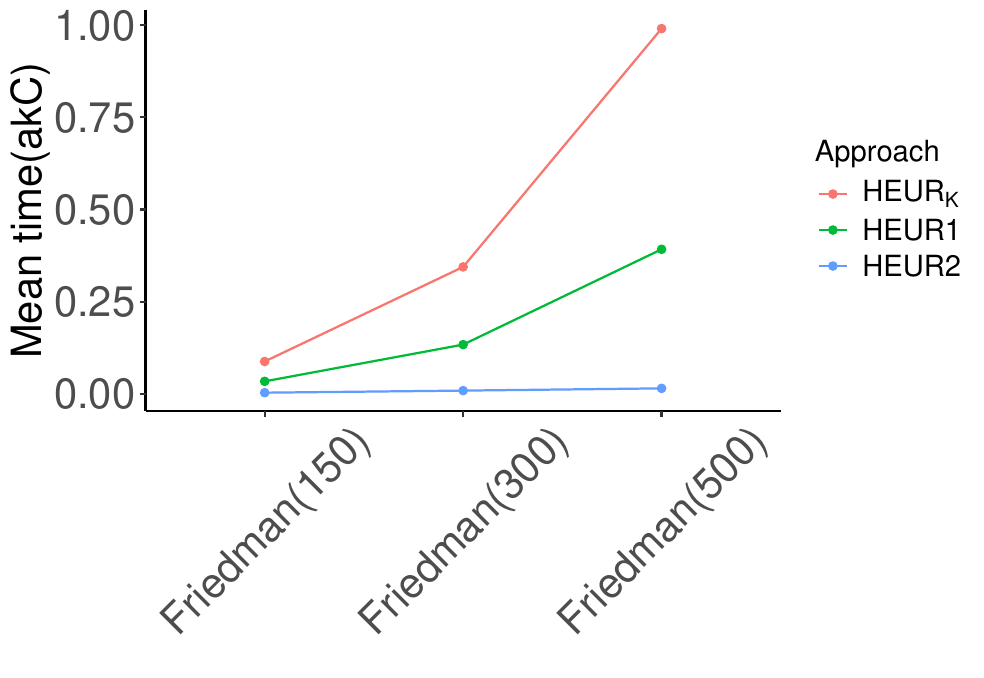}        
    \label{fig:FRIEDMAN_TIME1}
    \end{minipage}%
    \begin{minipage}{0.65\textwidth}
    \caption{Average computing time for MED weight}
        \centering
        \includegraphics[width=0.85\linewidth]{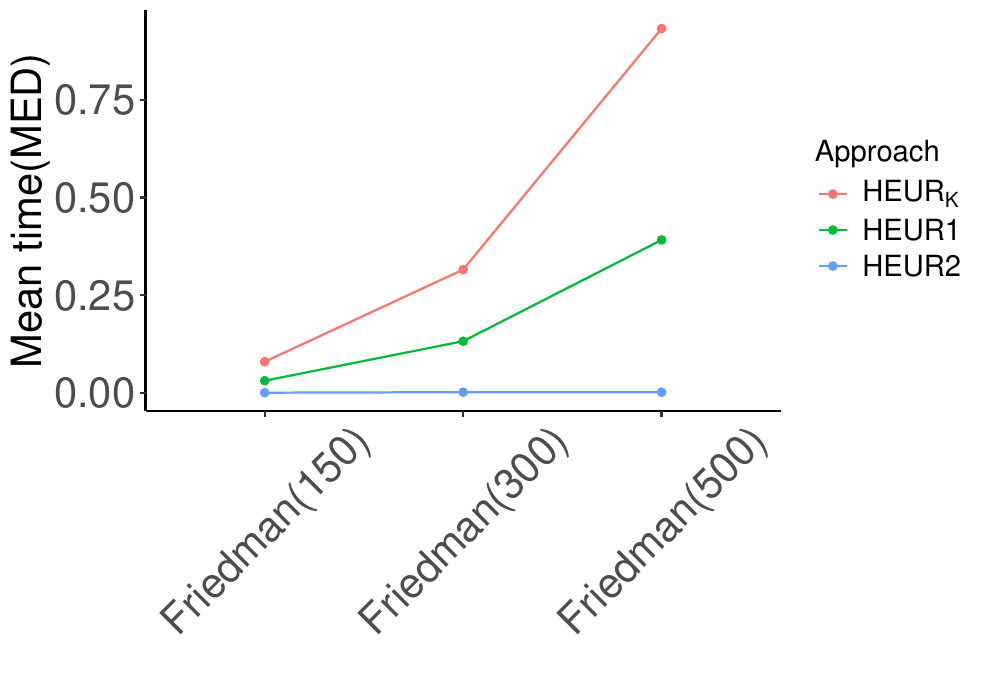}       
    \label{fig:FRIEDMAN_TIME2}
    \end{minipage}
\end{figure}

\begin{figure}[h!]
   \caption{Average computing time for TM weight}
    \centering
    \includegraphics[width=0.5\linewidth]{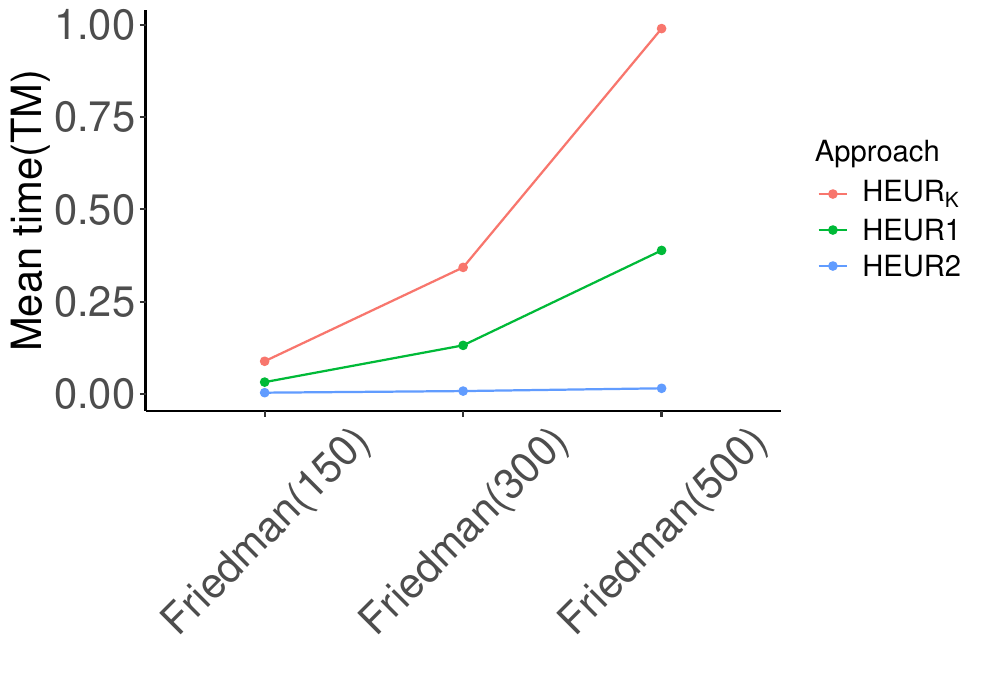}
    \label{fig:FRIEDMAN_TIME3}
\end{figure}

Since, in previous sections, we studied that non-convex operators ($\lambda$ weights) are the ones presenting bigger computing times, we restrict now to the study of akC, MED and TM times. Figures \ref{fig:FRIEDMAN_TIME1} to \ref{fig:FRIEDMAN_TIME3} report these times. We can see that, although the sample sizes are bigger, the computational time of the heuristics is still small.

\section{Conclusions}\label{sec:conclusions}

This paper introduces a support vector regression model where OWA operators can be included. This allows us to differently penalize each observation depending on the size associated with its deviation. Particularly, we provide a general model (NC-OWA-SVR) which consists in a quadratic mixed integer formulation. In addition, for the case of non-decreasing weights, a continuous quadratic formulation is developed (C-OWA-SVR).

Taking advantage of the dual formulations associated with these models, non-linear kernel functions can be accommodated allowing the use of non-linear regression functions. Besides, three heuristic approaches were proposed with the aim of providing good quality hyperplanes or regression functions in smaller times than the exact approaches.

In the computational results, we introduce a measure known as Goodness of  Fitting Index that is a generalization of an index previously used in regression; see \cite{BLANCO2018}. These results show that the exact approaches provide the best results of this index and they seem to be the most robust approaches for GoF. In general, we can observe that GoF is near 1 in many cases. 

Regarding the mean average error (MAE) and mean square error (MSE) values, we can also conclude that the exact formulations report the best results in most of the cases. In addition, the graphical insights show that HEUR2, OWA-SVR (exact formulation with linear kernel) and HEUR$_K$ are the approaches with best median results.

Although the exact approaches generally  provide the best results in our tests, they are partially limited by the required computing time to handle  large data sets.  In contrast,  we should mention that our heuristic approaches provide good alternatives in these cases. 

As future research, it would be interesting to explore enhancements to the proposed exact formulations to handle larger instances more efficiently.
Another interesting research line could be to study the weights selection depending on the specific application or dataset where this methodology is applied.
Besides,
the use of the $\ell_p$-norm in these formulations could provide new insights. In addition, the development of novel exact formulations for certain weights, based on their intrinsic properties, opens an interesting avenue for further study.

\section*{Acknowledgements}

The authors acknowledge the financial support of
the Spanish Ministry of Science and Innovation and ERDF via projects PID2020-114594GB-C21/C22. First and third author also thank the support of the European Union “NextGenerationEU”/PRTR through project TED2021-130875B-I00.


\appendix
\section{Algorithms flowcharts} \label{App:alg}
We present the flowcharts related to the heuristic approaches of Section \ref{sec:heur}.

\begin{figure}[!htb]
    \caption{Flowchart of Algorithm 1}
    \centering
    \includegraphics[width=0.95\linewidth]{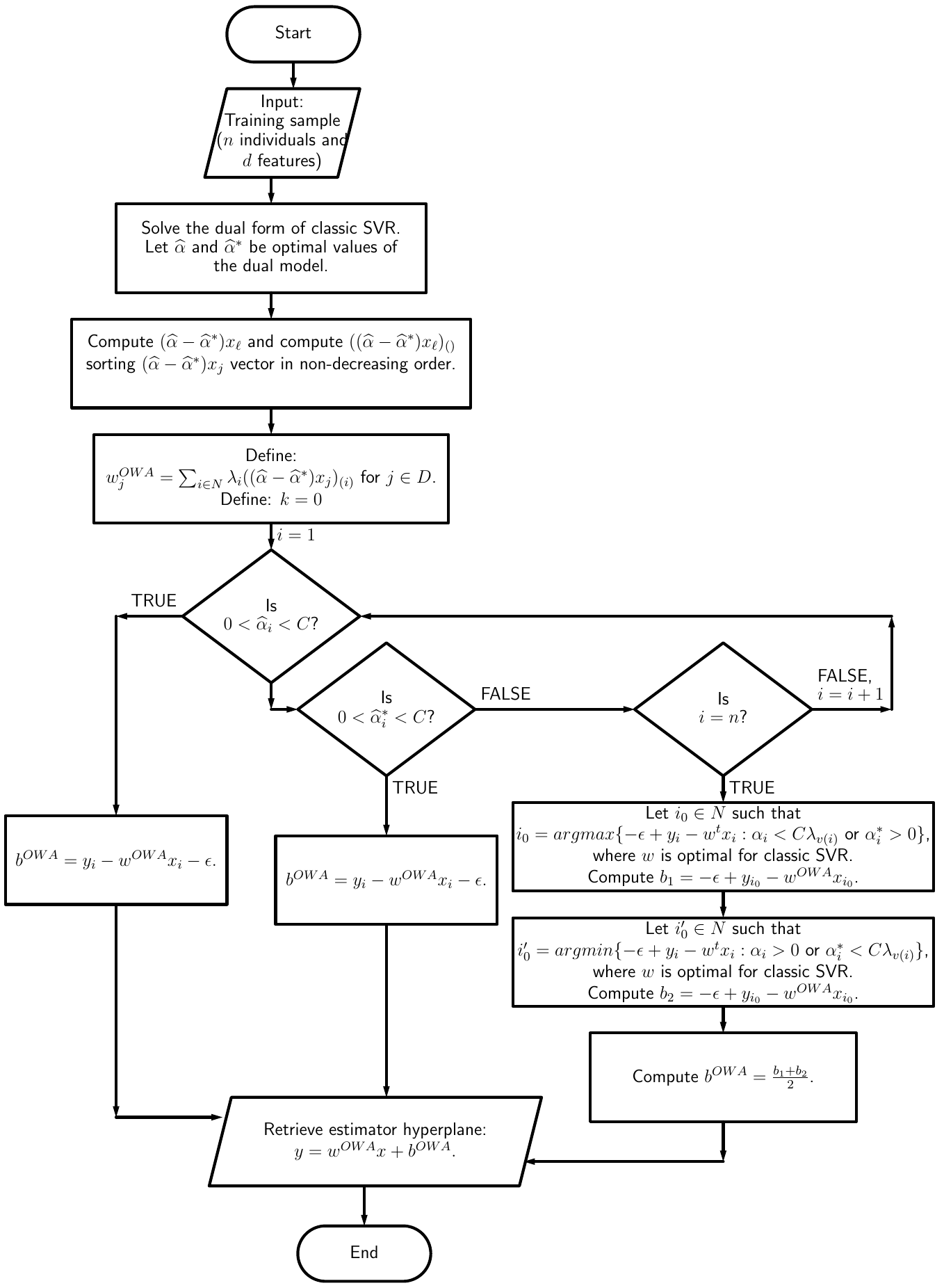}
    \label{fig:alg1}
\end{figure}

\begin{figure}[!htb]
\hspace*{-3cm}
    \centering
    \begin{subfigure}{.55\textwidth}
    \caption{Flowchart of Algorithm 2}
        \centering
        \includegraphics[width=0.85\linewidth]{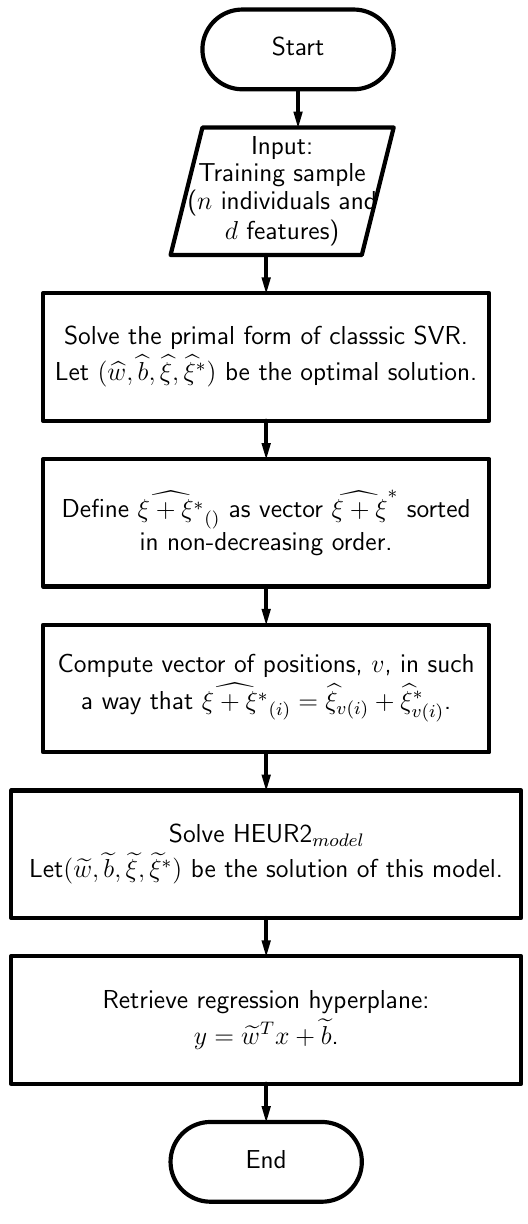}
            \label{fig:alg2}
    \end{subfigure}   
    \begin{subfigure}{0.55\textwidth}
     \caption{Flowchart of Algorithm 3}
        \centering
        \includegraphics[width=0.85\linewidth]{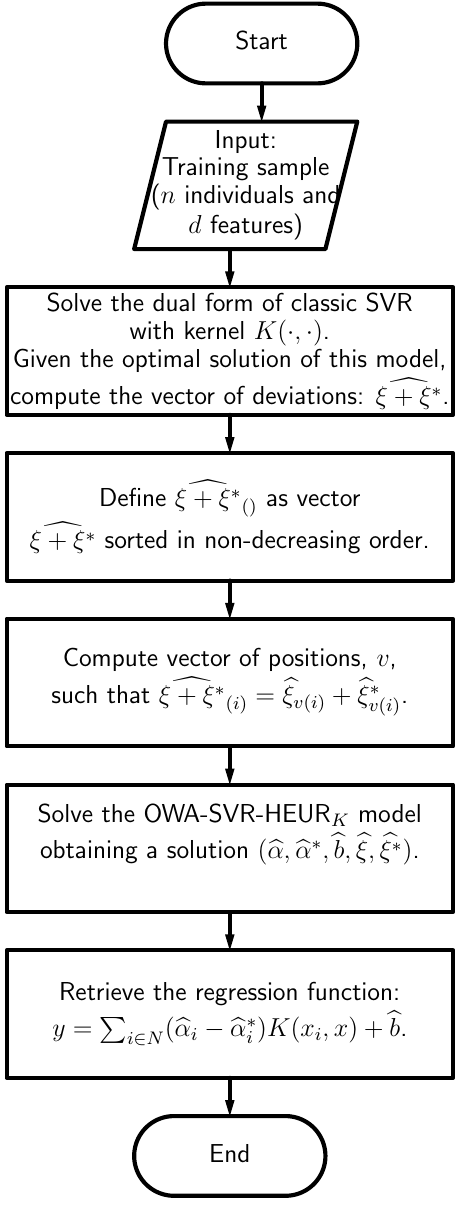}
            \label{fig:alg3}
        \end{subfigure}
\end{figure}
\clearpage
\section{ Statistical tests for evaluation metrics}
\label{app:test}

Table \ref{tab:Kruskal} reports the results of Kruskal-Wallis test for average MAE and average MSE obtained in the training of each dataset. This test allows us to decide whether there is a significant difference between the mean MAE values of the different approaches (HEUR1, HEUR2, OWA-SVR, HEUR$_K$, OWA-SVR$_K$) or not. The same test is carried out using the mean MSE results for each dataset.  Observe in Table \ref{tab:Kruskal} that all the $p$-values are really small. Consequently, there exists significant difference in all cases.

\begin{table}[htbp]
  \centering
  \caption{Kruskal-Wallis test for average MAE and MSE in different datasets.}
    \begin{tabular}{lrrrr}
    \hline
    \hline
          \multicolumn{5}{c}{Kruskal-Wallis test } \\
    \midrule
    \multicolumn{1}{c}{\multirow{2}[2]{*}{Data}} & \multicolumn{2}{c}{Average MAE} & \multicolumn{2}{c}{Average MSE} \\
          & \multicolumn{1}{l}{Statistic} & \multicolumn{1}{l}{$p$-value } & \multicolumn{1}{l}{Statistic} & \multicolumn{1}{l}{$p$-value } \\
    \hline
    \hline
    Quandt & 369.52 & \multicolumn{1}{l}{$<$ 2.2e-16} & 372.72 & \multicolumn{1}{l}{$<$2.2e-16} \\
    Boston & 109.45 & \multicolumn{1}{l}{$<$ 2.2e-16} & 154.79 & \multicolumn{1}{l}{$<$ 2.2e-16} \\
    CYGOB1 & 144.76 & \multicolumn{1}{l}{$<$ 2.2e-16} & 149.38 & \multicolumn{1}{l}{$<$ 2.2e-16} \\
    Eilon & 204.35 & \multicolumn{1}{l}{$<$ 2.2e-16} & 220.3 & \multicolumn{1}{l}{$<$ 2.2e-16} \\
    DW2   & 94.99 & \multicolumn{1}{l}{$<$ 2.2e-16} & 139.01 & \multicolumn{1}{l}{$<$ 2.2e-16} \\
    Alcohol & 418.1 & \multicolumn{1}{l}{$<$ 2.2e-16} & 450.16 & \multicolumn{1}{l}{$<$ 2.2e-16} \\
    DW3   & 128.06 & \multicolumn{1}{l}{$<$ 2.2e-16} & 213.38 & \multicolumn{1}{l}{$<$ 2.2e-16} \\
    HBK   & 148.96 & \multicolumn{1}{l}{$<$ 2.2e-16} & 231.55 & \multicolumn{1}{l}{$<$ 2.2e-16} \\
    Friedman (150) & 38.44 & 4.50E-09 & 387.77 & \multicolumn{1}{l}{$<$ 2.2e-16} \\
    Friedman (300) & 37.07 & 8.92E-09 & 391.12 & \multicolumn{1}{l}{$<$ 2.2e-16} \\
    Friedman (500) & 33    & 6.82E-08 & 392.77 & $<$2.20E-16 \\
    \hline
    \hline
    \end{tabular}%
  \label{tab:Kruskal}%
\end{table}%

Next tables report the results of pairwise Wilcoxon test for each evaluation metric and dataset. We can observe that, in most of the cases, there exists a significant difference in the mean MAE and MSE values between each pair of approaches.


\begin{table}[htbp]
  \centering
  \caption{Test on average MAE and MSE results for Quandt dataset}
  \scalebox{0.8}{
    \begin{tabular}{|l|rrrr|rrrr|}
    \hline
    \multicolumn{1}{|c|}{\multirow{2}[4]{*}{}} & \multicolumn{4}{c|}{Wilcoxon test on average MAE results} & \multicolumn{4}{c|}{Wilcoxon test on average MSE results} \\
          & \multicolumn{1}{l}{HEUR$_K$} & \multicolumn{1}{l}{HEUR1} & \multicolumn{1}{l}{HEUR2} & \multicolumn{1}{l|}{OWA-SVR$_K$} & \multicolumn{1}{l}{HEUR$_K$} & \multicolumn{1}{l}{HEUR1} & \multicolumn{1}{l}{HEUR2} & \multicolumn{1}{l|}{OWA-SVR$_K$} \\
     \hline
    HEUR1 & $<$2e-16 & -     & -     & -     & $<$2e-16 & -     & -     & - \\
    HEUR2 & 3.00E-12 & $<$2e-16 & -     & -     & 1.20E-12 & $<$2e-16 & -     & - \\
    OWA-SVR$_K$ & 0.0041 & $<$2e-16 & 5.60E-10 & -     & 0.00032 & $<$2e-16 & 3.40E-10 & - \\
    OWA-SVR & $<$2e-16 & $<$2e-16 & 6.60E-08 & $<$2e-16 & $<$2e-16 & $<$2e-16 & 4.60E-08 & $<$2e-16 \\
     \hline
    \end{tabular}
    }
  \label{tab:TESQUANDT}%
\end{table}%

\begin{table}[htbp]
  \centering
  \caption{Test on average MAE and MSE results for Boston dataset}
   \scalebox{0.8}{
    \begin{tabular}{|l|rrrr|rrrl|}
    \hline
    \multirow{2}[2]{*}{} & \multicolumn{4}{c|}{Wilcoxon test on average MAE results } & \multicolumn{4}{c|}{Wilcoxon test on average MSE results } \\
          & \multicolumn{1}{l}{HEUR$_K$} & \multicolumn{1}{l}{HEUR1} & \multicolumn{1}{l}{HEUR2} & \multicolumn{1}{l|}{OWA-SVR$_K$} & \multicolumn{1}{l}{HEUR$_K$} & \multicolumn{1}{l}{HEUR1} & \multicolumn{1}{l}{HEUR2} & OWA-SVR$_K$ \\
    \hline
    HEUR1 & 2.50E-13 & \multicolumn{1}{l}{-} & \multicolumn{1}{l}{-} & \multicolumn{1}{l|}{-} & 5.00E-15 & \multicolumn{1}{l}{-} & \multicolumn{1}{l}{-} & - \\
    HEUR2 & 0.0109 & 4.10E-11 & \multicolumn{1}{l}{-} & \multicolumn{1}{l|}{-} & 0.021 & 6.50E-12 & \multicolumn{1}{l}{-} & - \\
    K-OWA-SVR & 0.00096 & 1.10E-09 & 6.00E-05 & \multicolumn{1}{l|}{-} & 1.50E-10 & 8.00E-09 & 1.20E-07 & - \\
    OWA-SVR & 9.10E-06 & 4.80E-14 & 0.24356 & 7.70E-11 & 3.00E-06 & 1.90E-15 & 0.184 & $<$2e-16 \\
    \hline
    \end{tabular}%
    }
  \label{tab:TESTBOSTON}%
\end{table}%

\begin{table}[htbp]
  \centering
  \caption{Test on average MAE and MSE results for CYGOB1 dataset}
  \scalebox{0.8}{
    \begin{tabular}{|l|rrrr|rrrr|}
    \hline
    \multirow{2}[2]{*}{} & \multicolumn{4}{c|}{Wilcoxon test on average MAE results } & \multicolumn{4}{c|}{Wilcoxon test on average MSE results } \\
          & \multicolumn{1}{l}{HEUR$_K$} & \multicolumn{1}{l}{HEUR1} & \multicolumn{1}{l}{HEUR2} & \multicolumn{1}{l|}{OWA-SVR$_K$} & \multicolumn{1}{l}{HEUR$_K$} & \multicolumn{1}{l}{HEUR1} & \multicolumn{1}{l}{HEUR2} & \multicolumn{1}{l|}{OWA-SVR$_K$} \\
    \hline
    HEUR1 & $<$2e-16 & -     & -     & -     & $<$2e-16 & -     & -     & - \\
    HEUR2 & 0.00025 & 4.70E-12 & -     & -     & 0.00098 & 1.30E-12 & -     & - \\
    OWA-SVR$_K$ & 0.01169 & $<$2e-16 & 0.02577 & -     & 0.00484 & $<$2e-16 & 0.04716 & - \\
    OWA-SVR & 0.7847 & $<$2e-16 & 1.60E-05 & 0.00313 & 0.29476 & $<$2e-16 & 8.10E-06 & 0.00127 \\
    \hline
    \end{tabular}%
    }
  \label{tab:TESTCYGOB1}%
\end{table}%

\begin{table}[htbp]
  \centering
   \caption{Test on average MAE and MSE results for Eilon dataset}
  \scalebox{0.8}{
    \begin{tabular}{|l|rrrr|rrrr|}
    \hline
    \multirow{2}[2]{*}{} & \multicolumn{4}{c|}{Wilcoxon test on average MAE results } & \multicolumn{4}{c|}{Wilcoxon test on average MSE results } \\
          & \multicolumn{1}{l}{HEUR$_K$} & \multicolumn{1}{l}{HEUR1} & \multicolumn{1}{l}{HEUR2} & \multicolumn{1}{l|}{OWA-SVR$_K$} & \multicolumn{1}{l}{HEUR$_K$} & \multicolumn{1}{l}{HEUR1} & \multicolumn{1}{l}{HEUR2} & \multicolumn{1}{l|}{OWA-SVR$_K$} \\
    \hline
    HEUR1 & $<$2e-16 & -     & -     & -     & 3.00E-15 & -     & -     & - \\
    HEUR2 & 0.6176 & 2.10E-09 & -     & -     & 0.75  & 6.20E-09 & -     & - \\
    OWA-SVR$_K$ & $<$2e-16 & 8.20E-07 & 2.10E-07 & -     & $<$2e-16 & 1.70E-06 & 2.30E-08 & - \\
    OWA-SVR & 1.00E-05 & $<$2e-16 & 0.0001 & $<$2e-16 & 7.60E-07 & $<$2e-16 & 2.30E-05 & $<$2e-16 \\
    \hline
    \end{tabular}%
    }
  \label{tab:TESTEILON}%
\end{table}%

\begin{table}[htbp]
  \centering
     \caption{Test on average MAE and MSE results for DW2 dataset}
  \scalebox{0.8}{
    \begin{tabular}{|l|rrrr|rrrr|}
    \hline
    \multirow{2}[2]{*}{} & \multicolumn{4}{c|}{Wilcoxon test on average MAE results } & \multicolumn{4}{c|}{Wilcoxon test on average MSE results } \\
          & \multicolumn{1}{l}{HEUR$_K$} & \multicolumn{1}{l}{HEUR1} & \multicolumn{1}{l}{HEUR2} & \multicolumn{1}{l|}{OWA-SVR$_K$} & \multicolumn{1}{l}{HEUR$_K$} & \multicolumn{1}{l}{HEUR1} & \multicolumn{1}{l}{HEUR2} & \multicolumn{1}{l|}{OWA-SVR$_K$} \\
    \hline
    HEUR1 & $<$2e-16 & -     & -     & -     & \multicolumn{1}{l}{$<$2e-16} & \multicolumn{1}{l}{-} & \multicolumn{1}{l}{-} & \multicolumn{1}{l|}{-} \\
    HEUR2 & 0.6791 & 1.50E-12 & -     & -     & 0.151 & 4.70E-13 & \multicolumn{1}{l}{-} & \multicolumn{1}{l|}{-} \\
    OWA-SVR$_K$ & 0.8129 & 6.00E-16 & 0.8129 & -     & 4.80E-08 & 1.80E-10 & 1.40E-06 & \multicolumn{1}{l|}{-} \\
    OWA-SVR & 0.0028 & 3.10E-15 & 0.1209 & 0.016 & 0.001 & 4.00E-14 & 0.258 & 3.60E-13 \\
    \hline
    \end{tabular}%
    }
  \label{tab:TESTDW2}%
\end{table}%

\begin{table}[htbp]
  \centering
   \caption{Test on average MAE and MSE results for Alcohol dataset}
  \scalebox{0.8}{
    \begin{tabular}{|l|rrrr|rrrr|}
     \hline
    \multirow{2}[2]{*}{} & \multicolumn{4}{c|}{Wilcoxon test on average MAE results } & \multicolumn{4}{c|}{Wilcoxon test on average MSE results } \\
          & \multicolumn{1}{l}{HEUR$_K$} & \multicolumn{1}{l}{HEUR1} & \multicolumn{1}{l}{HEUR2} & \multicolumn{1}{l|}{OWA-SVR$_K$} & \multicolumn{1}{l}{HEUR$_K$} & \multicolumn{1}{l}{HEUR1} & \multicolumn{1}{l}{HEUR2} & \multicolumn{1}{l|}{OWA-SVR$_K$} \\
     \hline
    HEUR1 & $<$2e-16 & -     & -     & -     & $<$2e-16 & -     & -     & - \\
    HEUR2 & 1.80E-12 & $<$2e-16 & -     & -     & 1.20E-11 & $<$2e-16 & -     & - \\
    OWA-SVR$_K$ & 0.19  & $<$2e-16 & 5.60E-11 & -     & 1.40E-05 & $<$2e-16 & 2.10E-14 & - \\
    OWA-SVR & $<$2e-16 & $<$2e-16 & 1.70E-05 & $<$2e-16 & $<$2e-16 & $<$2e-16 & 5.80E-06 & $<$2e-16 \\
    \hline
    \end{tabular}%
    }
  \label{tab:TESTALCOHOL}%
\end{table}%

\begin{table}[htbp]
  \centering
  \caption{Test on average MAE and MSE results for DW3 dataset}
  \scalebox{0.8}{
    \begin{tabular}{|l|rrrr|rrrr|}
    \hline
    \multirow{2}[2]{*}{} & \multicolumn{4}{c|}{Wilcoxon test on average MAE results } & \multicolumn{4}{c|}{Wilcoxon test on average MSE results } \\
          & \multicolumn{1}{l}{HEUR$_K$} & \multicolumn{1}{l}{HEUR1} & \multicolumn{1}{l}{HEUR2} & \multicolumn{1}{l|}{OWA-SVR$_K$} & \multicolumn{1}{l}{HEUR$_K$} & \multicolumn{1}{l}{HEUR1} & \multicolumn{1}{l}{HEUR2} & \multicolumn{1}{l|}{OWA-SVR$_K$} \\
    \hline
    HEUR1 & $<$2e-16 & -     & -     & -     & $<$2e-16 & -     & -     & - \\
    HEUR2 & 0.19881 & $<$2e-16 & -     & -     & 0.044 & $<$2e-16 & -     & - \\
    OWA-SVR$_K$ & 0.59831 & $<$2e-16 & 0.0887 & -     & 9.10E-13 & $<$2e-16 & 8.10E-09 & - \\
    OWA-SVR & 0.00012 & $<$2e-16 & 0.0887 & 9.70E-05 & 8.30E-06 & $<$2e-16 & 0.095 & $<$2e-16 \\
    \hline
    \end{tabular}%
    }
  \label{tab:TESTDW3}%
\end{table}%

\begin{table}[htbp]
  \centering
   \caption{Test on average MAE and MSE results for HBK dataset}
  \scalebox{0.8}{
    \begin{tabular}{|l|rrrr|rrrr|}
    \hline
    \multirow{2}[2]{*}{} & \multicolumn{4}{c|}{Wilcoxon test on average MAE results } & \multicolumn{4}{c|}{Wilcoxon test on average MSE results } \\
          & \multicolumn{1}{l}{HEUR$_K$} & \multicolumn{1}{l}{HEUR1} & \multicolumn{1}{l}{HEUR2} & \multicolumn{1}{l|}{OWA-SVR$_K$} & \multicolumn{1}{l}{HEUR$_K$} & \multicolumn{1}{l}{HEUR1} & \multicolumn{1}{l}{HEUR2} & \multicolumn{1}{l|}{OWA-SVR$_K$} \\
    \hline
    HEUR1 & $<$2e-16 & -     & -     & -     & $<$2e-16 & -     & -     & - \\
    HEUR2 & 0.0408 & $<$2e-16 & -     & -     & 2.80E-10 & $<$2e-16 & -     & - \\
    OWA-SVR$_K$ & 0.0046 & $<$2e-16 & 0.9141 & -     & 4.20E-06 &$<$2e-16 & 0.0018 & - \\
    OWA-SVR & 0.2477 & $<$2e-16 & 0.6946 & 0.0408 & 5.20E-06 & $<$2e-16 & 0.0347 & 0.3458 \\
    \hline
    \end{tabular}%
    }
  \label{tab:TESTHBK}%
\end{table}%

\begin{table}[htbp]
  \centering
  \caption{Test on average MAE and MSE results for Friedman (150) dataset}
  \scalebox{0.8}{
    \begin{tabular}{|l|rr|rr|}
    \hline
    \multirow{2}[2]{*}{} & \multicolumn{2}{c|}{Wilcoxon test on average MAE results } & \multicolumn{2}{c|}{Wilcoxon test on average MSE results } \\
          & \multicolumn{1}{c}{HEUR$_K$} & \multicolumn{1}{c|}{HEUR1} & \multicolumn{1}{c}{HEUR$_K$} & \multicolumn{1}{c|}{HEUR1} \\
    \hline
    HEUR1 & 0.088 & -     & 0.072 & - \\
    HEUR2 & 4.10E-10 & 0.103 & $<$2e-16 & $<$2e-16 \\
    \hline
    \end{tabular}%
    }
  \label{tab:TESTF150}%
\end{table}%

\begin{table}[htbp]
  \centering
   \caption{Test on average MAE and MSE results for Friedman (300) dataset}
  \scalebox{0.8}{
    \begin{tabular}{|l|rr|rr|}
    \hline
    \multirow{2}[2]{*}{} & \multicolumn{2}{c|}{Wilcoxon test on average MAE results } & \multicolumn{2}{c|}{Wilcoxon test on average MSE results } \\
          & \multicolumn{1}{c}{HEUR$_K$} & \multicolumn{1}{c|}{HEUR1} & \multicolumn{1}{c}{HEUR$_K$} & \multicolumn{1}{c|}{HEUR1} \\
\hline    HEUR1 & 0.022 & -     & 0.0052 & - \\
    HEUR2 & 1.60E-09 & 0.592 & $<$2e-16 & $<$2e-16 \\
    \hline
    \end{tabular}%
    }
  \label{tab:TESTF300}%
\end{table}%

\begin{table}[htbp]
  \centering
   \caption{Test on average MAE and MSE results for Friedman (500) dataset}
  \scalebox{0.8}{
    \begin{tabular}{|l|rr|rr|}
    \hline
    \multirow{2}[2]{*}{} & \multicolumn{2}{c|}{Wilcoxon test on average MAE results } & \multicolumn{2}{c|}{Wilcoxon test on average MSE results } \\
          & \multicolumn{1}{c}{HEUR$_K$} & \multicolumn{1}{c|}{HEUR1} & \multicolumn{1}{c}{HEUR$_K$} & \multicolumn{1}{c|}{HEUR1} \\
    \hline
    HEUR1 & 0.022 & -     & 0.0052 & - \\
    HEUR2 & 1.60E-09 & 0.592 & $<$2e-16 & $<$2e-16 \\
   \hline
    \end{tabular}%
    }
  \label{tab:TESTF500}%
\end{table}%

\end{document}